\documentclass[reqno,11pt]{amsart}
\UseRawInputEncoding
\pdfoutput=1
\usepackage{amsfonts}
\usepackage{graphicx}
\usepackage{amsfonts,amsmath, amssymb,mathrsfs}
\usepackage{hyperref}
\usepackage{marginnote}
\usepackage{color}
\usepackage{tikz}
\usepackage{verbatim}
\allowdisplaybreaks

\numberwithin{equation}{section}

\newcommand{\dis}{\displaystyle}

\newcommand{\R}{\mathbb{R}}

\newtheorem{theorem}{Theorem}[section]
\newtheorem{corollary}[theorem]{Corollary}
\newtheorem{lemma}[theorem]{Lemma}

\newtheorem{remark}[theorem]{Remark}

\def\v{\varepsilon}

\def\x{\xi}

\def\d{\delta}

\def\dis{\displaystyle}

\def\f{\frac}

\begin{document}
	
\title[Uniqueness of composite wave of shock and rarefaction]{Uniqueness of composite wave of shock and rarefaction in the inviscid limit of Navier-Stokes equations}

\author[F.M. Huang]{Feimin Huang}
\address[Feimin Huang]{Institute of Applied Mathematics, Academy of Mathematics and Systems Science, Chinese Academy of Sciences, Beijing 100190, China}
\email{fhuang@amt.ac.cn}

\author[W.Q. Wang]{Weiqiang Wang}
\address[W.Q. Wang]{Institute of Applied Mathematics, Academy of Mathematics and Systems Science, Chinese Academy of Sciences, Beijing 100190, China}
\email{wangweiqiang@amss.ac.cn}

\author[Y. Wang]{Yi Wang}
\address[Yi Wang]{Institute of Applied Mathematics, Academy of Mathematics and Systems Science, Chinese Academy of Sciences, Beijing 100190, China}
\email{wangyi@amss.ac.cn}

\author[Y. Wang]{Yong Wang}
\address[Yong Wang]{Institute of Applied Mathematics, Academy of Mathematics and Systems Science, Chinese Academy of Sciences, Beijing 100190, China}
\email{yongwang@amss.ac.cn}

\begin{abstract}
	The uniqueness of entropy solution for the compressible Euler equations is a fundamental and challenging problem. In this paper, the uniqueness of a composite wave of shock and rarefaction of 1-d compressible Euler equations is proved in the inviscid limit of compressible Navier-Stokes equations.
	Moreover, the relative entropy around the original Riemann solution consisting of shock and rarefaction under the large perturbation is shown to be uniformly bounded by the framework developed in \cite{Kang-Vasseur-2021-Invent}. The proof contains two new ingredients: 1) a cut-off technique and the expanding property of rarefaction are used to overcome the errors generated by the viscosity related to inviscid rarefaction; 2)
	the error terms concerning the interactions between shock and rarefaction are controlled by the compressibility of shock, the decay of derivative of rarefaction and the separation of shock and rarefaction as time increases.
\end{abstract}

\date{\today}
\maketitle

\setcounter{tocdepth}{1}
\tableofcontents

\thispagestyle{empty}


\section{Introduction and main result}

Consider the one-dimensional compressible isentropic Navier-Stokes equations in Lagrangian coordinate
\begin{align}\label{1.1}
	\begin{cases}
		\dis \partial_t \mathfrak{v}^{\nu}-\partial_x \mathfrak{u}^{\nu}=0,\\[1mm]
		\dis \partial_t \mathfrak{u}^{\nu}+\partial_xp(\mathfrak{v}^{\nu}) =\nu\partial_x(\frac{b(\mathfrak{v}^{\nu})}{\mathfrak{v}^{\nu}} \partial_x\mathfrak{u}^{\nu}),
	\end{cases}\qquad (t,x)\in \R_{+}^2:=(0,+\infty)\times \R,
\end{align}
where $\mathfrak{v}^{\nu}>0$ is the specific volume, $\mathfrak{u}^{\nu}$ is the velocity of fluid, and $\nu>0$ is the inverse of Reynolds number. In the present paper, we consider the polytropic gas, that is, the pressure is given by
$$p(\mathfrak{v})=\mathfrak{v}^{-\gamma}\qquad \text{with }\gamma>1.$$
We assume the viscosity coefficient $b(\mathfrak{v})$ depends on the specific volume, and takes the form
$$b(\mathfrak{v})=b \mathfrak{v}^{-\alpha}\qquad  \text{with}\,\,  b,\,\,\alpha>0.$$
Without loss of generality, we assume $b=\gamma$ throughout this paper.

Formally, as $\nu\to0+$, the system \eqref{1.1} converges to the following compressible Euler system in Lagrangian coordinate
\begin{align}\label{1.2}
	\begin{cases}
		\partial_t \mathfrak{v}-\partial_x \mathfrak{u}=0,\\[1mm]
		\partial_t \mathfrak{u}+\partial_xp(\mathfrak{v}) =0,
	\end{cases}\qquad (t,x)\in \R_{+}^2:=[0,+\infty)\times \R,
\end{align}
For \eqref{1.2}, we impose the following Riemann initial data
\begin{align}\label{1.3}
	(\mathfrak{v},\mathfrak{u})|_{t=0}=
	\begin{cases}
		(v_-,u_-),\quad x<0,\\
		(v_+,u_+), \quad x>0.
	\end{cases}
\end{align}
It is well-known that the Euler system \eqref{1.2} is hyperbolic. This means that the matrix
\begin{equation}\nonumber
	\left(
	\begin{array}{cc}
		0 & \,\, -1  \\
		p'(\mathfrak{v}) & \,\, 0 \\
	\end{array}
	\right),
\end{equation}
has two real eigenvalues:
\begin{align}
	\lambda_1=-\sqrt{-p'(\mathfrak{v})},\quad\mbox{and}\quad \lambda_2=\sqrt{-p'(\mathfrak{v})}.\nonumber
\end{align}

For later use, we denote $U_-:=(v_-,u_-)$ and $U_+:=(v_+,u_+)$. The Riemann solution is determined by a combination of at most two elementary solutions from the following four families: 1-shock, 2-shock, 1-rarefaction and 2-rarefaction, see \cite{Dafermos-2016, Smoller-1994} for more details. We define integral curves (see Figure 1) across the point $U_-$
\begin{align}
	&RC_1:=\Big\{(\mathfrak{v},\mathfrak{u})\, |\, \mathfrak{v}\geq v_-,\,\, \mathfrak{u}=u_--\int_{v_-}^{\mathfrak{v}}\lambda_1(s)ds\Big\},\nonumber\\
	&RC_2:=\Big\{(\mathfrak{v},\mathfrak{u})\, |\, 0<\mathfrak{v}\leq v_-,\,\, \mathfrak{u}=u_--\int_{v_-}^{\mathfrak{v}} \lambda_2(s)ds\Big\},\nonumber\\
	&SC_1:=\left\{ (\mathfrak{v},\mathfrak{u})\, \bigg|\,  0<\mathfrak{v}\leq v_-,\,\, \exists \,\, \sigma \, s.t.\,\,\,
	\begin{array}{c}
		-\sigma (\mathfrak{v}-v_-)=(\mathfrak{u}-u_-),\\
		-\sigma (\mathfrak{u}-u_-)=(p(v_-)-p(\mathfrak{v}))
	\end{array}\right\},\nonumber\\
	&SC_2:=\left\{ (\mathfrak{v},\mathfrak{u})\, \bigg|\,\, \mathfrak{v}\geq v_-,\,\, \exists \,\, \sigma \, s.t.\,\,\,
	\begin{array}{c}
		-\sigma (\mathfrak{v}-v_-)=(\mathfrak{u}-u_-),\\
		-\sigma (\mathfrak{u}-u_-)=(p(\mathfrak{v})-p(v_-))
	\end{array}\right\}.\nonumber
\end{align}

The idea of constructing the solution of inviscid gases from viscous gases by vanishing viscosity limit can be dated back to Stokes \cite{Stokes}, and many works have been made on this direction. Based on $L^{\infty}$ compensated compactness method, the existence of bounded entropy solutions to Euler equations in Eulerian coordinate are proved by vanishing artificial viscosity, see \cite{Ding-Chen-Luo, Ding-Chen-Luo-1989, DiPerna-1983,  Lions-Perthame-Souganidis-1996, Lions-Perthame-Tadmor-1994}. Due to the lack of invariant region for compressible Naiver-Stokes equations with physical viscosity, the $L^{\infty}$ compactness framework fails. Instead, Chen and Perepelitsa \cite{Chen-Perepelista-2010} developed the $L^p$ compensated compactness framework introduced firstly in \cite{LeFloch-Westdickenberg-2007}, and proved the convergence of the solutions to Naiver-Stokes equations in Eulerian coordinate with constant viscosity (i.e. $\alpha=0$) towards a finite-energy entropy solution of Euler equations. We also refer to \cite{He-Wang-2022} for the case of viscosity depending on density. However, the uniqueness of these limit solutions is still open.

In the BV (bounded variation) setting, the celebrated progress was made by Glimm \cite{Glimm-1965}, in which he constructed the global weak entropy solutions with small BV norm for general strictly hyperbolic system. Later, the uniqueness of these weak entropy solutions was proved by Bressan, Crasta and Piccoli \cite{Bressan-Crasta-Piccoli-2000}, see also \cite{Bressan-2000,Bressan-Liu-Yang-1999,Chen-Krupa-Vasseur-2022,Liu-Yang-1999}. In 2005,
Bianchini and  Bressan \cite{Bianchini-Bressan-2005} showed the $BV$ solutions of artificial viscosity system converge to the unique BV solution of inviscid system as the viscosity goes to zero. Also in \cite{Bianchini-Bressan-2005}, they raised the problem of the inviscid limit of Navier-Stokes as an open problem. Recently, a remarkable progress towards this problem has been made by Kang and Vasseur \cite{Kang-Vasseur-2021-Invent}, in which they proved the single entropy shock solution of Euler equations \eqref{1.2} is stable and unique in the inviscid limit from a large family of Navier-Stokes system \eqref{1.1}. A bit later, they extend the result to the case of two shocks in \cite{Kang-Vasseur-2020}.
For the inviscid limit and uniqueness of contact discontinuity to Navier-Stokes-Fourier system, one can refer to \cite{Kang-Vasseur-Wang-CMP,Vasseur-Wang} and references therein. 

In present paper, we are interested in the stability and uniqueness of a composite wave of shock and rarefaction of compressible Euler equations \eqref{1.2} in the inviscid limit of Navier-Stokes equations \eqref{1.1}.
For clarity, we assume that the composite wave to \eqref{1.2}--\eqref{1.3} is a superposition of 1-shock wave and 2-rarefaction wave (see Figure 2), i.e., there exists a intermediate state $U_m:=(v_m,u_m)$ such that
\begin{align}
	\begin{split}
		&U_- \quad \& \quad U_m \quad \mbox{are connected by 1-shock wave},\nonumber\\
		&U_m \quad \& \quad U_+ \quad \mbox{are connected by 2-rarefaction wave}.\nonumber
	\end{split}
\end{align}
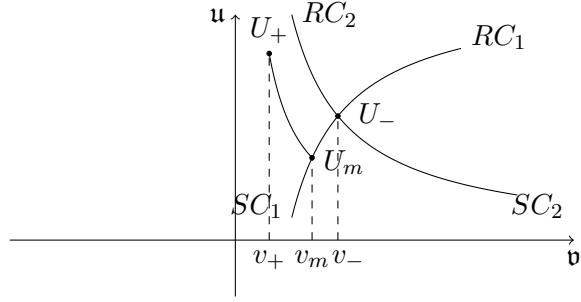
\begin{figure}
	\begin{tikzpicture}[scale=1.5]
		\draw [->] (-2,0)--(3,0);
		\draw [->] (0,-0.5)--(0,2);
		\draw [fill=black] (0.909,1.1) circle [radius=0.02];
        \draw [fill=black] (0.68,0.73) circle [radius=0.02];
        \draw [fill=black] (0.3,1.655) circle [radius=0.02];
        \draw[dashed] (0.909,0) -- (0.909,1.1);
        \draw[dashed] (0.68,0) -- (0.68,0.73);
        \draw[dashed] (0.3,0) -- (0.3,1.655);
       \draw[black,domain=0.5:2.5,samples=500] plot (\x,  {(1/\x)} );
       \draw[black,domain=0.3:0.68,samples=500] plot (\x,  {(0.4964/\x)} );
      \draw[black,domain=0.5:2,samples=500] plot (\x,  {(2.2-1/\x)} );
		\node [below] at (3,0) {$\mathfrak{v}$};
		\node [left] at (0,2) {$\mathfrak{u}$};
		\node [right] at (0.5,2) {$RC_{2}$};
		\node [right] at (2,1.8) {$RC_{1}$};
		\node [left] at (0.51,0.3) {$SC_{1}$};
		\node [left] at (3,0.3) {$SC_{2}$};
        \node [right] at (1,1.1) {$U_{-}$};
        \node [right] at (0.69, 0.69) {$U_{m}$};
        \node [above] at (0.3,1.655) {$U_{+}$};
        \node [below] at (1,0) {$v_{-}$};
        \node [below] at (0.68,0) {$v_{m}$};
        \node [below] at (0.3,0) {$v_{+}$};
	\end{tikzpicture}
	\caption{Riemann Integral Curves}
\end{figure}

\begin{figure}
	\begin{tikzpicture}[scale=1.5]
		\draw [->] (-3,0)--(3,0);
		\draw [->] (0,0)--(0,2);
		\node [below] at (3,0) {$x$};
		\node [left] at (0,2) {$t$};
		\draw [-,black] (0,0)--(-1.8,1.5);
		\node [above] at (-1.8,1.5) {$\mbox{Shock wave}$};
		\node [below] at (-1.7,1.1) {$x=\sigma_1 t$};
		\draw [-,dashed] (0,0)--(1.6,1.3);
		\draw [-,dashed] (0,0)--(1.7,1.2);
		\draw [-,dashed] (0,0)--(1.9,1.1);
		\draw [-] (0,0)--(2,1);
		\draw [-] (0,0)--(1.5,1.4);
		\node [above] at (2.8,1.3) {$\mbox{Rarefaction wave}$};
		\node [above] at (0.8,1.3) {$x=\lambda_2(v_m) t$};
		\node [below] at (2.5,0.9) {$x=\lambda_2(v_+) t$};
		\node [above] at (-1.3,0.1) {$(v_-,u_-)$};
		\node [above] at (0,0.4) {$(v_m,u_m)$};
		\node [above] at (1.4,0.1) {$(v_+,u_+)$};
	\end{tikzpicture}
	\caption{Riemann Solution}
\end{figure}
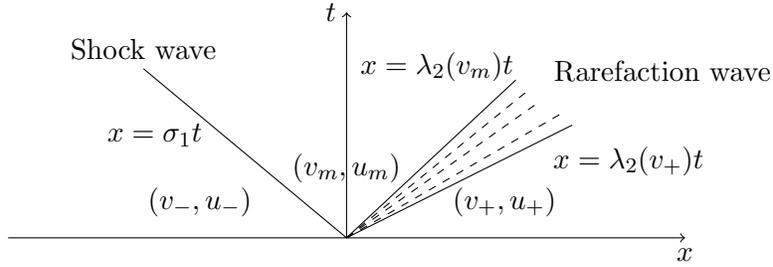

We denote
\begin{align}
	(\mathfrak{v}^s,\, \mathfrak{u}^s)(t,x)	=
	\begin{cases}
		(v_-,u_-),\quad x<\sigma_{1} t,\\
		(v_m,u_m), \quad x>\sigma_{1} t,
	\end{cases}\nonumber
\end{align}
with the shock speed $ \sigma_1=-\sqrt{\f{p(v_{m})-p(v_{-})}{v_{-}-v_{m}}}$. It is direct to know that $(\mathfrak{v}^s,\mathfrak{u}^s)(t,x)$ is a shock solution to \eqref{1.2} with the following Riemann initial data
\begin{align}\label{1.12}
	(\mathfrak{v},\, \mathfrak{u})|_{t=0}
	=
	\begin{cases}
		(v_-,u_-),\quad x<0,\\
		(v_m,u_m), \quad x>0.
	\end{cases}\nonumber
\end{align}
Let $(\mathfrak{v}^r,\, \mathfrak{u}^r)(\frac{x}{t})$ be the rarefaction wave of \eqref{1.2} with the Riemann initial data
\begin{align}
	(\mathfrak{v},\mathfrak{u})|_{t=0}
	=
	\begin{cases}
		(v_m,u_m),\quad x<0,\\
		(v_+,u_+), \quad\,\,x>0,
	\end{cases}\nonumber
\end{align}
that is,
$$
\lambda_{2}\left(\mathfrak{v}^{r}\left(\frac{x}{t}\right)\right)=w^{r}\left(\frac{x}{t}\right),\quad \mathfrak{u}^r=u_{m}-\int_{v_{m}}^{\mathfrak{v}^{r}}\lambda_{2}(s)ds,
$$
with
\begin{equation*}
	w^{r}\left(\frac{x}{t}\right)=\left\{
	\begin{aligned}
		&\lambda_{2}(v_{m}),\,\,\,\,x\leq\lambda_{2}(v_{m})t,\\
		&\frac{x}{t},\,\,\,\,\qquad \lambda_{2}(v_{m})t< x< \lambda_{2}(v_{+})t,\\
		&\lambda_{2}(v_{+}),\,\,\,\,x\geq\lambda_{2}(v_{+})t.
	\end{aligned}
	\right.
\end{equation*}
Now the Riemann solution of \eqref{1.2}-\eqref{1.3} can be written as
\begin{align}\label{1.13}
	\begin{cases}
		\bar{\mathfrak{v}}(t,x):=\mathfrak{v}^s(t,x)+\mathfrak{v}^r(\frac{x}{t})-v_m,\\
		\bar{\mathfrak{u}}(t,x):=\mathfrak{u}^s(t,x)+\mathfrak{u}^r(\frac{x}{t})-u_m.
	\end{cases}
\end{align}

Since both the shock wave $(\mathfrak{v}^s,\mathfrak{u}^s)(t,x)$ and rarefaction wave $(\mathfrak{v}^r,\mathfrak{u}^r)(\frac{x}{t})$ are not smooth functions, we have to mollify these two waves for later use. We define the scaling
\begin{equation}\label{scaling}
\tau=\frac{t}{\nu},\quad y=\frac{x}{\nu},
\end{equation}
and denote $\xi=y-\sigma_1 \tau$. Let $(\tilde{v}^s, \tilde{u}^s)(\xi)=(\tilde{v}_s, \tilde{u}_s)(y-\sigma_{1} \tau)$ be the viscous shock wave connected $(v_-,u_-) \,\, \& \,\, (v_m,u_m)$, i.e.,
\begin{align}\label{2.1}
	\begin{cases}
		\dis -\sigma_{1} (\tilde{v}^s)'-(\tilde{u}^s)'=0,\\
		\dis-\sigma_{1} (\tilde{u}^s)'+ (p(\tilde{v}^s))'=(\frac{b(\tilde{v}^s)}{\tilde{v}^s}(\tilde{u}^s)')',
	\end{cases}
\end{align}
where $\dis '=\f{d}{d\xi}$, see \cite{Matsumura-Wang} for the existence of $(\tilde{v}^{s},\tilde{u}^{s})(\xi)$. Define
\begin{align}\label{2.2}
	(\tilde{\mathfrak{v}}^s, \, \tilde{\mathfrak{u}}^s)(t,x)=(\tilde{v}^s, \tilde{u}^s)(y-\sigma_{1}\tau)=(\tilde{v}^s, \tilde{u}^s)(\frac{x-\sigma_{1} t}{\nu}),
\end{align}
then $(\tilde{\mathfrak{v}}^s, \, \tilde{\mathfrak{u}}^s)(t,x)$ is a viscous shock wave of \eqref{1.1}.

To define the smooth approximate rarefaction wave,  we consider the following Burgers equation
\begin{align}\label{2.3}
	\begin{cases}
		w_t+w w_x=0,\\
		w(0,x)=\frac{1}{2}(w_m+w_+) + \frac{1}{2}(w_+-w_m) \tanh (\f{x}{a}),
	\end{cases}
\end{align}
where $w_m=\lambda_2(v_m),\, w_+=\lambda_2(v_+)$ with $w_m<w_+$. Throughout the present paper, we take the positive parameter $a=a(\nu)>0$ small enough such that $\nu/(a(\nu))\to 0+$ as $\nu\to 0+$. Let $w_a(t,x)$ be the solution of \eqref{2.3}, then we can define the smooth approximate rarefaction wave $(\tilde{v}^r, \tilde{u}^r)(t,x)$ through
\begin{align}\label{2.5}
	\lambda_2(\tilde{v}^r)=w_a(t,x),\qquad \tilde{u}^{r}=u_{m}-\int_{v_{m}}^{\tilde{v}^{r}}\lambda_{2}(s)ds.
\end{align}
Clearly, it holds that
\begin{align}\label{2.6}
	\begin{cases}
		\tilde{v}^r_t-\tilde{u}^r_x=0,\\
		\tilde{u}^r_t+p(\tilde{v}^r)_x=0.
	\end{cases}
\end{align}
Combining \eqref{2.1} and \eqref{2.6}, the smooth approximate solution can be defined as
\begin{equation}\label{2.6-1}
\left\{
\begin{aligned}
	&\tilde{\mathfrak{v}}^{\nu}(t,x)=\tilde{\mathfrak{v}}^{s}(t,x)+\tilde{v}^{r}(t,x)-v_{m},\\ &\tilde{\mathfrak{u}}^{\nu}(t,x)=\tilde{\mathfrak{u}}^{s}(t,x)+\tilde{u}^{r}(t,x)-u_{m}.
\end{aligned}
\right.
\end{equation}

For any given pairs $(v_{1},u_{1})$ and $(v_2,u_2)$, we denote the potential energy by $Q(v):=\frac{v^{-\gamma+1}}{\gamma-1}$ and the relative potential energy by
\begin{equation*}
Q(v_1|v_{2}):=\frac{1}{\gamma-1}v_{1}^{-\gamma+1}-\frac{1}{\gamma-1}v_{2}^{-\gamma+1}+(v_{2})^{-\gamma}(v_{1}-v_2),
\end{equation*}
and the relative total energy as
\begin{equation}\label{1.18}
	\eta\big((v_{1},u_{1})|(v_{2},u_{2})\big):=\frac{1}{2}|u_{1}-u_{2}|^2+Q(v_{1}|v_{2}).
\end{equation}

We denote the shock and rarefaction wave strength respectively
\begin{align}\label{1.16}
	\v_1:=|p(v_-)-p(v_m)|\quad \mbox{and} \quad \v_2:=|p(v_m)-p(v_+)|,
\end{align}
and assume that
\begin{align}\label{1.17}
	\v_2\leq C(v_-) \v_1,
\end{align}
which means that the rarefaction wave strength $\v_2$ is at most the same order as the shock wave strength $\v_1$, and $\v_2$ can be vanished.

The main theorem is stated as follows.
\begin{theorem}\label{thm1.1}
	Assume $\gamma>1$, $0<\alpha\leq \gamma\leq \alpha+ 1$ and \eqref{1.17}. Given an initial data $(\mathfrak{v}_0,\, \mathfrak{u}_0)$ of the compressible Euler system \eqref{1.2} such that
	\begin{align}\label{1.17-1}
		\mathcal{E}_0:=\int_{\R} \eta\big((\mathfrak{v}_0,\, \mathfrak{u}_0)(x)\,|\, (\bar{\mathfrak{v}},\bar{\mathfrak{u}})(0,x)\big) \,dx <\infty.
	\end{align}
	Then, for each $0<\nu \ll 1$, the following statements hold:
	
	\smallskip
	
	\noindent{\rm 1)} There exists a sequence of smooth functions $\{(\mathfrak{v}_0^\nu,\, \mathfrak{u}_0^\nu)\}_{\nu>0}$ such that
	\begin{align}\label{1.17-2}
		\begin{split}
			&\mathfrak{v}_{0}^{\nu}\rightarrow \mathfrak{v}_{0}\quad \text{in }L_{\mathrm{loc}}^{1}(\R)\quad \text{ as }\nu \to 0,\\
			&\lim_{\nu\to 0} (\mathfrak{v}_0^\nu,\, \mathfrak{u}_0^\nu)=(\mathfrak{v}_0,\, \mathfrak{u}_0),\quad \mbox{a.e.},\quad \mbox{and} \quad \mathfrak{v}_0^\nu>0,\\
			&\lim_{\nu\to 0} \int_{\R} \left\{Q(\mathfrak{v}_0^\nu\,| \tilde{\mathfrak{v}}_{0}^\nu) + \f12 \Big[\mathfrak{u}_0^\nu + \nu\f{\gamma}{\alpha} ((p(\mathfrak{v}_0^{\nu}))^{\f{\alpha}{\gamma}})_x -\tilde{\mathfrak{u}}_{0}^{\nu}-\nu\f{\gamma}{\alpha} ((p(\tilde{\mathfrak{v}}_{0}^{\nu}))^{\f{\alpha}{\gamma}})_x\Big]^2\right\}\,dx
			=\mathcal{E}_0,
		\end{split}
	\end{align}
where $(\tilde{\mathfrak{v}}_{0}^{\nu},\tilde{\mathfrak{v}}_{0}^{\nu})(x):=(\tilde{\mathfrak{v}}^{\nu},\tilde{\mathfrak{u}}^{\nu})(0,x)$.
	\smallskip
	
	\noindent{\rm 2)} For any given $T>0$, there exists a smooth solution $(\mathfrak{v}^{\nu},\mathfrak{u}^{\nu})$ to the compressible Navier-Stokes equations \eqref{1.1} with initial data $(\mathfrak{v}_{0}^{\nu},\mathfrak{u}_{0}^{\nu})$. Moreover, there exist limits $\mathfrak{v}_{\infty}\in L^\infty(0,T; L^\infty(\R)+\mathcal{M}(\R))$ and $\mathfrak{u}_{\infty}\in L^{\infty}(0,T;L_{\rm loc}^2(\R)$ such that
	\begin{align}\label{1.17-3}
		&\mathfrak{v}^{\nu}\rightharpoonup\mathfrak{v}_{\infty}, \quad\mathfrak{u}^{\nu} \rightharpoonup\mathfrak{u}_{\infty},\quad \mbox{in}\,\,\, \mathcal{M}_{\mathrm{loc}}((0,T)\times \R),\quad \mbox{as}\,\,\, \nu\to 0+,
	\end{align}
	where $\mathcal{M}_{\mathrm{loc}}((0,T)\times \R)$ is the space of locally bounded Radon measures.
	
	\smallskip
	
	\noindent{\rm 3)} By Radon-Nikodym's theorem, let $\mathfrak{v}_{\infty}=\mathfrak{v}_{a}+\mathfrak{v}_{s}$, where $\mathfrak{v}_{s}$ is the singular part, and $\mathfrak{v}_{a}$ is the regular part with $d\mathfrak{v}_{a}=\mathfrak{m}_{a}dtdx$ and $\mathfrak{m}_{a}\in L_{\mathrm{loc}}^1(\R_{+}^2)$ being the Radon-Nikodym derivative of $d\mathfrak{v}_{a}$ with respect to the Lebesgue measure.
	Then there exist a shift $X_{\infty}\in BV(0,T)$ and a constant $C>0$ such that
	$$
dQ(\mathfrak{v}_{\infty}\,|\, \bar{\mathfrak{v}}_{X_{\infty}}):=Q(\mathfrak{m}_{a}\,|\,\bar{\mathfrak{v}}_{X_{\infty}})\,dtdx+|Q'(\overline{V}(t,x))|\,d\mathfrak{v}_{s}\in L^\infty(0,T; \mathcal{M}(\R)),
	$$
	and for almost every $t \in (0,T)$,
	\begin{align}\label{1.21}
		&\int_{\R} \f12 |\mathfrak{u}_{\infty}-\bar{\mathfrak{u}}_{X_{\infty}}|^2 dx+ \left(\int_{x\in\R}Q(\mathfrak{v}_{\infty}\,|\, \bar{\mathfrak{v}}_{X_{\infty}})(t,dx) \right)\nonumber\\
		&:=\int_{\R} \f12 |\mathfrak{u}_{\infty}-\bar{\mathfrak{u}}_{X_{\infty}}|^2 dx+\int_{\R}Q(\mathfrak{m}_{a}\,|\,\bar{\mathfrak{v}}_{X_{\infty}})\,dx+\int_{\R}|Q'(\overline{V}(t,x))|\mathfrak{v}_{s}(t,dx)\leq C \mathcal{E}_0,
	\end{align}
	where
	\begin{align}\label{1.23}
		\begin{split}
			\bar{\mathfrak{v}}_{X_{\infty}}(t,x)&:=\mathfrak{v}^s(x-\sigma_1 t-X_\infty(t))+\mathfrak{v}^r(t,x)-v_m,\\
			\bar{\mathfrak{u}}_{X_{\infty}}(t,x)&:=\mathfrak{u}^s(x-\sigma_1 t-X_\infty(t))+\mathfrak{u}^r(t,x)-u_m,
		\end{split}
	\end{align}
and
\begin{equation}\label{7.31}
	\overline{V}(t,x):=\left\{
	\begin{aligned}
		&v_{-},\quad \text{for }(t,x)\in \overline{\Omega},\\
		&\mathfrak{v}^{r},\quad\,\,\text{for }(t,x)\in \left(\overline{\Omega}\right)^{c},
	\end{aligned}
	\right.
\end{equation}
with $\Omega:=\{(t,x)\in (0,T)\times \R\,|\,x<\sigma_{1}t+X_{\infty}(t)\}$. Moreover, there exists a constant $C(T)$ depending only on $T$ such that for almost everywhere $t\in (0,T)$,
\begin{align}\label{1.24}
|X_{\infty}(t)|\leq C(T)(\mathcal{E}_{0}+\mathcal{E}_{0}^{\frac{1}{3}}),
\end{align}
and
\begin{align}\label{1.21-1}
		&\int_{\R} \f12 |\mathfrak{u}_{\infty}-\bar{\mathfrak{u}}|^2 dx+\int_{\R}Q(\mathfrak{m}_{a}\,|\,\bar{\mathfrak{v}})\,dx+\int_{\R}|Q'(\overline{V}(t,x))|\mathfrak{v}_{s}(t,dx)\leq C(T)(\mathcal{E}_0+\mathcal{E}_{0}^{\frac{1}{3}}),
\end{align}
where $(\bar{\mathfrak{v}},\bar{\mathfrak{u}})$ is the Riemann solution given in \eqref{1.13}.

\end{theorem}

\begin{corollary}
If $\mathcal{E}_{0}=0$ in \eqref{1.17-1}, i.e. the initial perturbation is zero, then it follows from \eqref{1.21-1} that
$$
\mathfrak{u}_{\infty}(t,x)=\bar{\mathfrak{u}}(t,x),\quad \mathfrak{v}_{\infty}(t,x)=\bar{\mathfrak{v}}(t,x)\qquad \text{for a.e. }(t,x)\in (0,T)\times \R,$$
which implies that
the composite wave of shock and rarefaction to the compressible Euler system \eqref{1.2} is unique in the class of inviscid limit of compressible Navier-Stokes equations \eqref{1.1}.
\end{corollary}

Now, we briefly explain two key points in the proof. Compared with \cite{Kang-Vasseur-2021-Invent,Kang-Vasseur-2020}, the appearance of  rarefaction wave is a new difficulty.
Since it is hard to construct an  approximate rarefaction wave satisfying exactly the Navier-Stokes equations \eqref{1.1}, the error term $\int_{\R}w_{X}\cdot (p(v)-p(\tilde{v}))\cdot F_{1}\,dy$ (see \eqref{3.29}) caused by viscosity related to the inviscid approximate rarefaction $(\tilde{v}^{r},\tilde{u}^{r})(t,x)$ \eqref{2.5} should be analyzed carefully, where $w_{X}$ is the weight function defined in \eqref{3.15-3}. To overcome this difficulty, we first decompose $F_{1}$ as
\begin{align*}
F_{1}&=\partial_{yy}\tilde{v}_X^s\, [b_1(\tilde{v})-b_1(\tilde{v}_X^s)] + |\partial_y \tilde{v}_X^s|^2 \, [b_2(\tilde{v})-b_2(\tilde{v}_X^s)]+ 2b_2(\tilde{v}) \, \partial_y \tilde{v}_X^s\, \partial_y \tilde{v}^r\\
	&\quad +b_1(\tilde{v}) \, \partial_{yy}\tilde{v}^r
	+ b_2(\tilde{v})\, |\partial_y \tilde{v}^r|^2=:\sum\limits_{j=1}^{5}F_{1j},
\end{align*}
where $b_1(s)=s^{\beta} p'(s)$, $b_2(s)=[\beta-(\gamma+1)] s^{\beta-1} p'(s)$. The first three terms $\int_{\R}w_{X}\cdot (p(v)-p(\tilde{v}))\cdot \sum_{j=1}^3F_{1j}$ involve the interaction of viscous shock wave and the approximate rarefaction wave, which can be controlled by decay property of the interaction of shock and rarefaction, see \eqref{5.26}--\eqref{5.41} for details.

For the last two terms $\int_{\R}w_{X}\cdot (p(v)-p(\tilde{v}))\cdot \sum_{j=4}^5F_{1j}$, noting that $p(v)-p(\tilde{v})\sim v^{-\gamma}$ as $v\to 0$, which cannot be controlled by the relative potential energy $Q(v|\tilde{v})$, since it behaves like $v^{-\gamma+1}$ as $v\to 0$. Hence we shall apply a cut-off technique to $p(v)-p(\tilde{v})$, that is, we define \begin{align}
	p(v_k)-p(\tilde{v})
	:=
	\begin{cases}
		-k_1, \qquad\qquad\,  \mbox{if}\,\, p(v)-p(\tilde{v})<-k_1,\\
		p(v)-p(\tilde{v}),\quad\mbox{if} -k_1\leq p(v)-p(\tilde{v})\leq k_2,\\
		k_2, \qquad\qquad\quad \mbox{if}\,\, p(v)-p(\tilde{v})> k_2.
	\end{cases}\nonumber
\end{align}
where $k_{1}:=\frac{1}{2}p(v_{-})$ and $k_{2}=2p(v_{-})+1$. For $-k_{1}\leq p(v)-p(\tilde{v})\leq k_{2}$, by using Sobolev inequality and temporal decay of $\partial_{y}\tilde{v}^{r}$ and $\partial_{yy}\tilde{v}^{r}$, one has
\begin{align*}
&\int_{R}w_{X}\cdot |p(v)-p(\tilde{v})|\cdot (|F_{14}+F_{15}|)\mathbf{1}_{\{-k_{1}\leq p(v)-p(\tilde{v})\leq k_{2}\}}\,dy\\
&\leq C\|p(v_{k})-p(\tilde{v})\|_{L^{\infty}}\int_{\R}|\partial_{y}\tilde{v}^{r}|^2+|\partial_{yy}\tilde{v}^{r}|^2\,dy\\
&\leq C\v_{2}\min\{\frac{\nu \v_{2}}{a(\nu)},\tau^{-1}\}\|\partial_{y}(p(v_{k}-p(\tilde{v}))\|_{L^2}^{\frac{1}{2}}\cdot \|p(v_{k})-p(\tilde{v})\|_{L^2}^{\frac{1}{2}}\\
&\leq \f{C}{\delta_0} \v_2^{\f43} \left(\min\{\f{\nu\v_2}{a},\, \tau^{-1}\}\right)^{\f43}\Big(1+\int_{\R} \eta(U|\tilde{U})dy\Big)\\
	&\quad
	+ C \delta_0 \int_{\R} w_X\cdot  v^{\beta}\cdot  |\partial_y(p(v)-p(\tilde{v}))|^2dy,
\end{align*}
where $\delta_{0}$ is a small enough parameter such that the last term will be absorbed by $D(v)$ in the good term $G_{\delta}(U)$ in \eqref{3.30}, see \eqref{5.46}--\eqref{5.47} and \eqref{5.52}--\eqref{5.53} for details. Here, we emphasize that $\min\{\f{\nu\v_2}{a},\, \tau^{-1}\}^{\frac{4}{3}}\in L^{1}(0,\infty)$ uniformly in $\nu$, which is crucial for us to apply Gr\"{o}nwall inequality to obtain global uniform energy estimate. For $p(v)-p(\tilde{v})\geq k_2$ and $p(v)-p(\tilde{v})\leq -k_{1}$, noting that $|p(v)-p(\tilde{v})|\leq Cp(v|\tilde{v})$, we use the expanding property of rarefaction, i.e., the good term $G^{r}(v)$ brought by the rarefaction wave in \eqref{3.30}, to obtain
$$
\int_{R}w_{X}\cdot |p(v)-p(\tilde{v})|\cdot (|F_{14}+F_{15}|)\mathbf{1}_{\{p(v)-p(\tilde{v})\leq -k_{1}\text{ or } p(v)-p(\tilde{v})\geq k_{2}\}}\,dy\leq C\frac{\nu}{a(\nu)}G^{r}(v),
$$
see \eqref{5.48} and \eqref{5.53-1} for details.

On the other hand, the error term $\int_{\R}w_{X}\cdot (h-\tilde{h})\cdot F_{2}\,dy$ due to the interactions of shock and rarefaction waves is another difficulty in our proof. Noting from \eqref{3.7} that
\begin{align*}
\int_{\R}w_{X}\cdot (h-\tilde{h})\cdot F_{2}\,dy&\leq \Big\vert\int_{\R}w_{X}\cdot (h-\tilde{h})\cdot \partial_{y}\tilde{v}^{r}\cdot [p'(v)-p'(\tilde{v}^{r})]\,dy\Big\vert\\
&\quad +\Big\vert\int_{\R}w_{X}\cdot (h-\tilde{h})\cdot \partial_{y}\tilde{v}_{X}^{s}\cdot [p'(\tilde{v})-p'(\tilde{v}_{X}^{s})]\,dy\Big\vert\\
&=:E_{1}+E_{2}.
\end{align*}
For $E_{1}$, using the good decay property of $\partial_{y}\tilde{v}^{r}$ and the separation of shock and rarefaction waves, one has
\begin{align*}
E_{1}&\leq C\left(\int_{\R} |h-\tilde{h}|^2\, dy\right)^{\f12}
	\left(\int_{\R} |\partial_y \tilde{v}^r|^2\cdot |v^s_X-v_m|^2dy\right)^{\f12}\nonumber\\
	&\leq \Big[\f{\nu}{a}\exp\Big\{-\f{\nu}{C a}\tau\Big\} + \v_1 \exp\Big\{-\f{\v_1}{C}\tau\Big\}\Big]\int_{\R} |h-\tilde{h}|^2 dy\nonumber\\
	&\quad + C\v_1^2\v_2^2 \exp\Big\{-\f{\nu}{C a}\tau\Big\}
	+ C \Big(\min\{\f{\nu\v_2}{a},\, \tau^{-1}\} \Big)^2,
\end{align*}
see \eqref{5.56}--\eqref{5.58} for details. We also point out that $[\f{\nu}{a}\exp\{-\f{\nu}{C a}\tau\} + \v_1 \exp\{-\f{\v_1}{C}\tau\}]\in L^{1}(0,\infty)$ uniformly in $\nu$, which is crucial for us to apply Gr\"{o}nwall inequality to get the global uniform energy estimate.

For $E_{2}$, since $\partial_{y}\tilde{v}_{X}^{s}$ does not have the temporal decay property as $\partial_{y}\tilde{v}^{r}$, if we use the similar argument as $E_{1}$, then we can only get
$$
E_{2}\leq \Big[\exp\Big\{-\f{\nu}{C a}\tau\Big\} + \v_1 \exp\Big\{-\f{\v_1}{C}\tau\Big\}\Big]\int_{\R} |h-\tilde{h}|^2 dy+ C\v_1^2\v_2^2 \exp\Big\{-\f{\nu}{C a}\tau\Big\}
+ C \v_{1}\exp\Big\{-\f{\v_1}{C}\tau\Big\}.
$$
Noting that $\int_{0}^{\infty}e^{-\frac{\nu}{C a}}\tau\,d\tau\sim \frac{a}{\nu}\to +\infty$ as $\nu\to 0$, then one fails to get the global uniform energy estimate by Gr\"{o}nwall inequality. To overcome this difficulty, we introduce a new term $\v G^{h}$ from the good term $G_{\delta}(U)$ in \eqref{3.30}, that is,
$$
\v G^{h}=\v \frac{\sigma_{1}}{2}\int_{\R}\partial_{y}w_{X}\cdot |h-\tilde{h}|^2dy\qquad \text{where }\v=\v_{1}\v_{2},
$$
Then, by using \eqref{3.15-3} and careful analysis on the interaction of shock and rarefaction, we obtain
\begin{align*}
	E_{2}
	&\leq \left(\int_{\R} |\partial_y \tilde{v}^s_X|\cdot |h-\tilde{h}|^2 dy\right)^{\f12}\cdot \left(\int_{\R} |\partial_y \tilde{v}^s_X|\cdot |\tilde{v}^r-v_m|^2 dy \right)^{\f12}\nonumber\\
	&\leq C\sqrt{\f{\v_1\v_2}{\lambda}} \left(\int_{\R} |\partial_yw_X|\cdot |h-\tilde{h}|^2 dy\right)^{\f12}\cdot \left(\int_{\R} |\partial_y \tilde{v}^s_X|\cdot |\tilde{v}^r-v_m| dy \right)^{\f12}\nonumber\\
	&\leq \delta_0 \v_1\v_2 \int_{\R} \f{|\sigma_1|}{2} |\partial_yw_X|\cdot |h-\tilde{h}|^2 dy
	+C \f{\v_1\v_2}{\delta_0\lambda} \left(\exp\big\{-\f{\nu}{C_\ast a}\tau\big\}
	+   \exp\big\{-\f{\v_1}{C}\tau\big\}\right)\nonumber\\
	&\leq \delta_0 \v\, G^h(h)  +C \f{\v_1\v_2}{\delta_0\lambda} \left(\exp\big\{-\f{\nu}{C a}\tau\big\}
	+   \exp\big\{-\f{\v_1}{C}\tau\big\}\right),
\end{align*}
where $\delta_{0}$ is a small enough constant such that $\delta_{0} \v G^{h}$ can be absorbed by $\v G^{h}(h)$, see \eqref{5.58-1}--\eqref{5.60} for details. We point out that if the strength of rarefaction $\v_2=0$, the $\v G^{h}=0$ and it turns back to the case of single shock \cite{Kang-Vasseur-2021-Invent}.

As a conclusion, in the scaling coordinate $(\tau,y)$, one can obtain the global energy inequality
$$
\frac{d}{d\tau}\int_{\R}w_{X}\cdot \eta(U|\tilde{U})\,dy+\text{ some positive terms }\leq 2f(\nu,\tau)\int_{\R}\eta(U|\tilde{U})\,dy+3g(\nu,\tau).
$$
Here, we emphasize $\int_{0}^{\infty}f(\nu,\tau)\,d\tau<C$ is bounded uniformly in $\nu$, which is crucial for us to apply the Gr\"{o}nwall's inequality to obtain uniform energy estimate in the original coordinate $(t,x)$. While,  $\int_{0}^{\infty}g(\nu,\tau)\,d\tau\sim Ca(\nu)/\nu$ is not uniformly bounded in $\nu$ since $a(\nu)/\nu\to \infty$ as $\nu \to 0$. Fortunately,
it is  enough for us to take inviscid limit in the original coordinate $(t,x)$ since $\nu \int_{0}^{\infty}g(\nu,\tau)\,d\tau\sim Ca(\nu)\to 0$ as $\nu\to 0$. See \eqref{6.26}-\eqref{6.30} for details.

The rest of the paper are organized as follows. The Sections 2--5 are devoted to establish the global energy estimate Theorem \ref{thm3.1}.  In Section 2, we first recall some basic properties about the viscous shock wave and approximate rarefaction wave, and then reformulate the problem and present some basic estimates around composite wave. In Section 3, {\it a priori} estimates in a general setting are stated, and in Section 4, the estimates for large perturbation part are established. Then the proof of global energy estimate Theorem \ref{thm3.1} is accomplished in Section 5. Finally, based on Theorem \ref{thm3.1}, we complete the proof of Theorem \ref{thm1.1} in Section 6.

			%


\section{Reformulation and basic energy estimate}

\subsection{Viscous shock wave and approximate rarefaction wave}
For later use, we present some useful properties for the viscous shock wave and approximate rarefaction wave.

\begin{lemma}[\cite{Goodman,Kang-Vasseur-2021-JEMS}]\label{lem2.3}
For the viscous shock wave $(\tilde{v}^{s},\tilde{u}^{s})(y-\sigma_1\tau)$ defined in \eqref{2.1}, it holds that
	\begin{align}\label{2.3-3}
		\begin{split}
			&|\tilde{v}^s(y-\sigma_1 \tau)-v_m|+|\tilde{u}^s(y-\sigma_1 \tau)-u_m|\leq C_1\v_1 e^{-c_1\v_1 |y-\sigma_1 \tau|},\quad \mbox{if}\,\, y\geq \sigma_1 \tau,\\
			&|\partial_y \tilde{v}^s(y-\sigma_1 \tau)| + |\partial_y \tilde{u}^s(y-\sigma_1 \tau)|\leq C_1\v_1^2 e^{-c_1\v_1 |y-\sigma_1 \tau|},\qquad\qquad \forall y\in \R, \, \tau\geq0,\\
			&|\partial_y^2 \tilde{v}^s(y-\sigma_1 \tau)| + |\partial_y^2 \tilde{u}^s(y-\sigma_1 \tau)|\leq C_1\v_1^3 e^{-c_1\v_1 |y-\sigma_1 \tau|},\qquad\qquad \forall y\in \R, \, \tau\geq0,
		\end{split}
	\end{align}
and
\begin{equation}\label{2.3-2}
	\inf_{\xi\in [-\varepsilon_{1}^{-1},\varepsilon_{1}^{-1}]}|\partial_{\xi}\tilde{v}^{s}|\geq C_{1}\varepsilon_{1}^2.
\end{equation}
	where $\xi=y-\sigma_{1}\tau$, $c_1>0$ and $C_1>0$ are some positive constants depending only on $U_-$ and $U_m$.
\end{lemma}

\begin{lemma}[\cite{Xin-1993-CPAM}]\label{lem2.2}
	For $(\tilde{v}^{r},\tilde{u}^{r})(t,x)$ defined in \eqref{2.5}, it holds that
	
	\noindent 1) $\partial_x \tilde{v}^r<0,\,\, \partial_x \tilde{u}^r>0$, and $\|\tilde{v}^r(t,\cdot)-\mathfrak{v}^r(\f{\cdot}{t})\|_{L^1} + \|\tilde{u}^r(t,\cdot)-\mathfrak{u}^r(\f{\cdot}{t})\|_{L^1}\leq C_2\v_2\cdot a$.
	\smallskip
	
	\noindent 2) There exists $a_{0}\in (0,1)$ such that if $a\in(0,a_0)$, then for any $p\in (1,+\infty]$,
	\begin{align}
		&\|\tilde{v}^r(t,\cdot)-\mathfrak{v}^r(\f{\cdot}{t})\|_{L^p} + \|\tilde{u}^r(t,\cdot)-\mathfrak{u}^r(\f{\cdot}{t})\|_{L^p}\nonumber\\
		&\leq C_2 \min\Big\{a^{\f1p}\v_2,\,\, a\,\v_2^{\f1p} t^{-1+\f1p}  (\ln(\f{1+t}{a}))^{1-\f1p}\Big\};\nonumber
	\end{align}
	
	\smallskip
	
	\noindent 3) For all $t>0$, $p\in [1,+\infty]$ and $a(\nu)>0$, it holds
	\begin{align}\label{2.12}
		\begin{split}
			&\|\partial_x(\tilde{v}^r,\tilde{u}^r)(t,\cdot)\|_{L^p}\leq C_2 \min\Big\{\v_2 a^{-1+\f1p},
			\, \v_2^{\f1p} t^{-1+\f1p}\Big\},\\
			&\|\partial_{x}^2(\tilde{v}^r,\tilde{u}^r)(t,\cdot)\|_{L^p}\leq C_2a^{-1+\f1p}\min\Big\{\f{\v_2}{a},
			\, \f{2}{t}\Big\}.
		\end{split}
	\end{align}
	
	\smallskip
	
	\noindent 4) It holds that
	\begin{align}\label{2.13}
		&|\tilde{v}^r(t,x)-v_m|+ |\tilde{u}^r(t,x)-u_m|\leq C_2\v_2 \exp\Big\{-\f{1}{a}|x-\lambda_2(v_m)t|\Big\},\quad \mbox{for}\,\, x\leq \lambda_2(v_m)t,
	\end{align}
	and
	\begin{align}\label{2.14}
		\begin{split}
			&|(\partial_x \tilde{v}^r,\partial_x \tilde{u}^r)(t,x)|\leq C_2 \f{\v_2}{a} \, \exp\Big\{-\f{1}{a}|x-\lambda_2(v_m)t|\Big\}, \\
			&|(\partial_x^2 \tilde{v}^r,\partial_x^2 \tilde{u}^r)(t,x)|\leq C_2 \f{\v_2}{a^2}\, \exp\Big\{-\f{1}{a}|x-\lambda_2(v_m)t|\Big\},
		\end{split}
		\qquad \mbox{for}\,\, x\leq \lambda_2(v_m)t,
	\end{align}
	where $C_2>0$ is a constant depending only on $U_m$ and $U_+$.
\end{lemma}

\subsection{Reformulation of \eqref{1.1}}
Due to the BD entropy ({\it {cf.}} \cite{Bresch-Desjardins-2003,Bresch-Desjardins-2006}), we define
\begin{align}
	\mathfrak{h}^{\nu}:=\mathfrak{u}^{\nu}+ \nu\f{\gamma}{\alpha} ((p(\mathfrak{v}^{\nu}))^{\f{\alpha}{\gamma}})_x,\nonumber
\end{align}
then for smooth solution, the Navier-Stokes equations \eqref{1.1} can be rewritten as
\begin{align}\label{1.15}
	\begin{cases}
		\dis \partial_t \mathfrak{v}^{\nu}- \partial_x \mathfrak{h}^{\nu}=-\nu \dis \partial_x((\mathfrak{v}^{\nu})^\beta \partial_x p(\mathfrak{v}^{\nu})),\\
		\partial_t\mathfrak{h}^{\nu}+\partial_x p(\mathfrak{v}^{\nu})=0,
	\end{cases}
\end{align}
where $\beta=\gamma-\alpha$. Recalling the scaling \eqref{scaling}, we define
\begin{align}
	(v,\, u)(\tau,y):=(\mathfrak{v}^{\nu},\, \mathfrak{u}^{\nu})(t,x),\quad h(\tau,y):=\mathfrak{h}^{\nu}(t,x),\nonumber
\end{align}
then in the coordinate of $(\tau,y)$, the compressible Navier-Stokes equations \eqref{1.15} becomes
\begin{align}\label{1.15-1}
	\begin{cases}
		\partial_\tau v- \partial_y h=- \partial_y(v^\beta \partial_y p),\\
		\partial_\tau h+\partial_y p(v)=0.
	\end{cases}
\end{align}
Also, for the viscous shock wave, we denote
\begin{equation*}
	\tilde{h}^s:= \tilde{u}^s+\f{\gamma}{\alpha} ((p(\tilde{v}^s))^{\f{\alpha}{\gamma}})_y,
\end{equation*}
with $(\tilde{v}^{s},\tilde{u}^{s})$  being the viscous shock wave given in \eqref{2.1}, then it is direct to know
\begin{align}\label{3.15-2}
	\begin{cases}
		\partial_\tau \tilde{v}^s -\partial_y \tilde{h}^s=-\partial_y((\tilde{v}^s)^\beta \partial_y p(\tilde{v}^s)),\\
		\partial_\tau \tilde{h}^s + \partial_y p(\tilde{v}^s)=0.
	\end{cases}
\end{align}

Let $X:=X(\tau)$ be the shift, we define
\begin{align}
	(\tilde{v}^s_X, \tilde{h}^{s}_X):=(\tilde{v}^s, \tilde{h}^{s})(y-\sigma \tau-X(\tau)),\nonumber
\end{align}
which yields that
\begin{align}\label{3.4}
	\begin{cases}
		\partial_\tau \tilde{v}^s_{X} -\partial_y \tilde{h}^s_X + \dot{X}(\tau) \cdot \partial_y\tilde{v}^s_X=-\partial_y((\tilde{v}^s_X)^\beta \partial_y p(\tilde{v}^s_X)),\\
		\partial_\tau \tilde{h}^s_X + \partial_y p(\tilde{v}^s_X) + \dot{X}(\tau)\cdot \partial_y \tilde{h}^s_X =0.
	\end{cases}
\end{align}

Since the rarefaction wave does not have viscous term, we define the composite wave as
\begin{align}
	\tilde{v}:= \tilde{v}^s_{X} + \tilde{v}^r - v_m,\qquad \tilde{h}:=\tilde{h}^s_X + \tilde{u}^r-u_m,\nonumber
\end{align}
which, together with \eqref{2.6} and \eqref{3.4}, yields that
\begin{align}\label{3.6}
	\begin{cases}
		\partial_\tau\tilde{v} - \partial_y\tilde{h} + \dot{X}(\tau)\cdot \partial_y \tilde{v}^s_X = - \partial_y(\tilde{v}^\beta \partial_y p(\tilde{v})) + F_1,\\
		\partial_\tau \tilde{h} + \partial_y p(\tilde{v}) + \dot{X}(\tau)\cdot \partial_y\tilde{h}^s_X= F_2,
	\end{cases}
\end{align}
where
\begin{align}\label{3.7}
	\begin{split}
		F_1&:=\partial_y\big(\tilde{v}^\beta \partial_y p(\tilde{v})-(\tilde{v}^s_X)^\beta \partial_y p(\tilde{v}^s_X)\big),\\
		F_2&:=\partial_y\big(p(\tilde{v})-p(\tilde{v}^r)-p(\tilde{v}^s_X)\big).
	\end{split}
\end{align}

\

We rewrite \eqref{1.15-1} into the viscous hyperbolic system of conservation laws
\begin{align}\label{3.8}
	\partial_\tau U + \partial_y A(U)
	=
	\left(
	\begin{array}{c}
		-(v^{\beta} \partial_y p)_y \\
		0
	\end{array}
	\right),
\end{align}
where
\begin{align}
	U:=\left(
	\begin{array}{c}
		v\\
		h
	\end{array}\right),
	\qquad
	A(U):=
	\left(
	\begin{array}{c}
		-h\\
		p
	\end{array}\right).\nonumber
\end{align}
The entropy of \eqref{3.8} is
\begin{align}
	\eta(U)=\f12 h^2 + Q(v)\quad \mbox{with}\quad Q(v):=\frac{v^{-\gamma+1}}{\gamma-1},\nonumber
\end{align}
and the corresponding entropy flux is
\begin{align}
	G(U):=p(v)\, h.\nonumber
\end{align}
Clearly, it holds
\begin{align}
	\nabla_{U}\eta(U)=
	\left(
	\begin{array}{c}
		-p(v)\\
		h
	\end{array}
	\right).\nonumber
\end{align}
We denote the matrix
\begin{align}
	M(U):=
	\left(
	\begin{array}{cc}
		v^{\beta}\,\, &0\\
		0 & 0
	\end{array}
	\right),\nonumber
\end{align}
Hence \eqref{3.8} can be rewritten as
\begin{align}
	\partial_\tau U + \partial_y A(U)=\partial_y (M(U)\partial_y \nabla_{U} \eta(U)).\nonumber
\end{align}
Similarly, we can rewrite \eqref{3.6} as
\begin{align}
	\partial_\tau \tilde{U} + \partial_y A(\tilde{U}) = \partial_y(M(\tilde{U})\partial_y \nabla_{U}\eta(\tilde{U})) - \dot{X}(\tau)\partial_y(\tilde{U}^s_X) +
	\left(\begin{array}{c}
		F_1\\
		F_2
	\end{array}
	\right)\nonumber
\end{align}
where $\tilde{U}:=(\tilde{v},\tilde{h})^{t}$, $\tilde{U}^s:=(\tilde{v}^{s}, \tilde{h}^{{\normalsize {\small }}s})^{t}$ and $\tilde{U}^s_X:=\tilde{U}^s(y-\sigma_1 \tau- X(\tau))$.

For later use, we also define the following relative functions

\begin{align}
	\begin{split}
		\eta(U|V):&=\eta(U)-\eta(V)-\nabla_{U}\eta(V)\cdot (U-V),\\
		A(U|V):&=A(U)-A(V)-\nabla_{U} A(V)\cdot (U-V),\\
		Q(v|\tilde{v}):&=Q(v)-Q(\tilde{v})-Q'(\tilde{v}) (v-\tilde{v}),\\
		p(v|\tilde{v}):&=p(v)-p(\tilde{v})-p'(\tilde{v}) (v-\tilde{v}),\nonumber
	\end{split}
\end{align}
and the relative entropy-flux
\begin{align}
	G(U;V):= G(U)-G(V)-\nabla_{U}\eta(V)\cdot (A(U)-A(V)).\nonumber
\end{align}
A direct calculation shows that
\begin{align}
	\begin{split}
		\eta(U|\tilde{U})&\equiv\f12 |h-\tilde{h}|^2 + Q(v|\tilde{v}),\\
		G(U;\tilde{U})&\equiv(p(v)-p(\tilde{v}))(h-\tilde{h}),
	\end{split}\nonumber
\end{align}
and
\begin{align}
	A(U|\tilde{U})\equiv
	\left(
	\begin{array}{c}
		0\\
		p(v|\tilde{v})
	\end{array}
	\right).\nonumber
\end{align}

Let us review some properties of $Q(v|w)$ and $p(v|w)$ for later use.
\begin{lemma}[{\cite[Lemma 2.4]{Kang-Vasseur-2021-JEMS}}]\label{lem2.4}
	For given constants $\gamma>1$, and $v_{-}>0$, there exist constants $\tilde{c}_{1}, \tilde{c}_{2}>0$ such that the following inequalities hold:
	\begin{itemize}
	\item[(1)] For any $w \in\left(0, 2v_{-}\right)$,
	\begin{align}\label{A.1}
		\begin{split}
			&Q(v | w) \geq \tilde{c}_{1}|v-w|^{2} \quad \text { for all } 0<v \leq 3 v_{-}, \\
			&Q(v | w) \geq \tilde{c}_{2}|v-w| \quad\,\,\text { for all } v \geq 3 v_{-}.
		\end{split}
	\end{align}
	\item[(2)] Moreover, if $0<w \leq u \leq v$ or $0<v \leq u \leq w$ then
	\begin{align}\label{A.1-1}
		Q(v | w) \geq Q(u | w),
	\end{align}
	and for any $\delta_{*}>0$, there exists a constant $C>0$ such that if, in addition, $v_{-}>w>v_{-}-\delta_{*} / 2$ and $|w-u|>\delta_{*}$, we have
	\begin{equation}\label{A.3}
		Q(v | w)-Q(u | w) \geq C|u-v|.
	\end{equation}
\end{itemize}
\end{lemma}

\begin{remark}
	By similar arguments as in the proof of \cite[Lemma 2.4]{Kang-Vasseur-2021-JEMS}, we obtain that \eqref{A.1}-\eqref{A.3} also hold for $p(v|w)$.
\end{remark}

\subsection{Basic energy estimate}
Now, we consider the uniform energy estimate around the composition wave. Following \cite{Kang-Vasseur-2021-Invent,Kang-Vasseur-2020}, we define the weight function
\begin{align}\label{3.15-4}
	w(\tau,y)=w(y-\sigma_1\tau):=1-\f{\lambda}{\v_1} [p(\tilde{v}^s(y-\sigma_1\tau))-p(v_-)],
\end{align}
where $\v_{1}$ is defined \eqref{1.16} and $\lambda$ is a small positive parameter such that $0<w<1$. Clearly, it follows from $\eqref{3.15-2}_{2}$ and \eqref{3.15-4} that
\begin{align}\label{3.15-3}
	\partial_yw=-\f{\lambda}{\v_1} \partial_yp(\tilde{v}^s)=-\f{\lambda}{\v_1}\sigma_1\partial_{y}\tilde{h}^{s}<0.
\end{align}
For later use, we denote $w_X:=w(y-\sigma_1 \tau-X(\tau))$.

\smallskip

The main context of  present paper is to establish a global estimate, that is, Theorem \ref{thm3.1}. The proof of Theorem \ref{thm3.1} is very complicate, in fact, its proofs are included in sections 3-6.

\smallskip

Firstly, by similar arguments as in \cite{Kang-Vasseur-2021-Invent,Kang-Vasseur-Wang-2021}, we have
\begin{lemma}\label{lem3.1}
	Let $\gamma>1$, it holds that
	\begin{align}
		\frac{d}{d\tau}\int_{\R}  w_{X} \cdot \eta(U|\tilde{U})(\tau,y)\, dy=\dot{X}(\tau)\cdot Y(U) + J^{\mbox{bad}}- J^{\mbox{good}},\nonumber
	\end{align}
	where
	\begin{align}\label{3.22}
		Y(U)&:=-\int_{\R} \partial_yw_{X} \cdot \eta(U|\tilde{U})\, dy + \int_{\R} w_{X}\cdot  (U-\tilde{U})^{t} \, \nabla^2\eta(\tilde{U}) \, \partial_y (\tilde{U}^{s}_{X}) dy\nonumber\\
		&\equiv -\int_{\R}\f12 \partial_y w_X\, |h-\tilde{h}|^2 \, dy  -  \int_{\R}  \partial_y w_X\,Q(v|\tilde{v}) \,dy\nonumber\\
		&\quad - \int_{\R} w_X\, p'(\tilde{v})\, \partial_y\tilde{v}^s_X\, (v-\tilde{v})\, dy + \int_{\R} w_X\,  \partial_y\tilde{h}^s_X\, (h-\tilde{h})\, dy,
	\end{align}
	\begin{align}
		J^{\mbox{bad}}:&=\int_{\R} w_{X}\, (p(v)-p(\tilde{v}))\cdot (h-\tilde{h}) dy+\sigma_1 \int_{\R} w_{X}\cdot  \partial_y\tilde{v}^s_X\cdot p(v|\tilde{v}) dy\nonumber\\
		&\quad-\int_{\R} \partial_yw_{X}\cdot v^{\beta}\cdot (p(v)-p(\tilde{v}))\cdot \partial_y(p(v)-p(\tilde{v})) \, dy\nonumber\\
		&\quad - \int_{\R} w_{X}\cdot  \partial_y(p(v)-p(\tilde{v}))\cdot (v^{\beta}-\tilde{v}^{\beta}) \cdot \partial_yp(\tilde{v}) \, dy\nonumber\\
		&\quad + \int_{\R} w_{X}\, \Big[ (p(v)-p(\tilde{v})) \cdot F_1 - (h-\tilde{h})\cdot F_2\Big] \, dy,\nonumber
	\end{align}
	and
	\begin{align}
		J^{\mbox{good}}:&=\sigma_1 \int_{\R}  \partial_yw_{X}\cdot \Big[\f12 |h-\tilde{h}|^2 + Q(v|\tilde{v})\Big]\, dy + \int_{\R} w_X\cdot \partial_y \tilde{u}^r\cdot p(v|\tilde{v}) \, dy\nonumber\\
		&\quad + \int_{\R} \big|\partial_yw_{X}\cdot \partial_yp(\tilde{v}) \cdot   (p(v)-p(\tilde{v}))\cdot (v^{\beta}-\tilde{v}^{\beta})\big|  \, dy\nonumber\\
		&\quad + \int_{\R} w_X\cdot v^{\beta}\cdot |\partial_y(p(v)-p(\tilde{v}))|^2 dy.\nonumber
	\end{align}
\end{lemma}


\begin{lemma}\label{lem3.3}
	Denote $\v:=\v_1\v_2$, for any $\delta>0$, it holds that
	\begin{align}
		J^{\mbox{bad}}- J^{\mbox{good}}=B_{\delta}(U)-G_{\delta}(U),\label{3.7-1}
	\end{align}
	with
	\begin{align}\label{3.29}
		B_{\delta}(U):&= \sigma_1  \int_{\R}  w_X\cdot \partial_y \tilde{v}^s_X \cdot p(v|\tilde{v})\, dy\nonumber\\
		&\quad +\int_{\R} \partial_y w_X\cdot (p(v)-p(\tilde{v}))\, (h-\tilde{h})\cdot {\bf 1}_{\{p(v)-p(\tilde{v})>\delta\}} \, dy\nonumber\\
		&\quad + \frac{1}{1-\v} \f{1}{2\sigma_1}  \int_{\R} \partial_y w_X \, |p(v)-p(\tilde{v})|^2  {\bf 1}_{\{p(v)-p(\tilde{v})\leq \delta\}} \, dy\nonumber\\
		&\quad-\int_{\R} \partial_yw_{X}\cdot v^{\beta}\cdot (p(v)-p(\tilde{v}))\cdot \partial_y(p(v)-p(\tilde{v})) \, dy\nonumber\\
		&\quad - \int_{\R} w_{X}\cdot  \partial_y(p(v)-p(\tilde{v}))\cdot (v^{\beta}-\tilde{v}^{\beta}) \cdot \partial_yp(\tilde{v}) \, dy\nonumber\\
		&\quad + \int_{\R} w_{X}\, \big[ (p(v)-p(\tilde{v})) \cdot F_1 - (h-\tilde{h})\cdot F_2\big] \, dy\nonumber\\
		&=: B_1(v)+B_{2}^-(U) + \frac{1}{1-\v} \, B_{2}^+(v) + B_3(v) +B_4(v) + B_5(U),
	\end{align}
	and
	\begin{align}\label{3.30}
		G_{\delta}(U):&=(1-\v)\frac{\sigma_1}{2} \int_{\R}  \partial_yw_{X}\cdot   \Big|h-\tilde{h}-\frac{p(v)-p(\tilde{v})}{(1-\v)\sigma_1}\Big|^2 \cdot {\bf 1}_{\{p(v)-p(\tilde{v})\leq \delta\}} \, dy \nonumber\\
		&\quad + (1-\v)\frac{\sigma_1}{2} \int_{\R}  \partial_yw_{X}\cdot   |h-\tilde{h}|^2 \cdot {\bf 1}_{\{p(v)-p(\tilde{v})>\delta\}} \, dy
		\nonumber\\
		&\quad + \sigma_1 \int_{\R}  \partial_yw_{X}\cdot  Q(v|\tilde{v}) \, dy+ \int_{\R} w_X\cdot \partial_y \tilde{u}^r\cdot p(v|\tilde{v}) \, dy\nonumber\\
		&\quad + \v\frac{\sigma_1}{2} \int_{\R}  \partial_yw_{X}\cdot   |h-\tilde{h}|^2\, dy + \int_{\R} w_X\cdot v^{\beta}\cdot |\partial_y(p(v)-p(\tilde{v}))|^2 dy\nonumber\\
		&\quad + \int_{\R} \big|\partial_yw_{X}\cdot \partial_yp(\tilde{v}) \cdot   (p(v)-p(\tilde{v}))\cdot (v^{\beta}-\tilde{v}^{\beta})\big|  \, dy\nonumber\\
		&=: (1-\v) \, G_1^+(U) + (1-\v) \, G_1^-(U) + G_2(v) + G^r(v) + \v \, G^h(h) + D(v) + G_3(v).
	\end{align}
\end{lemma}

\begin{remark}\label{rem3.3}
	The decomposition of the good term $G_{\delta}(U)$ is slightly different from the one in \cite{Kang-Vasseur-2021-Invent,Kang-Vasseur-2020}. In fact, due to appearance of rarefaction wave, as mentioned in previous section, we shall distribute a $\varepsilon-$part in $G^{h}$ to absorb the some bad terms in $B_{5}(U)$, see the second step in the proof of Lemma \ref{lem5.2} for detailed calculations.
\end{remark}

\begin{lemma}\label{lem3.4}
	Assume $Y(U)\leq \v_1^2$, then there exist positive constants $\tilde{C}_0>0$ and $\tilde{C}_1>0$ such that
	\begin{align}\label{3.31}
		\int_{\R}|\partial_y w_X| \, \eta(U|\tilde{U}) \, dy \leq \tilde{C}_0 \frac{\v_1^2}{\lambda},
	\end{align}
	and
\begin{align}\label{3.32}
	|Y(U)|\leq \tilde{C}_1 \frac{\v_1^2}{\lambda},
\end{align}
where $\tilde{C}_0$ and $\tilde{C}_1$ depend only on $U_-$.
\end{lemma}

\noindent{\bf Proof.} Using \eqref{A.1}, we have
\begin{align}\label{3.33}
	\int_{\R}|\partial_y w_X| \, \eta(U|\tilde{U}) \, dy
	&=\f12 \int_{\R}|\partial_y w_X| \, |h-\tilde{h}|^2 dy + \int_{\R}|\partial_y w_X| \, Q(v|\tilde{v}) \, dy \nonumber\\
	&\geq \f12 \int_{\R}|\partial_y w_X| \, |h-\tilde{h}|^2 dy + c_1 \int_{0<v<3v_-}|\partial_yw_X| \, |v-\tilde{v}|^2 \, dy\nonumber\\
	&\quad + c_2 \int_{v\geq 3v_-}|\partial_y w_X| \cdot |v-\tilde{v}|  \, dy.
\end{align}
Noting \eqref{3.15-3} and $Y(U)\leq \v_1^2$, we have
\begin{align}\label{3.34}
	&\int_{\R}|\partial_y w_X| \, \eta(U|\tilde{U}) \, dy=Y(U) - \int_{\R} w_{X}\cdot  (U-\tilde{U})^{t} \, \nabla^2\eta(\tilde{U}) \, \partial_y (\tilde{U}^{s}_{X}) dy\nonumber\\
	&\leq \v_1^2 + C\frac{\v_1}{\lambda} \int_{\R}|\partial_y w_X| (|h-\tilde{h}| + |v-\tilde{v}|) dy\nonumber\\
	&\leq \v_1^2 + C\frac{\v_1}{\lambda}  \int_{v\geq 3v_-}|\partial_y w_X| \cdot |v-\tilde{v}|  \, dy \nonumber\\
	&\quad  + C\frac{\v_1}{\lambda}\Big(\int_{\R} |\partial_y w_X|\, dy\Big)^{\f12}\Big(\int_{0<v<3v_-}  |\partial_y w_X|\cdot|v-\tilde{v}|^2 \, dy+\int_{\R} |\partial_y w_X|\cdot |h-\tilde{h}|^2 dy\Big)^{\f12} \nonumber\\
	&\leq \v_1^2 + C\frac{\v_1}{\lambda}  \int_{v\geq 3v_-}|\partial_y w_X| \cdot |v-\tilde{v}|  \, dy\nonumber\\
	&\quad + C \frac{\v_1}{\sqrt{\lambda}}\Big(\int_{0<v<3v_-}  |\partial_y w_X|\cdot|v-\tilde{v}|^2 \, dy+\int_{\R} |\partial_y w_X|\cdot |h-\tilde{h}|^2 \, dy\Big)^{\f12}\nonumber\\
	&\leq \v_1^2 + C\frac{\v_1^2}{\lambda} + C\frac{\v_1}{\lambda}  \int_{v\geq 3v_-}|\partial_y w_X| \cdot |v-\tilde{v}|  \, dy+ \frac{c_1}{4}  \int_{0<v<3v_-}|\partial_y w_X| \, |v-\tilde{v}|^2 \, dy\nonumber\\
	&\quad + \f14 \int_{\R} |\partial_y w_X|\cdot |h-\tilde{h}|^2 \, dy\nonumber\\
	& \leq \v_1^2 + C\frac{\v_1^2}{\lambda} + \frac{1}{2}\int_{\R}|\partial_y w_X| \, \eta(U|\tilde{U}) \, dy,
\end{align}
where we have used \eqref{3.33} and chosen $C\frac{\v_1}{\lambda}\leq \delta_0\leq \frac{c_2}{8}$ in the last inequality. Then we have from \eqref{3.34} that
\begin{align*}
	\int_{\R}|\partial_y w_X| \, \eta(U|\tilde{U}) \, dy
	&\leq  \tilde{C}_0\frac{\v_1^2}{\lambda},
\end{align*}
which yields \eqref{3.31}.

For \eqref{3.32}, using \eqref{3.15-3} and \eqref{3.34}, we have
\begin{align*}
	|Y(U)|&\leq  \int_{\R}|\partial_y w_X| \, \eta(U|\tilde{U}) \, dy + \left|\int_{\R} w_{X}\cdot  (U-\tilde{U})^{t} \, \nabla^2\eta(\tilde{U}) \, \partial_y (\tilde{U}^{s}_{X})\, dy \right|\nonumber\\
	&\leq  \int_{\R}|\partial_y w_X| \, \eta(U|\tilde{U}) \, dy + C\frac{\v_1}{\lambda} \int_{\R}|\partial_y w_X| (|h-\tilde{h}| + |v-\tilde{v}|) dy\nonumber\\
	&\leq C \int_{\R}|\partial_y w_X| \, \eta(U|\tilde{U}) \, dy+C\frac{\v_{1}^2}{\lambda}\leq  \tilde{C}_1\frac{\v_1^2}{\lambda},\nonumber
\end{align*}
which yields \eqref{3.32}. Therefore the proof of Lemma \ref{lem3.4} is complete. $\hfill\Box$

\

It follows from \eqref{3.22} that
\begin{align}
	Y(U)
	&=\frac{1}{1-\v} \f{1}{\sigma_1}  \int_{\R} w_X\,  \partial_y\tilde{h}^s_X\, [p(v)-p(\tilde{v})]\, dy-\frac{1}{ (1-\v)^2} \frac{1}{2\sigma_1^2} \int_{\R} \partial_yw_X \, |p(v)-p(\tilde{v})|^2 dy\nonumber\\
	&\quad + \int_{\R} w_X\,  \partial_y\tilde{h}^s_X\, \Big(h-\tilde{h}-\frac{p(v)-p(\tilde{v})}{\sigma_1 (1-\v)}\Big)\, dy
	-  \int_{\R}  \partial_y w_X\,Q(v|\tilde{v}) \,dy\nonumber\\
	&\quad -\int_{\R} w_X\, [p'(\tilde{v})-p'(\tilde{v}_X^s)]\, \partial_y\tilde{v}^s_X\, (v-\tilde{v})\, dy - \int_{\R} w_X\, p'(\tilde{v}_X^s)\, \partial_y\tilde{v}^s_X\, (v-\tilde{v})\, dy\nonumber\\
	&\quad -\f12\int_{\R} \partial_yw_X\,   \Big(h-\tilde{h}-\frac{p(v)-p(\tilde{v})}{\sigma_1 (1-\v)}\Big)\,\Big(h-\tilde{h}+\frac{p(v)-p(\tilde{v})}{\sigma_1 (1-\v)}\Big) dy.\nonumber
\end{align}
We define
\begin{align}\label{3.36}
	\mathcal{Y}_g(v):=&-\frac{1}{ (1-\v)^2} \frac{1}{2\sigma_1^2} \int_{\R} \partial_yw_X \, |p(v)-p(\tilde{v})|^2 dy - \int_{\R}  \partial_y w_X\,Q(v|\tilde{v}) \,dy\nonumber\\
	&-\int_{\R} w_X\, p'(\tilde{v}_X^s)\, \partial_y\tilde{v}^s_X\, (v-\tilde{v})\, dy +\frac{1}{1-\v} \f{1}{\sigma_1} \int_{\R}  w_X\,  \partial_y\tilde{h}^s_X\, [p(v)-p(\tilde{v})]\, dy\nonumber\\
	=:&\f{1}{(1-\v)^2} \mathcal{Y}_{g1}(v) + \mathcal{Y}_{g2}(v) + \mathcal{Y}_{g3}(v) + \f{1}{1-\v} \mathcal{Y}_{g4}(v).
\end{align}
Then we can rewrite $Y(U)$ as
\begin{align}
	Y(U)&=\mathcal{Y}_g(v) +  \int_{\R} w_X\,  \partial_y\tilde{h}^s_X\, \Big(h-\tilde{h}-\frac{p(v)-p(\tilde{v})}{\sigma_1 (1-\v)}\Big)\, dy\nonumber\\
	&\quad -\f12\int_{\R} \partial_yw_X\,   \Big(h-\tilde{h}-\frac{p(v)-p(\tilde{v})}{\sigma_1 (1-\v)}\Big)\,\Big(h-\tilde{h}+\frac{p(v)-p(\tilde{v})}{\sigma_1 (1-\v)}\Big) dy\nonumber\\
	&\quad -\int_{\R} w_X\, [p'(\tilde{v})-p'(\tilde{v}_X^s)]\, \partial_y\tilde{v}^s_X\, (v-\tilde{v})\, dy.\nonumber
\end{align}

For later use, we consider a truncation of $|p(v)-p(\tilde{v})|$ with a parameter $\delta_1>0$. Let $\bar{\psi}(y)$ be a continuous function defined by
\begin{equation}
	\bar{\psi}(y)=\inf(\delta_{1},\sup(-\delta_{1},y)).\nonumber
\end{equation}
Since the function $p$ is one to one, then $\bar{v}$ can be uniquely defined by
\begin{equation}
p(\bar{v})-p(\tilde{v})=\bar{\psi}(p(v)-p(\tilde{v})).\nonumber
\end{equation}

By using similar arguments in \cite[Lemma 3.2]{Kang-Vasseur-2021-JEMS}, we have following lemma.

\begin{lemma}[{\cite[Lemma 3.2]{Kang-Vasseur-2021-JEMS}}]\label{lem3.4-1}
	For fixed $v_{-}>0$, there exist constants $\delta_{2}>0$ and $\tilde{C}_{2}>0$ which may depend on $\delta_{2}$ such that whenever $|Y(v)|\leq \v_{1}^2$, it holds
	\begin{equation}
		|\mathcal{Y}_{g}(\bar{v})|\leq \tilde{C}_{2}\frac{\v_{1}^2}{\lambda},\qquad \text{for }0<\delta_1\leq \delta_{2}.\nonumber
	\end{equation}
\end{lemma}

The following lemma plays a crucial role in establishing in Theorem \ref{thm3.1}. Although the definition of $\mathcal{Y}_g(v)$ is slightly different from the ones in  \cite{Kang-Vasseur-2021-JEMS,Kang-Vasseur-Wang-2021}, the proof is almost the same as \cite[Proposition 3.4]{Kang-Vasseur-2021-JEMS} and \cite[Lemma 4.5]{Kang-Vasseur-Wang-2021}. Therefore, we omit the proof for simplicity of presentation.

\begin{lemma}[\cite{Kang-Vasseur-2021-JEMS}]\label{lem3.5}
	For any given positive constant $\tilde{C}_2>0$,
	there exist small positive constants $\v_0, \d_1>0$ such that for any $\v_1, \v_2\in(0,\v_0]$ with \eqref{1.17}, and $\lambda,\delta\in(0,\d_1]$, the following is true:  For any function $v:\R\to \R_+$ such that $\mathcal{D}(v)+ \mathcal{G}_2(v)$ is finite, if
	\begin{align}
		|\mathcal{Y}_g(v)|\leq \tilde{C}_2 \frac{\v_1^2}{\lambda}\quad \mbox{and}\quad \|p(v)-p(\tilde{v})\|_{L^\infty}\leq 2\delta_1,\nonumber
	\end{align}
then we have
	\begin{align}
		\mathcal{R}_{\delta}(v):&=-\f{1}{\v_1\delta} |\mathcal{Y}_g(v)|^2 + \mathcal{I}_1(v) +\delta |\mathcal{I}_1(v)| + \mathcal{I}_2(v) + \delta \frac{\v_1}{\lambda} |\mathcal{I}_2(v)|\nonumber\\
		&\quad - (1-\delta \f{\v_1}{\lambda})\, \mathcal{G}_2(v) -(1-\delta) \mathcal{D}(v)\leq 0,\nonumber
	\end{align}
	where
	\begin{align}\label{3.40}
		\begin{split}
			&\mathcal{I}_1(v):=\sigma_1  \int_{\R}  w_X\cdot \partial_y \tilde{v}^s_X \cdot p(v|\tilde{v})\, dy,\quad
			\mathcal{I}_2(v):=\f{1}{2\sigma_1}  \int_{\R} \partial_y w_X \, |p(v)-p(\tilde{v})|^2\, dy,\\
			&\mathcal{G}_2(v):= \sigma_1   \int_{\R} \partial_y w_X \, Q(v|\tilde{v})\, dy,\qquad \quad
			\mathcal{D}(v):= \int_{\R} w_X \, v^{\beta}\, |\partial_y(p(v)-p(\tilde{v}))|^2\, dy,
		\end{split}
	\end{align}
	We point out  that the constants $\v_0$ and $\delta_1$ depend on $\tilde{C}_2$.
\end{lemma}

\section{A priori estimate for a general setting}
The following lemma \ref{lem4.1} is a simple generalization of \cite[Lemmas 4.5-4.8]{Kang-Vasseur-2021-Invent}, see also \cite[Proposition 6.3]{Kang-Vasseur-2020}. Since the proof is almost the same, we omit it here for simplicity.

\begin{lemma}\label{lem4.1}
	Let $U_{-}:=\left(v_{-}, u_{-}\right) \in \mathbb{R}^{+} \times \mathbb{R}$ be  a given constant state with $v_{-}>0$. Let $\tilde{U}_{0}:=(\tilde{v}_{0}, \tilde{h}_{0}): \mathbb{R} \rightarrow \mathbb{R}^{+} \times \mathbb{R}$, and $\mathbf{w}: \mathbb{R} \rightarrow \mathbb{R}$ be any functions such that
	\begin{align}
		&\left|\tilde{U}_{0}(y)-U_{-}\right| \leq C_\ast \delta_{0}, \quad \tilde{v}_{0}(y) \geq C_\ast^{-1}, \quad \forall y \in \mathbb{R},\label{4.1-1}\\
		&\left|(\tilde{v}_{0})_{y}\right| \leq \delta_{0}^{2}, \quad|\mathbf{w}(y)| \leq C_\ast \v_1 \lambda \exp \left(-C_\ast^{-1} \varepsilon_1|y|\right), \quad \forall y \in \mathbb{R},\label{4.2}\\
		&\inf_{-\v_1^{-1} \leq y \leq \v_1^{-1}}|\mathbf{w}(y)| \geq C_\ast^{-1}\v_1 \lambda, \quad \int_{\mathbb{R}}|\mathbf{w}| d y=\lambda,\label{4.3}
	\end{align}
	where $C_\ast>0$ is a positive constant depending only on $U_-$ and hereafter $\delta_{0}$ and $\lambda$ are positive two parameters.
	
	\medskip
	
	Let $U:=(v, h): \mathbb{R} \rightarrow \mathbb{R}^{+} \times \mathbb{R}$ be any function such that
	\begin{align}\label{4.3-1}
      \int_{\mathbb{R}}|\mathbf{w}|\cdot |h-\tilde{h}_{0}|^{2} dy
	+\int_{\mathbb{R}}|\mathbf{w}| Q\left(v|\tilde{v}_{0}\right) d y \leq C \frac{\v_1^{2}}{\lambda}.
	\end{align}
	Let $\bar{\bf v}$ be a $\delta_{1}$-truncation of $v$ defined by
	\begin{align}\label{4.4}
		p(\bar{\bf v})-p(\tilde{v}_{0})
		=\bar{\psi}\left(p(v)-p(\tilde{v}_{0}\right)), \quad \text { where } \quad \bar{\psi}(y)=\inf(\delta_{1}, \sup \left(-\delta_{1}, y\right)).
	\end{align}
$\bar{\mathbf{v}}$ is well-defined since the function $p$ is one to one. Let $\mathbf{\Omega}:=\left\{y\in\R \,|\, p(v)-p(\tilde{v}_{0}) \leq \delta_{1}\right\}$, and
	\begin{align}
		&\mathbf{G}_{1}^{-}(U):=\int_{\mathbf{\Omega}^c}|\mathbf{w}|\cdot |h-\tilde{h}_{0}|^{2} d y,\quad \mathbf{D}(v):=\int_{\mathbb{R}} v^{\beta} \left|\partial_{y} \big(p(v)-p(\tilde{v}_{0}) \big)\right|^{2} d y,\label{4.6}\\
		&\mathbf{G}_{2}(v):=\int_{\mathbb{R}}|\mathbf{w}| Q(v|\tilde{v}_{0}) d y,\quad \tilde{\mathbf{G}}_{2}(v):=\int_{\mathbb{R}} |(\tilde{v}_{0})_{y}| \, \min\{Q(v | \tilde{v}_{0}), \, p(v|\tilde{v}_0)\} d y.\label{4.7}
	\end{align}
	
	Let $q=\f{2\gamma}{\gamma+\alpha}$ with $0<\alpha\leq \gamma\leq \alpha+1$. Then there exist positive constants $\delta_{0}, \delta_{1}$, and $C$ $\mathrm{(}$depending on $\delta_{1}$$\mathrm{)}$ such that for any $\varepsilon_1, \lambda>0$ satisfying $\dis\f{\varepsilon_1}{\lambda}<\delta_{0}^4$, $\v_1\leq \lambda^{\f{2}{2-q}}$ and $\lambda<\delta_{0}$, the following estimates hold
	\begin{align}
		&\int_{\bf\Omega}|\mathbf{w}|\cdot|p(v)-p(\bar{\mathbf{v}})|^{2} dy
		+\int_{\bf\Omega}|\mathbf{w}|\cdot|p(v)-p(\bar{\mathbf{v}})| d y
		\leq C\sqrt{\frac{\v_1}{\lambda}} \mathbf{D}(v), \label{4.9}\\
		&\int_{\bf\Omega}|\mathbf{w}|\cdot\big||p(v)-p(\tilde{v}_{0})|^{2}-|p(\bar{\mathbf{v}})-p(\tilde{v}_{0})|^{2}\big|d y \leq C\sqrt{\frac{\v_1}{\lambda}} \mathbf{D}(v),\label{4.10}
	\end{align}
	and
	\begin{align}
		&\int_{\mathbb{R}}|\mathbf{w}|^{2}\cdot v^{\beta}\cdot |p(v)-p(\bar{\mathbf{v}})|^{2} d y+\int_{\mathbb{R}}|\mathbf{w}|^{2}\cdot v^{\beta}\cdot |p(v)-p(\bar{\mathbf{v}})| dy,\nonumber\\
	&\qquad\qquad\qquad\leq C \lambda^{2}\Big(\mathbf{D}(v)+\delta_{0} \tilde{\mathbf{G}}_{2}(v)\Big)
		+ C\lambda \v_1 \, \mathbf{G}_2(v),\label{4.11}\\[2mm]
		&\int_{\mathbb{R}}|\mathbf{w}|^{2}\cdot \big|v^{\beta}\cdot| p(v)-p(\tilde{v}_{0})|^{2}-\bar{\mathbf{v}}^{\beta}\cdot |p(\bar{\mathbf{v}})-p(\tilde{v}_{0})|^{2}\big| d y \nonumber\\
		&\qquad\qquad\qquad\leq C \lambda^{2}\Big(\mathbf{D}(v)+\delta_{0} \tilde{\mathbf{G}}_{2}(v)\Big)
		+ C\lambda \v_1 \, \mathbf{G}_2(v),\label{4.12}
	\end{align}
	and
	\begin{align}
		&\int_{\mathbb{R}}|\mathbf{w}|\cdot |Q(v |\tilde{v}_{0})-Q(\bar{\mathbf{v}}|\tilde{v}_{0})|dy
    +\int_{\mathbb{R}}|\mathbf{w}|\cdot |v-\bar{\mathbf{v}}| dy
		\leq C\Big(\mathbf{G}_{2}(v)-{\mathbf{G}}_{2}(\bar{\mathbf{v}})\Big),\label{4.14}\\
		&\int_{\mathbb{R}}|\mathbf{w}|\cdot|p(v | \tilde{v}_{0})-p(\bar{\mathbf{v}} | \tilde{v}_{0})| d y  \leq C \sqrt{\frac{\v_1}{\lambda}}\Big[\mathbf{D}(v)+\delta_{0}^2 \tilde{\mathbf{G}}_{2}(v)\Big]+C\Big(\mathbf{G}_{2}(v)-{\mathbf{G}}_{2}(\bar{\mathbf{v}})\Big), \label{4.13}
	\end{align}
	and
	\begin{align}
		&\int_{\mathbf{\Omega}^{c}}|\mathbf{w}|\cdot |p(v)-p(\bar{\mathbf{v}})|^{2} d y
		\leq C \sqrt{\f{\v_1}{\lambda}} \left(\mathbf{D}(v) +\delta_{0}^2 \tilde{\mathbf{G}}_{2}(v)\right)+ C \sqrt{\f{\v_1^{2-q}}{\lambda}}\left(\mathbf{D}(v) +\delta_{0}^2 \tilde{\mathbf{G}}_{2}(v)\right)^q,  \label{4.15}\\
		&\int_{\mathbf{\Omega}^{c}}|\mathbf{w}|\cdot |p(v)-p(\tilde{v}_{0})|\cdot |h-\tilde{h}_{0}| d y \leq \delta_{0}\left(\mathbf{D}(v)+\delta_{0}^2 \tilde{\mathbf{G}}_{2}(v)\right)+\left(\delta_{1}+C \delta_{0}\right) \mathbf{G}_{1}^{-}(U),\label{4.16}\\
		&\int_{\mathbf{\Omega}^{c}}|\mathbf{w}|\Big[Q(\bar{\mathbf{v}} | \tilde{v}_{0})+|\bar{\mathbf{v}}-\tilde{v}_{0}|\Big] d y
		\leq C \sqrt{\frac{\v_1}{\lambda}}\Big(\mathbf{D}(v)+\delta_{0} \tilde{\mathbf{G}}_{2}(v)\Big),\label{4.17}
	\end{align}
	and
	\begin{align}
		&\int_{\mathbb{R}}|\mathbf{w}|^{2}\cdot  \frac{|v^{\beta}-\bar{\mathbf{v}}^{\beta}|^{2}}{v^{\beta}} dy
		\leq C \lambda\left(\mathbf{D}(v)+\delta_{0} \tilde{\mathbf{G}}_{2}(v)\right)+C \lambda\Big(\mathbf{G}_{2}(v)-\mathbf{G}_{2}(\bar{v})\Big),\label{4.18}\\
		&\int_{\mathbb{R}}|\mathbf{w}|^{2}\cdot \Big|\frac{|v^{\beta}-\tilde{v}_{0}^{\beta}|^{2}}{v^{\beta}}-\frac{|\bar{\mathbf{v}}^{\beta}-\tilde{v}_{0}^{\beta}|^{2}}{\bar{\mathbf{v}}^{\beta}}\Big| d y\nonumber \\
		&\qquad \leq C \lambda\Big(\mathbf{D}(v)+\delta_{0} \tilde{\mathbf{G}}_{2}(v)\Big)
		+C \lambda\Big(\mathbf{G}_{2}(v)-\mathbf{G}_{2}(\bar{\mathbf{v}})\Big).\label{4.19}
	\end{align}
\end{lemma}

\begin{lemma}\label{lem4.2}
	Let $q=\f{2\gamma}{\gamma+\alpha}$ with $0<\alpha\leq \gamma\leq \alpha+1$.
	Assume the same hypotheses as in Lemma \ref{lem4.1}. In addition, let $\mathbf{v}_i:\, \R\to \R$ be any functions such that
	\begin{align}\label{4.94}
		|\mathbf{v}_i(y)|\leq \f{\v_1}{\lambda} |\mathbf{w}(y)|,\quad \mbox{for}\quad y\in\R,\, \, i=1,2.
	\end{align}
	Consider the following functionals:
	\begin{align}
		\begin{split}
			&\mathbf{B}_{1}(v) :=\sigma_1 \int_{\mathbb{R}} \mathbf{v}_{1} \cdot p\left(v | \tilde{v}_{0}\right) d y,\quad \mathbf{B}_{2}^{+}(v):=\frac{1}{2 \sigma_1} \int_{\mathbf{\Omega}} \mathbf{w}\cdot \left|p(v)-p\left(\tilde{v}_{0}\right)\right|^{2} d y, \\
			&\mathbf{B}_{2}^{-}(U) :=\int_{\mathbf{\Omega}^{c}} \mathbf{w}\cdot \big(p(v)-p\left(\tilde{v}_{0}\right)\big)\cdot \big(h-\tilde{h}_{0}\big) d y, \\
			&\mathbf{B}_{3}(v) :=-\int_{\mathbb{R}} \mathbf{w}\cdot  v^{\beta}\cdot \left(p(v)-p\left(\tilde{v}_{0}\right)\right)\cdot  \partial_{y}\left(p(v)-p\left(\tilde{v}_{0}\right)\right) d y, \\
			&\mathbf{B}_{4}(v) :=\int_{\mathbb{R}}|\mathbf{w}\cdot |\left|\mathbf{v}_1\right|\cdot \left|p(v)-p\left(\tilde{v}_{0}\right)\right|\cdot |v^{\beta}-\tilde{v}_{0}^{\beta}| d y, \\
			&\mathbf{B}_{5}(v) :=\int_{\mathbb{R}}\left|\mathbf{v}_1\right|\cdot \left|\partial_{y}\left(p(v)-p\left(\tilde{v}_{0}\right)\right)\right|\cdot |v^{\beta}-\tilde{v}_{0}^{\beta}| d y,\nonumber
		\end{split}
	\end{align}
	and
	\begin{align}
		\begin{split}
			\mathbf{Y}^{g}(v) &:=-\frac{1}{2 \sigma^{2}_1(1-\v)^2} \int_{\mathbf{\Omega}} \mathbf{w}\cdot \left|p(v)-p\left(\tilde{v}_{0}\right)\right|^{2} d y-\int_{\mathbf{\Omega}} \mathbf{w}\cdot  Q\left(v|\tilde{v}_{0}\right) d y \\
			&\qquad -\int_{\mathbf{\Omega}} \mathbf{v}_1\cdot  p'(\tilde{v}_0) \cdot \left(v-\tilde{v}_{0}\right) d y+\frac{1}{\sigma_1(1-\v)} \int_{\mathbf{\Omega}} \mathbf{v}_{2}\cdot \left(p(v)-p\left(\tilde{v}_{0}\right)\right) d y, \\
			\mathbf{Y}^{l}(U) &:=\int_{\mathbf{\Omega}} \mathbf{v}_{2}\cdot \left(h-\tilde{h}_{0}-\frac{p(v)-p\left(\tilde{v}_{0}\right)}{\sigma_1(1-\v)}\right) d x, \\
			\mathbf{Y}^{s}(U) &:=-\int_{\mathbf{\Omega}^{c}} \mathbf{w}\cdot  Q\left(v | \tilde{v}_{0}\right) d y-\int_{\mathbf{\Omega}^c} \mathbf{v}_{1}\cdot  \left(v-\tilde{v}_{0}\right) d y\\
			&\qquad-\int_{\mathbf{\Omega}^{c}}\f12 \mathbf{w}\cdot  \left|h-\tilde{h}_{0}\right|^{2} d y +\int_{\mathbf{\Omega}^c} \mathbf{v}_{2}\cdot (h-\tilde{h}_{0}) d y \\
			&=: \mathbf{Y}^{s}_1(v)+\mathbf{Y}^{s}_2(v)+\mathbf{Y}^{s}_3(U)+\mathbf{Y}^{s}_4(U),\\[2mm]
			\mathbf{Y}^b(U)&:=-\frac{1}{2} \int_{\mathbf{\Omega}} \mathbf{w}\cdot \left(h-\tilde{h}_{0}-\frac{p(v)-p\left(\tilde{v}_{0}\right)}{\sigma_1(1-\v)}\right)^{2} d y\\
			&\quad-\frac{1}{\sigma_1(1-\v)} \int_{\mathbf{\Omega}} \mathbf{w}\cdot \left(p(v)-p\left(\tilde{v}_{0}\right)\right)\cdot \left(h-\tilde{h}_{0}-\frac{p(v)-p\left(\tilde{v}_{0}\right)}{\sigma_1(1-\v)}\right) d y,\nonumber
		\end{split}
	\end{align}
	and
	\begin{align}
		\mathbf{G}_{1}^{+}(U):=\int_{\mathbf{\Omega}}|\mathbf{w}|\cdot \left|h-\tilde{h}_{0}-\frac{p(v)-p\left(\tilde{v}_{0}\right)}{\sigma_1(1-\v)}\right|^{2} d y,\nonumber
	\end{align}
	where $\tilde{v}_0$ is any given functions such that $|\tilde{v}_0-v_-|\leq \f12v_-$. We denote
	
	\smallskip
	
	Then there exist positive constants $\delta_{0}, \delta_{1}$, $C^\ast$ and $C$ (in particular, $C^\ast$ depends  on $U_-$, $C$ depends on $U_-$ and $\delta_{1}$, both $C^\ast$ and $C$ are independent of $\delta_0$) such that for any $\varepsilon_1, \lambda>0$ satisfying $\dis\f{\varepsilon_1}{\lambda}<\delta_{0}^4$, $\v_1\leq \lambda^{\f{2}{2-q}}$ and $\lambda<\delta_{0}$, the following estimates hold
	\begin{align}
			&\left|\mathbf{B}_{1}(v)-\mathbf{B}_{1}(\bar{\mathbf{v}})\right| \leq C \delta_{0}\left(\mathbf{D}(v)+\tilde{\mathbf{G}}_{2}(v)+\left(\mathbf{G}_{2}(v)-\mathbf{G}_{2}(\bar{\mathbf{v}})\right)\right),\label{4.99} \\
			&\left|\mathbf{B}_{2}^{-}(U)\right| \leq \delta_{0}\left(\mathbf{D}(v)+\tilde{\mathbf{G}}_{2}(v)\right)+\left(\delta_{1}+C \delta_{0}\right) \mathbf{G}_{1}^{-}(U), \label{4.99-1}\\
			&\left|\mathbf{B}_{2}^{+}(v)-\mathbf{B}_{2}^{+}(\bar{\mathbf{v}})\right| \leq C\delta_{0} \mathbf{D}(v),\label{4.99-2}\\
		&\left|\mathbf{B}_{1}(\bar{\mathbf{v}})\right|+\left|\mathbf{B}_{2}^{+}(\bar{\mathbf{v}})\right| \leq C(v_-) \int_{\mathbb{R}}|\mathbf{w}| Q\left(\bar{\mathbf{v}} | \tilde{v}_{0}\right) d y \leq C^{*} \frac{\v_1^{2}}{\lambda},\label{4.100}\\[2mm]
		&\left|\mathbf{B}_{3}(v)\right|+\left|\mathbf{B}_{4}(v)\right|+\left|\mathbf{B}_{5}(v)\right|\nonumber\\
		& \leq C \delta_{0}\left(\mathbf{D}(v)+\delta_0\tilde{\mathbf{G}}_{2}(v)+\left(\mathbf{G}_{2}(v)-\mathbf{G}_{2}(\bar{\mathbf{v}})\right)+\frac{\v_1}{\lambda} \mathbf{G}_{2}(\bar{\mathbf{v}})\right). \label{4.101}
	\end{align}
	If, in addition, we assume that
	\begin{align}\label{4.95}
		\tilde{\mathbf{G}}_2(v)
		+\mathbf{D}(v) \leq \frac{3 C^{\ast}}{\sqrt{\delta_{0}}} \frac{\v_1^{2}}{\lambda},
	\end{align}
	then it holds that
	\begin{align}\label{4.102}
		&\left|\mathbf{Y}^{g}(v)-\mathbf{Y}^{g}(\bar{\mathbf{v}})\right|^{2}+|\mathbf{Y}^{b}(U)|^{2}+|\mathbf{Y}^{l}(U)|^{2}+\left|\mathbf{Y}^{s}(U)\right|^{2}\nonumber\\
        &\leq C \frac{\v_1^{2}}{\lambda}\bigg\{\delta_{0} \mathbf{D}(v)+\left(\mathbf{G}_{2}(v)-\mathbf{G}_{2}(\bar{\mathbf{v}})\right)+\delta_{0} \sqrt{\f{\v_1}{\lambda}} \tilde{\mathbf{G}}_{2}(v)\nonumber\\
        &+\left(\frac{\v_1}{\lambda}\right)^{\f14} \mathbf{G}_{2}(\bar{\mathbf{v}})+\mathbf{G}_{1}^{-}(U)+\left(\frac{\lambda}{\v_1}\right)^{\f14} \mathbf{G}_{1}^{+}(U)\bigg\}.
	\end{align}
\end{lemma}

\noindent{\bf Proof.}  Since the proof is complicate, we divide it into two steps.

\smallskip

\noindent{\it Step 1. Proofs of \eqref{4.99}-\eqref{4.101}.} Using \eqref{4.13}, \eqref{4.16} and \eqref{4.94}, one has
\begin{align}
	\begin{split}
		&\left|\mathbf{B}_{1}(v)-\mathbf{B}_{1}(\bar{\mathbf{v}})\right|
		\leq C \frac{\v_1}{\lambda}\left(\mathbf{D}(v)+\tilde{\mathbf{G}}_{2}(v)+\left(\mathbf{G}_{2}(v)-\mathbf{G}_{2}(\bar{\mathbf{v}})\right)\right), \nonumber\\
		&\left|\mathbf{B}_{2}^{-}(U)\right| \leq \delta_{0}\left(\mathbf{D}(v)+\tilde{\mathbf{G}}_{2}(v)\right)+\left(\delta_{1}+C \delta_{0}\right) \mathbf{G}_{1}^{-}(U).\nonumber
	\end{split}
\end{align}
It follows from \eqref{4.10} that
\begin{align}
	\left|\mathbf{B}_{2}^{+}(v)-\mathbf{B}_{2}^{+}(\bar{\mathbf{v}})\right|&=\frac{1}{2\sigma_{1}}\int_{\mathbf{\Omega}}|\mathbf{w}|\cdot \Big|| p(v)-\left.p\left(\tilde{v}_{0}\right)\right|^{2}-\left|p(\overline{\mathbf{v}})-p\left(\tilde{v}_{0}\right)\right|^{2} \Big| d y\leq C\sqrt{\frac{\v_1}{\lambda}} \mathbf{D}(v).\nonumber
\end{align}
Thus, for any $\v_1, \lambda$ satisfying $\v_1 / \lambda<\delta_{0}^4$, we obtain the desired estimates \eqref{4.99}-\eqref{4.99-2}.

\smallskip

For \eqref{4.100},  we have from \eqref{4.3-1} that
\begin{align}
	\left|\mathbf{B}_{1}(\bar{\mathbf{v}})\right|+\left|\mathbf{B}_{2}^{+}(\bar{\mathbf{v}})\right| &\leq C(v_-) \int_{\mathbb{R}}|\mathbf{w}| \cdot Q\left(\bar{\mathbf{v}} | \tilde{v}_{0}\right) d y\leq C(v_-) \int_{\mathbb{R}}|\mathbf{w}|\cdot  Q\left(v|\tilde{v}_{0}\right) dy \leq C^{\ast} \frac{\v_1^{2}}{\lambda}.\nonumber
\end{align}

\smallskip

For $\mathbf{B}_3(v)$, using Young's inequality, then one obtains
\begin{align}\label{4.106}
	\left|\mathbf{B}_{3}(v)\right|
	&\leq \delta_{0} \mathbf{D}(v)+\frac{C}{\delta_{0}} \underbrace{\int_{\mathbb{R}}|\mathbf{w}|^{2}\cdot  v^{\beta}\cdot \left|p(v)-p\left(\tilde{v}_{0}\right)\right|^{2} d y}_{=: J_{1}(v)}\nonumber\\
	&\leq \delta_{0} \mathbf{D}(v) + \frac{C}{\delta_{0}}J_{1}(\bar{\mathbf{v}}) + \frac{C}{\delta_{0}}|J_{1}(v)-J_{1}(\bar{\mathbf{v}})|.
\end{align}
For $\mathbf{B}_4(v)$ and $\mathbf{B}_5(v)$, we have
\begin{align}\label{4.107}
	|\mathbf{B}_4(v)|
	&\leq C J_1(v) + \underbrace{\int_{\mathbb{R}}|\mathbf{v}_1|^{2}\cdot  \frac{|v^{\beta}-\tilde{v}_0^{\beta}|^{2}}{v^{\beta}} d y}_{=: J_{2}(v)}\nonumber\\
	&\leq C (J_1(\bar{\mathbf{v}})+ J_2(\bar{\mathbf{v}})) + C\big(|J_1(v)-J_1(\bar{\mathbf{v}})| + |J_2(v)-J_2(\bar{\mathbf{v}})| \big),
\end{align}
and
\begin{align}\label{4.108}
	|\mathbf{B}_5(v)|&\leq \delta_{0} \mathbf{D}(v) + \frac{C}{\delta_{0}}J_2(v)\leq \delta_{0} \mathbf{D}(v) + \frac{C}{\delta_{0}}J_2(\bar{\mathbf{v}})
	+ \frac{C}{\delta_{0}} |J_2(v)-J_2(\bar{\mathbf{v}})|.
\end{align}
Using \eqref{4.12}, one obtains
\begin{align}\label{4.109}
	|J_1(v)-J_1(\bar{\mathbf{v}})|\leq C \lambda^{2}\Big(\mathbf{D}(v)+\delta_{0} \tilde{\bf G}_{2}(v)\Big)
	+ C\lambda \v_1 \, \mathbf{G}_2(v).
\end{align}
Using \eqref{4.94} and \eqref{4.19}, one has
\begin{align}\label{4.110}
	|J_2(v)-J_2(\bar{\mathbf{v}})|
	&\leq \frac{\v_1^2}{\lambda^2} \int_{\mathbb{R}}|\mathbf{w}|^{2}\cdot \Big|\frac{|v^{\beta}-\tilde{v}_{0}^{\beta}|^{2}}{v^{\beta}}-\frac{|\bar{\mathbf{v}}^{\beta}-\tilde{v}_{0}^{\beta}|^{2}}{\bar{\mathbf{v}}^{\beta}}\Big| d y\nonumber\\
	&\leq   C\frac{\v_1^2}{\lambda} \left\{ \mathbf{D}(v)+\delta_{0} \tilde{\mathbf{G}}_{2}(v)
	+ \Big(\mathbf{G}_{2}(v)-\mathbf{G}_{2}(\bar{\mathbf{v}})\Big) \right\}.
\end{align}
For the remaining terms $J_1(\bar{\mathbf{v}})$ and $J_2(\bar{\mathbf{v}})$, it follows from \eqref{4.94} and \eqref{4.2}  that
\begin{align}\label{4.111}
	\frac{C}{\delta_{0}}[J_1(\bar{\mathbf{v}})+J_2(\bar{\mathbf{v}})]
	&\leq \frac{C}{\delta_{0}} \Big[ \v_1 \lambda \int_{\mathbb{R}}|\mathbf{w}|\cdot  \bar{\mathbf{v}}^{\beta}\cdot \left|p(\bar{\mathbf{v}})-p\left(\tilde{v}_{0}\right)\right|^{2} d y +
	\frac{\v_1^3}{\lambda} \int_{\mathbb{R}}|\mathbf{w}|\cdot  \frac{|\bar{\mathbf{v}}^{\beta}-\tilde{v}_{0}^{\beta}|^{2}}{\bar{\mathbf{v}}^{\beta}}  d y\Big]\nonumber\\
	&\leq C\delta_0 \f{\v_1}{\lambda}\Big[\int_{\mathbb{R}}|\mathbf{w}|\cdot  \bar{\mathbf{v}}^{\beta}\cdot \left|p(\bar{\mathbf{v}})-p\left(\tilde{v}_{0}\right)\right|^{2} d y +
	\int_{\mathbb{R}}|\mathbf{w}|\cdot  \frac{|\bar{\mathbf{v}}^{\beta}-\tilde{v}_{0}^{\beta}|^{2}}{\bar{\mathbf{v}}^{\beta}}  d y\Big].
\end{align}
Noting $C^{-1}\leq \bar{\mathbf{v}}\leq C$ and $\left|p(\bar{\mathbf{v}})-p\left(\tilde{v}_{0}\right)\right|^{2}+|\bar{\mathbf{v}}^{\beta}-\tilde{v}_{0}^{\beta}|^{2}\leq C |\bar{\mathbf{v}}-\tilde{v}_{0}|^{2}\leq CQ(\bar{\mathbf{v}}|\tilde{v}_0)$, we get
\begin{align}
	&\int_{\mathbb{R}}|\mathbf{w}|\cdot  \bar{\mathbf{v}}^{\beta}\cdot \left|p(\bar{\mathbf{v}})-p\left(\tilde{v}_{0}\right)\right|^{2} d y + \int_{\mathbb{R}}|\mathbf{w}|\cdot  \frac{|\bar{\mathbf{v}}^{\beta}-\tilde{v}_{0}^{\beta}|^{2}}{\bar{\mathbf{v}}^{\beta}}  d y\nonumber\\
	&\leq C(v_-)\int_{\mathbb{R}} |\mathbf{w}| Q(\bar{\mathbf{v}}|\tilde{v}_0) dy\leq C(v_-)\mathbf{G}_2(\bar{\mathbf{v}}),\nonumber
\end{align}
which, together with \eqref{4.106}-\eqref{4.111}, yields that
\begin{align}
	&|\mathbf{B}_{3}(v)| +|\mathbf{B}_{4}(v)| + |\mathbf{B}_{5}(v)| \leq C \delta_{0} \Big[\mathbf{D}(v) +  \delta_0\tilde{\mathbf{G}}_{2}(v)
	+ \big(\mathbf{G}_{2}(v)-\mathbf{G}_{2}(\bar{\mathbf{v}})\big)+  \f{\v_1}{\lambda}\mathbf{G}_2(\bar{\mathbf{v}})\Big].\nonumber
\end{align}

\medskip

\noindent{\it Step 2. Proof of \eqref{4.102}.} Using \eqref{4.9}-\eqref{4.10}, \eqref{4.14} and \eqref{4.94}, one obtains
\begin{align}\label{4.114}
	&|\mathbf{Y}^g(v)-\mathbf{Y}^g(\bar{\mathbf{v}})|+ |\mathbf{Y}_1^s(v)-\mathbf{Y}_1^s(\bar{\mathbf{v}})| +|\mathbf{Y}_2^s(v)-\mathbf{Y}_2^s(\bar{\mathbf{v}})|\nonumber\\
	&\leq C\int_{\mathbf{\Omega}} |\mathbf{w}|\cdot \big||p(v)-p(\tilde{v}_{0})|^{2}-|p(\bar{\mathbf{v}})-p(\tilde{v}_{0})|^{2}\big| d y +C \int_{\mathbf{\Omega}} |\mathbf{w}|\cdot  |p(v)-p(\bar{\mathbf{v}})| d y\nonumber\\
	&\quad +\int_{\R} |\mathbf{w}| \cdot \Big(| Q\left(v|\tilde{v}_{0}\right)- Q\left(\bar{\mathbf{v}}|\tilde{v}_{0}\right)|+ |v-\bar{\mathbf{v}}| \Big) d y \nonumber\\
	&\leq C\sqrt{\f{\v_1}{\lambda}}\mathbf{D}(v) + C \Big(\mathbf{G}_{2}(v)-\mathbf{G}_{2}(\bar{\mathbf{v}})\Big).
\end{align}
It follows from \eqref{4.17} and \eqref{4.94} that
\begin{align}\label{4.115}
	\sum\limits_{i=1}^2|\mathbf{Y}_i^s(\bar{\mathbf{v}})|
	&\leq \int_{\mathbf{\Omega}^{c}} |\mathbf{w}|\cdot  Q(\bar{\mathbf{v}} | \tilde{v}_{0}) d y+\int_{\mathbf{\Omega}^c} |\mathbf{v}_{1}|\cdot  |\bar{\mathbf{v}}-\tilde{v}_0| d y\leq C\sqrt{\f{\v_1}{\lambda}} \Big\{\mathbf{D}(U) + \delta_0 \tilde{\mathbf{G}}_2(U)\Big\}.
\end{align}
By using \eqref{4.99}-\eqref{4.100}, one can obtain
\begin{align}\label{4.116}
	|\mathbf{Y}_3^s(U)| + |\mathbf{Y}^b(U)|&\leq C\int_{\mathbf{\Omega}^{c}} |\mathbf{w}|\cdot |h-\tilde{h}_{0}|^{2} d y
	+ C  \int_{\mathbf{\Omega}} |\mathbf{w}|\cdot \left( |h-\tilde{h}_{0}|^2 + |p(v)-p(\tilde{v}_{0}|^2\right)  d y\nonumber\\
	&\leq C\int_{\R} |\mathbf{w}| \cdot |h-\tilde{h}_{0}|^{2} d y
	+ C   \int_{\mathbf{\Omega}}  |\mathbf{w}|\cdot  |p(v)-p(\tilde{v}_{0}|^2   d y  \nonumber\\
	&\leq C\int_{\R} |\mathbf{w}| \cdot |h-\tilde{h}_{0}|^{2} d y + C |\mathbf{B}_2^+(v)|\nonumber\\
	&\leq  C\int_{\R} |\mathbf{w}| \cdot |h-\tilde{h}_{0}|^{2} d y + C \delta_{0}  \mathbf{D}(v) + C\f{\v_1^2}{\lambda}.
\end{align}
Thus, under the assumption \eqref{4.95}, and using \eqref{4.7}, \eqref{4.114}-\eqref{4.116} and \eqref{4.3-1}, then one has
\begin{align}\label{4.117}
	&|\mathbf{Y}^g(v)-\mathbf{Y}^g(\bar{\mathbf{v}})|+ |\mathbf{Y}_1^s(v)| +|\mathbf{Y}_2^s(v)| + |\mathbf{Y}_3^s(U)| + |\mathbf{Y}^b(U)|\nonumber\\
	&\leq |\mathbf{Y}^g(v)-\mathbf{Y}^g(\bar{\mathbf{v}})| +  |\mathbf{Y}_1^s(v)-\mathbf{Y}_1^s(\bar{\mathbf{v}})| +  |\mathbf{Y}_2^s(v)-\mathbf{Y}_2^s(\bar{\mathbf{v}})|\nonumber\\
	&\quad + |\mathbf{Y}_1^s(\bar{\mathbf{v}})| + |\mathbf{Y}_2^s(\bar{\mathbf{v}})| + |\mathbf{Y}_3^s(U)| + |\mathbf{Y}^b(U)|\nonumber\\
	&\leq C\left(\sqrt{\f{\v_1}{\lambda}}+ \delta_0\right)\mathbf{D}(v) + C \Big(\mathbf{G}_{2}(v)-\mathbf{G}_{2}(\bar{\mathbf{v}})\Big)
	+ C\sqrt{\f{\v_1}{\lambda}}   \delta_0 \tilde{\mathbf{G}}_2(v) \nonumber\\
	&\quad + C\int_{\R} |\mathbf{w}| \cdot |h-\tilde{h}_{0}|^{2} d y +  C\f{\v_1^2}{\lambda} \leq C\f{\v_1^2}{\lambda}.
\end{align}

\medskip

Using Young's inequality, \eqref{4.99}-\eqref{4.100}, we have
\begin{align*}
	&\left|\int_{\mathbf{\Omega}} \mathbf{w}\cdot \left(p(v)-p\left(\tilde{v}_{0}\right)\right)\cdot \left(h-\tilde{h}_{0}-\frac{p(v)-p\left(\tilde{v}_{0}\right)}{\sigma_1(1-\v)}\right) d y\right|\nonumber \\
	&\leq\left(\frac{\lambda}{\v_1}\right)^{\f14} \mathbf{G}_{1}^{+}(U)+C\left(\frac{\v_1}{\lambda}\right)^{\f14} \int_{\mathbf{\Omega}}|\mathbf{w}|\left|p(v)-p\left(\tilde{v}_{0}\right)\right|^{2} d y \nonumber \\
	& \leq\left(\frac{\lambda}{\v_1}\right)^{\f14} \mathbf{G}_{1}^{+}(U)+C\left(\frac{\v_1}{\lambda}\right)^{\f14}\Big(\mathbf{B}_{2}^{+}(\bar{\mathbf{v}})+\left(\mathbf{B}_{2}^{+}(v)-\mathbf{B}_{2}^{+}(\bar{\mathbf{v}})\right)\Big)\nonumber \\
	& \leq\left(\frac{\lambda}{\v_1}\right)^{\f14} \mathbf{G}_{1}^{+}(U)+C\left(\frac{\v_1}{\lambda}\right)^{\f14}\Big(\mathbf{G}_{2}(\bar{\mathbf{v}})+ \delta_{0} \mathbf{D}(v)\Big),
\end{align*}
which, together with the fact $\frac{\lambda}{\v_1}\gg 1$, yields
\begin{align}\label{4.118}
	\left|\mathbf{Y}^{b}(U)\right| \leq C\left(\frac{\lambda}{\v_1}\right)^{1 / 4} \mathbf{G}_{1}^{+}(U)+C\left(\frac{\v_1}{\lambda}\right)^{1 / 4}\Big(\mathbf{G}_{2}(\bar{\mathbf{v}})+ \delta_{0} \mathbf{D}(v)\Big).
\end{align}
Combining \eqref{4.118}, \eqref{4.114}-\eqref{4.115} and   $|\mathbf{Y}_{3}^{s}(U)| \leq C \mathbf{G}_{1}^{-}(U)$, we get
\begin{align}\label{4.119}
	&\left|\mathbf{Y}^{g}(v)-\mathbf{Y}^{g}(\bar{\mathbf{v}})\right|+\left|\mathbf{Y}_{1}^{s}(v)\right|+\left|\mathbf{Y}_{2}^{s}(v)\right|+\left|\mathbf{Y}_{3}^{s}(U)\right|+\left|\mathbf{Y}^{b}(U)\right| \nonumber\\
	&\leq |\mathbf{Y}^g(v)-\mathbf{Y}^g(\bar{\mathbf{v}})| +  |\mathbf{Y}_1^s(v)-\mathbf{Y}_1^s(\bar{\mathbf{v}})| +  |\mathbf{Y}_2^s(v)-\mathbf{Y}_2^s(\bar{\mathbf{v}})|\nonumber\\
	&\quad + |\mathbf{Y}_1^s(\bar{\mathbf{v}})| + |\mathbf{Y}_2^s(\bar{\mathbf{v}})| + |\mathbf{Y}_3^s(U)| + |\mathbf{Y}^b(U)|\nonumber\\
	&\leq  C\sqrt{\f{\v_1}{\lambda}}\mathbf{D}(v) + C \Big(\mathbf{G}_{2}(v)-\mathbf{G}_{2}(\bar{\mathbf{v}})\Big)
	+ C\sqrt{\f{\v_1}{\lambda}} \Big\{\mathbf{D}(v) + \delta_0 \tilde{\mathbf{G}}_2(v)\Big\}\nonumber\\
	&\quad + C \mathbf{G}_{1}^{-}(U) + C\left(\frac{\lambda}{\v_1}\right)^{\f14} \mathbf{G}_{1}^{+}(U)+C\left(\frac{\v_1}{\lambda}\right)^{\f14}\Big(\mathbf{G}_{2}(\bar{\mathbf{v}})+ \delta_{0} \mathbf{D}(v)\Big)\nonumber\\
	&\leq C\bigg( \delta_{0} \mathbf{D}(v)+\left(\mathbf{G}_{2}(v)-\mathbf{G}_{2}(\bar{\mathbf{v}})\right)+\delta_{0}\sqrt{\f{\v_1}{\lambda}} \tilde{\mathbf{G}}_{2}(v)+\left(\frac{\v_1}{\lambda}\right)^{\f14} \mathbf{G}_{2}(\bar{\mathbf{v}})\nonumber \\
	&\qquad+\mathbf{G}_{1}^{-}(U)+\left(\frac{\lambda}{\v_1}\right)^{\f14} \mathbf{G}_{1}^{+}(U)\bigg) .
\end{align}

For the remaining terms $\mathbf{Y}_{4}^{s}(U)$ and $\mathbf{Y}^{l}(U)$, it follows from the Holder's inequality, \eqref{4.94} and \eqref{4.3} that
\begin{align}\label{4.120}
	\begin{split}
		&\left|\mathbf{Y}_{4}^{s}(U)\right|^{2} \leq C\left(\frac{\v_1}{\lambda}\right)^{2}\left(\int_{\mathbb{R}}|\mathbf{w}| d y\right) \int_{\mathbf{\Omega}^{c}}|\mathbf{w}|\cdot |h-\tilde{h}_{0}|^{2} d y\leq C \frac{\v_1^{2}}{\lambda} \mathbf{G}_{1}^{-}(U), \\
		&\left|\mathbf{Y}^{l}(U)\right|^{2} \leq C\left(\frac{\v_1}{\lambda}\right)^{2}\left(\int_{\mathbb{R}}|\mathbf{w}| d y\right) \int_{\mathbf{\Omega}}|\mathbf{w}|\cdot \left|h-\tilde{h}_{0}-\frac{p(v)-p\left(\tilde{v}_{0}\right)}{\sigma_1(1-\v)}\right|^{2} d y \leq C \frac{\v_1^{2}}{\lambda} \mathbf{G}_{1}^{+}(U).
	\end{split}
\end{align}
Hence, combining \eqref{4.117}, \eqref{4.119} and \eqref{4.120}, we obtain
\begin{align*}
	&|\mathbf{Y}^{g}(v)-\mathbf{Y}^{g}(\bar{\mathbf{v}})|^2+|\mathbf{Y}^{b}(U)|^2+|\mathbf{Y}^{l}(U)|^2+|\mathbf{Y}^{s}(U)|^2 \nonumber\\
	&\leq 2 \Big(|\mathbf{Y}^{g}(v)-\mathbf{Y}^{g}(\bar{\mathbf{v}})|+|\mathbf{Y}^{b}(U)|+|\mathbf{Y}_1^{s}(v)|+|\mathbf{Y}_2^{s}(v)|+|\mathbf{Y}_3^{s}(h)|\Big)^2 \nonumber\\
	&\quad + 2|\mathbf{Y}_{4}^{s}(U)|^{2}+ |\mathbf{Y}^{l}(U)|^2\nonumber\\
	&\leq C\f{\v_1^2}{\lambda} \bigg( \delta_{0} \mathbf{D}(v)+\left(\mathbf{G}_{2}(v)-\mathbf{G}_{2}(\bar{\mathbf{v}})\right)+\delta_{0} \sqrt{\f{\v_1}{\lambda}} \tilde{\mathbf{G}}_{2}(v)+\left(\frac{\v_1}{\lambda}\right)^{\f14} \mathbf{G}_{2}(\bar{\mathbf{v}})\nonumber \\
	&\qquad+\mathbf{G}_{1}^{-}(U)+\left(\frac{\lambda}{\v_1}\right)^{\f14} \mathbf{G}_{1}^{+}(U)\bigg),
\end{align*}
which yields \eqref{4.102}.
Therefore the proof of Lemma \ref{lem4.2} is complete. $\hfill\Box$

\section{Estimates for large perturbation part}
In this section, we shall consider the estimates for large perturbation part of the solution $U=(v,h)$ to \eqref{1.15-1}. The main new difficulty comes from the  interaction of shock and rarefaction wave with the solution in the case of large perturbation. Recall the definitions of $B_1(v), B_{2}^-(U),B_{2}^+(v),B_3(v),  B_4(v),B_5(U)$ in \eqref{3.29} and $G_1^{\pm}(U),G_2(v),  G^r(v),G^h(h),D(v),G_3(v)$ in \eqref{3.30}.

\smallskip

For later use, we define
\begin{align}
	\begin{split}
		&p(\bar{v})-p(\tilde{v}_{0}):=\bar{\psi}_{\delta_1}(p(v)-p(\tilde{v}_{0})) \quad  \mbox{and}\quad \bar{\psi}_{\delta_{1}}(y)=\inf(\delta_{1},\sup(-\delta_{1},y)),\\
		&\bar{U}:=(\bar{v},h),\quad \Omega:=\{y\in\R\, |\, p(v)-p(\tilde{v})\leq \delta_1\}.
	\end{split}\nonumber
\end{align}
Noting \eqref{3.22}, we can rewrite $Y(U)$ as
\begin{align}\label{5.2-1}
	Y(U)
&=-\int_{\R}\f12 \partial_y w_X\, |h-\tilde{h}|^2 \, dy  -  \int_{\R}  \partial_y w_X\,Q(v|\tilde{v}) \,dy\nonumber\\
	&\quad - \int_{\R} w_X\,  p'(\tilde{v}_X^s) \, \partial_y\tilde{v}^s_X\, (v-\tilde{v})\, dy + \int_{\R} w_X\,  \partial_y\tilde{h}^s_X\, (h-\tilde{h})\, dy\nonumber\\
	&\quad - \int_{\R} w_X\, [ p'(\tilde{v}) - p'(\tilde{v}_X^s) ]\, \partial_y\tilde{v}^s_X\, (v-\tilde{v})\, dy\nonumber\\
	&=: \mathbb{Y}^g(v) + \mathbb{Y}^b(U) + \mathbb{Y}^l(U) + \mathbb{Y}^s(U),
\end{align}
where
\begin{align}\label{5.2-3}
	\mathbb{Y}^g(v):= &-\f1{2\sigma_1^2(1-\v)^2}\int_{\Omega} \partial_y w_X\, |p(v)-p(\tilde{v})|^2 \, dy  -  \int_{\Omega}  \partial_y w_X\,Q(v|\tilde{v}) \,dy\nonumber\\
	&- \int_{\Omega} w_X\,  p'(\tilde{v}_X^s) \, \partial_y\tilde{v}^s_X\, (v-\tilde{v})\, dy + \f{1}{\sigma_1(1-\v)}\int_{\Omega} w_X\,  \partial_y\tilde{h}^s_X\, (p(v)- p(\tilde{v}))\, dy,
\end{align}
\begin{align}\label{5.2-4}
	\mathbb{Y}^b(U):=&-\frac{1}{2} \int_{\Omega} \partial_y w_X\,\left(h-\tilde{h}-\frac{p(v)-p\left(\tilde{v}\right)}{\sigma_1(1-\v)}\right)^{2} d y\nonumber\\
	&-\frac{1}{\sigma_1(1-\v)} \int_{\Omega} \partial_y w_X\,\left(p(v)-p\left(\tilde{v}\right)\right)\left(h-\tilde{h}-\frac{p(v)-p\left(\tilde{v}\right)}{\sigma_1(1-\v)}\right) d y,
\end{align}

\begin{align}\label{5.2-5}
	\mathbb{Y}^l(U)&:=  \int_{\Omega}  w_X\,\partial_y\tilde{h}_X^s\, \left(h-\tilde{h}-\frac{p(v)-p\left(\tilde{v}\right)}{\sigma_1(1-\v)}\right) dy,
\end{align}
and
\begin{align}\label{5.2-2}
	\mathbb{Y}^s(U):=&-\int_{\Omega^c}  \partial_y w_X Q(v|\tilde{v})dy- \int_{\Omega^c} w_X\,  p'(\tilde{v}_X^s) \, \partial_y\tilde{v}^s_X\, (v-\tilde{v})\, dy
	\nonumber\\
	&-\int_{\Omega^c} \f12 \partial_y w_X |h-\tilde{h}|^2 dy+ \int_{\Omega^c} w_X\,  \partial_y\tilde{h}^s_X\, (h-\tilde{h})\, dy\nonumber\\
	&-\int_{\R} w_X\, [ p'(\tilde{v}) - p'(\tilde{v}_X^s) ]\, \partial_y\tilde{v}^s_X\, (v-\tilde{v})\, dy\nonumber\\
	=:& \mathbb{Y}^s_1(v) + \mathbb{Y}^s_2(v) +\mathbb{Y}^s_3(U) + \mathbb{Y}^s_4(U) + \mathbb{Y}^s_5(v).
\end{align}

\smallskip

Motivated by \cite{Kang-Vasseur-2020}, we define the equation of shift function
\begin{align}\label{5.7-1}
	\dot{X}(\tau)=
	\begin{cases}
		\dis	-\f12 \sigma_1,\qquad \qquad \qquad \qquad \quad \mbox{if}\,\, Y(U)\leq -\v_1^2,\\[3mm]
		\dis	\f12 \sigma_1 \cdot \f{Y(U)}{\v_1^2}, \qquad \qquad \qquad\,\,\,\mbox{if}\,\, -\v_1^2\leq Y(U)\leq 0,\\[3mm]
		\dis	-\f{1}{\v_1^4} Y(U)\cdot (2 |J^{bad}|+1),\quad \,\mbox{if}\,\, 0\leq Y(U)\leq \v_1^2,\\[3mm]
		\dis	-\f{1}{\v_1^2} \cdot (2 |J^{bad}|+1),\qquad\quad \,\,\,\,\mbox{if}\,\, Y(U)\geq \v_1^2,
	\end{cases}
\end{align}
with initial data $X(0)=0$.
It follows from \eqref{5.7-1} that
\begin{align*}
	\dot{X}(\tau)\leq
	\begin{cases}
		\dis	-\f12 \sigma_1,\quad \mbox{if}\,\, Y(U)\leq -\v_1^2,\\[3mm]
		\dis	-\f12 \sigma_1, \quad \mbox{if}\,\, -\v_1^2\leq Y(U)\leq 0,\\[3mm]
		\dis	0,\qquad \quad\mbox{if}\,\, 0\leq Y(U)\leq \v_1^2,\\[3mm]
		\dis	0,\qquad \quad\mbox{if}\,\, Y(U)\geq \v_1^2.
	\end{cases}
\end{align*}
Hence we have
\begin{align}
	\dot{X}(\tau)\leq \f12 |\sigma_1|\qquad \text{for all }\tau>0,\nonumber
\end{align}
which yields that
\begin{align}\label{5.7-3}
	X(\tau)\leq \f{1}{2} |\sigma_1| \tau =-\f12 \sigma_1 \tau \qquad \text{for all }\tau>0.
\end{align}

\begin{lemma}\label{lem5.1}
	For $X(\tau)$ defined above, the following interaction estimates hold:
	\begin{align}\label{5.10-1}
		&|\tilde{v}^r-v_m|\cdot |\partial_y \tilde{v}^s_X|
		\leq
		\begin{cases}
			\dis  C(v_-)\v_2 |\partial_y \tilde{v}^s_X|\exp\Big(-\f{\nu}{a}(|y|+\lambda_2(v_m)\tau)\Big),\quad\,\,y\leq 0,\\[2mm]
			\dis C(v_-)\v_1^2 \v_2 \exp\Big(-c_1\v_1(|y|+\f12|\sigma_1|\tau)\Big),\qquad\quad y\geq0,
		\end{cases}
	\end{align}
	\begin{align}\label{5.11-1}
		&|\tilde{v}^s_X-v_m|\cdot |\partial_y \tilde{v}^r|
		\leq
		\begin{cases}
			\dis C(v_-)\f{\nu}{a}\v_1 \v_2 \exp\Big(-\f{\nu}{a}(|y|+\lambda_2(v_m)\tau)\Big),\qquad y\leq 0,\\[2mm]
			\dis  C(v_-)\v_1  |\partial_y \tilde{v}^r| \exp\Big(-c_1\v_1(|y|+\f12|\sigma_1|\tau)\Big),\quad y\geq0,
		\end{cases}
	\end{align}
	and
	\begin{align}\label{5.12-2}
		|\partial_y\tilde{v}_X^s\cdot \partial_y\tilde{v}^r|
		\leq
		\begin{cases}
			\dis C(v_-)  |\partial_y\tilde{v}_X^s|\f{\nu \v_2}{a}  \exp\Big(-\f{\nu}{a}(|y|+\lambda_2(v_m)\tau)\Big),\qquad\,y\leq 0,\\[2mm]
			\dis C(v_-) \v_1^2 \exp\Big(-c_1\v_1 (|y|+\f12|\sigma_1|\tau) \Big) \cdot |\partial_y\tilde{v}^r| ,\qquad\,\,y\geq 0,
		\end{cases}
	\end{align}
	where $C(v_-)$ is a positive constant depending only on $v_-$.
\end{lemma}

\noindent{\bf Proof.} Noting \eqref{5.7-3}, one has
\begin{align*}
	\sigma_1\tau + X(\tau)\leq \f12 \sigma_1\tau \leq -\f12 |\sigma_1| \tau,
\end{align*}
which yields that
\begin{align}\label{5.12-1}
	y-\sigma_1\tau - X(\tau)\geq y+\f12 |\sigma_1| \tau.
\end{align}
It follows from \eqref{5.12-1} and Lemmas \ref{lem2.3}--\ref{lem2.2} that
\begin{align}
	|\tilde{v}^r-v_m|\cdot |\partial_y \tilde{v}^s_X|
	&\leq
	\begin{cases}
		\dis C(v_-) \v_2 |\partial_y \tilde{v}^s_X| \exp\Big(-\f{\nu}{a}|y-\lambda_2(v_m)\tau|\Big),\quad\,\,\,y\leq 0,\\[2mm]
		\dis C(v_-)\v_1^2 \v_2  \exp\Big(-c_1\v_1 |y-\sigma_1\tau-X(\tau)|\Big),\quad y\geq 0,
	\end{cases}\nonumber\\
	&\leq
	\begin{cases}
		\dis C(v_-)\v_2 |\partial_y \tilde{v}^s_X| \exp\Big(-\f{\nu}{a}(|y|+\lambda_2(v_m)\tau)\Big),\quad y\leq 0,\\[2mm]
		\dis C(v_-)\v_1^2 \v_2 \exp\Big(-c_1\v_1(|y|+\f12|\sigma_1|\tau)\Big),
		\qquad\,\,\,\, y\geq0,
	\end{cases}\nonumber
\end{align}
which yields \eqref{5.10-1}.

\smallskip

For \eqref{5.11-1}, using \eqref{5.12-1} and Lemmas \ref{lem2.3}--\ref{lem2.2}, one obtains
\begin{align}
	|\tilde{v}^s_X-v_m|\cdot |\partial_y \tilde{v}^r|
	&\leq
	\begin{cases}
		\dis C(v_-)\f{\nu}{a}\v_1 \v_2 \exp\Big(-\f{\nu}{a}(|y|+\lambda_2(v_m)\tau)\Big),\qquad \quad\,\,y\leq 0,\\[2mm]
		\dis  C(v_-)\v_1  |\partial_y \tilde{v}^r| \exp\Big(-c_1\v_1 |y-\sigma_1\tau-X(\tau)|\Big),\quad y\geq0,
	\end{cases}\nonumber\\
	&\leq
	\begin{cases}
		\dis C(v_-)\f{\nu}{a}\v_1 \v_2 \exp\Big(-\f{\nu}{a}(|y|+\lambda_2(v_m)\tau)\Big),\qquad \quad\,\,\,y\leq 0,\\[2mm]
		\dis  C(v_-)\v_1  |\partial_y \tilde{v}^r| \exp\Big(-c_1\v_1(|y|+\f12|\sigma_1|\tau)\Big),\qquad \,\,\,\,y\geq0.\nonumber
	\end{cases}
\end{align}

Finally, using \eqref{2.3-3} and \eqref{2.14}, one can obtain \eqref{5.12-2} directly. Therefore the proof of Lemma \ref{lem5.1} is complete.  $\hfill\Box$

\begin{lemma}\label{lem5.2}
	Let $q=\f{2\gamma}{\gamma+\alpha}$ with $0<\alpha\leq \gamma\leq \alpha+1$. let $U=(v,h)$ be the solution of Navier-Stokes equations \eqref{1.15-1}. Then there exist positive constants $\delta_{0}, \delta_{1}$, $C_\ast,\,C_1^\ast$ and $C$ $\mathrm{(}$in particular, $C_\ast,\, C_1^\ast$ depend  only on $U_-$; $C$ depends on $U_-$ and $\delta_{1}$;  all $C_\ast, C_1^\ast$ and $C$ are independent of $\delta_0$$\mathrm{)}$ such that for any $\varepsilon_1, \lambda>0$ satisfying $\dis\f{\varepsilon_1}{\lambda}<\delta_{0}^4$, $\v_1\leq \lambda^{\f{2}{2-q}}$, $\lambda<\delta_{0}$ and $\f{\nu}{a}\leq \delta_0$, the following estimates hold:
	
	\smallskip
	
	\noindent{\bf Part I.} For $U$ satisfying $Y(U)\leq \v_1^2$, one obtains
	\begin{align}
         &\begin{aligned}
			&|B_{1}(v)-B_{1}(\bar{v})| \leq C  \delta_{0}\Big[D(v)+[G_{2}(v)-G_{2}(\bar{v})] +\f{\v_1}{\lambda}G_{2}(v) + G^r(v)\Big], \\
			&|B_{2}^{-}(U)| \leq C\delta_{0} \Big[D(v)+\f{\v_1}{\lambda}G_{2}(v) + G^r(v) \Big]+\left(\delta_{1}+C \delta_{0}\right) G_{1}^{-}(U), \\
			&|B_{2}^{+}(v)-B_{2}^{+}(\bar{v})| \leq  C\delta_{0}  D(v),\\
	&|B_{1}(\bar{v})|+|B_{2}^{+}(\bar{v})| \leq C(v_{-}) \int_{\mathbb{R}}|\partial_yw_X|\cdot Q\left(\bar{{v}} | \tilde{v}_{0}\right) d y \leq C_1^{*} \frac{\v_1^{2}}{\lambda},
	 \end{aligned}\label{5.1}\\
	&|B_{3}(v)|+|B_{4}(v)|
		\leq C \delta_{0}\left[ D(v)+G^{r}(v)+\left(G_{2}(v)-G_{2}(\bar{v})\right)+\frac{\v_1}{\lambda} G_{2}(\bar{v})\right],\label{5.3}\\
	&\begin{aligned}
		|B_5(U)| &\leq C\delta_0 [D(v) + G^r(v)]+ C(v_-)\delta_0 \f{\v_1}{\lambda} G_2(v) +\delta_0 \v_1\v_2\, G^h(h)\\
		&\quad +  f(\nu,a,\v_1,\v_2,\tau) \int_{\R} \eta(U|\tilde{U}) dy  +g(\nu,a,\v_1,\v_2,\tau),\label{5.10}
	\end{aligned}\\
&\begin{aligned}
|B_{\delta_1}(U)|&\leq C \delta_0   \Big[D(v) + G^r(v)\Big]  +  f(\nu,a,\v_1,\v_2,\tau) \int_{\R} \eta(U|\tilde{U}) dy\\
&\quad+4 C_1^\ast \f{\v_1^2}{\lambda} + g(\nu,a,\v_1,\v_2,\tau),\label{5.11}
\end{aligned}
\end{align}
where
\begin{align}\label{5.12}
		\begin{split}
			f(\nu,a,\v_1,\v_2,\tau)&:=C(v_-)\Big[ \f{1}{\delta_0} (\min\{\f{\nu\v_2}{a},\tau^{-1}\})^{\f43}
			+ \f{\nu}{a} \exp\{-\f{\nu}{C_\ast a}\tau\} + \v_1 \exp\{-\f{\v_1}{C_\ast}\tau\} \Big],\\
			g(\nu,a,\v_1,\v_2,\tau)&:=C(v_{-})\Big[\f{1}{\delta_0} (\min\{\f{\nu\v_2}{a},\tau^{-1}\})^{\f43}
			+ \f{\v_1\v_2}{\delta_0\lambda}\Big(\exp\{-\f{\nu}{C_\ast a}\tau\} +   \exp\{-\f{\v_1}{C_\ast}\tau\}\Big) \Big].
		\end{split}
\end{align}
	
\smallskip
	
\noindent{\bf Part II.} For $U$ satisfying $Y(U)\leq \v_1^2$, $\dis G^r(v)+  D(v)\leq 2\f{C_1^\ast}{\sqrt{\delta_0}} \f{\v_1^2}{\lambda}$, one gets
	\begin{align}\label{5.13}
	&|\mathbb{Y}^{g}(v)-\mathbb{Y}^{g}(\bar{v})|^2+|\mathbb{Y}^{b}(U)|^2+|\mathbb{Y}^{l}(U)|^2+|\mathbb{Y}^{s}(U)|^2 \nonumber\\
		&\leq C\f{\v_1^2}{\lambda}  \bigg\{ \delta_{0} \Big[D(v) + G^r(v)\Big]+\Big[G_{2}(v)-G_{2}(\bar{v})\Big] +\left(\frac{\v_1}{\lambda}\right)^{\f14} G_{2}(\bar{v})\nonumber \\
		&\quad+G_{1}^{-}(U)+\left(\frac{\lambda}{\v_1}\right)^{\f14} G_{1}^{+}(U)\bigg\} + C\v_1^2\v_2^2 \Big[ \exp\Big\{-\f{\nu}{C_\ast  a}\tau\Big\}
		+ \exp\Big\{-\f{\v_1}{C_\ast}\tau\Big\}\Big].
	\end{align}
\end{lemma}

\noindent{\bf Proof.} Noting from \eqref{3.15-3} that
\begin{equation*}
	|\partial_y \tilde{v}_X^s|\leq C(v_-)\f{\v_1}{\lambda}|\partial_y w_X| \quad
	\mbox{and} \quad |\partial_y \tilde{h}_X^s|\leq C(v_-)\f{\v_1}{\lambda}|\partial_y w_X|,
\end{equation*}
we shall use the a priori estimates in Lemma \ref{lem4.2} to prove our lemma by taking
\begin{align}\label{5.22-1}
	(\tilde{v}_0,\tilde{h}_0)=(\tilde{v},\tilde{h}),\quad |\mathbf{v}_i|\cong |\partial_y \tilde{v}_X^s|\text{ or }|\partial_y \tilde{h}_X^s| \quad \text{ and }\quad |\mathbf{w}|=|\partial_y w_X|,
\quad \text{for }i=1,2.
\end{align}
Since the proof is very complicate, we divide it into four steps. Here we remark that the shift does not affect the validity of Lemmas \ref{lem4.1}--\ref{lem4.2} since we can change the variable $y \to y-X(\tau)$.
Noting $Y(U)\leq \v_1^2$, Lemma \ref{lem3.4} holds.

\smallskip

\noindent{\it Step 1. Proof of \eqref{5.1}-\eqref{5.3}.}
From the definition of $\tilde{\mathbf{G}}_{2}(v)$ in \eqref{4.7}, one can obtain
\begin{align}\label{5.16}
	\tilde{\mathbf{G}}_2(v)&\leq \int_{\R}  |\partial_y \tilde{v}_X^s| Q(v|\tilde{v})dy + \int_{\R}  |\partial_y \tilde{u}^r| p(v|\tilde{v})dy\nonumber\\
	&\leq C(v_-)\Big\{\f{\v_1}{\lambda} G_2(v) + G^r(v)\Big\}.
\end{align}
It follows from \eqref{5.16}, \eqref{4.99}-\eqref{4.100} and Lemma \ref{lem3.4} that
\begin{align}\label{5.17}
	\begin{split}
		&|B_{1}(v)-B_{1}(\bar{v})| \leq C \delta_{0}\Big[D(v)+\big(G_{2}(v)-G_{2}(\bar{v})\big) + \f{\v_1}{\lambda} G_2(v) + G^r(v)\Big], \\
		&\left|B_{2}^{-}(U)\right| \leq C \delta_{0}\left(D(v)+\f{\v_1}{\lambda} G_2(v) + G^r(v)\right)+\left(\delta_{1}+C \delta_{0}\right) G_{1}^{-}(U), \\
		&\left|B_{2}^{+}(v)-B_{2}^{+}(\bar{v})\right| \leq C\delta_{0} D(v),
	\end{split}
\end{align}
and
\begin{align}\label{5.18}
	|B_{1}(\bar{v})|+\left|B_{2}^{+}(\bar{v})\right| \leq C(v_-) \int_{\mathbb{R}}|\partial_y w_X| \cdot Q\left(\bar{v} | \tilde{v}_{0}\right) d y \leq C_1^{\ast} \frac{\v_1^{2}}{\lambda}.
\end{align}

\smallskip

Using \eqref{4.101} and \eqref{5.22-1}--\eqref{5.16}, one has
\begin{align}\label{5.19}
	|B_3(v)|&\leq |\mathbf{B}_3(v)|\leq C \delta_{0} \Big[\mathbf{D}(v) +  \delta_0 \tilde{\mathbf{G}}_{2}(v)
	+ \Big(\mathbf{G}_{2}(v)-\mathbf{G}_{2}(\bar{v})\Big)+  \f{\v_1}{\lambda}\mathbf{G}_2(\bar{v})\Big]\nonumber\\
	&\leq C \delta_{0} \Big[D(v) + \delta_0 \f{\v_1}{\lambda} G_2(v) + \delta_0 G^r(v)
	+ \big(G_{2}(v)-G_{2}(\bar{v})\big)+  \f{\v_1}{\lambda}G_2(\bar{v})\Big]\nonumber\\
	&\leq C \delta_{0} \Big[D(v)  + \delta_0 G^r(v)
	+ \big(G_{2}(v)-G_{2}(\bar{v})\big)+  \f{\v_1}{\lambda}G_2(\bar{v})\Big],
\end{align}
and
\begin{align}\label{5.20}
	|B_4(v)|
	&\leq \int_{\R} \big|w_{X}\cdot  \partial_y(p(v)-p(\tilde{v}))\cdot (v^{\beta}-\tilde{v}^{\beta})\big| \cdot |p'(\tilde{v})|\cdot (|\partial_y\tilde{v}^s_X|+|\partial_y\tilde{v}^r|) \, dy\nonumber\\
	&\leq  C(v_-)\int_{\R} \big|\partial_y(p(v)-p(\tilde{v}))\cdot (v^{\beta}-\tilde{v}^{\beta})\big| \cdot |\partial_y\tilde{v}^s_X| \, dy\nonumber\\
	&\quad +  C(v_-)\int_{\R} \big|\partial_y(p(v)-p(\tilde{v}))\cdot (v^{\beta}-\tilde{v}^{\beta})\big| \cdot    |\partial_y\tilde{v}^r| \, dy\nonumber\\
	&\leq C \delta_{0} \Big[D(v)  + G^r(v)
	+ \big(G_{2}(v)-G_{2}(\bar{v})\big)+  \f{\v_1}{\lambda}G_2(\bar{v})\Big]\nonumber\\
	&\quad +  \f{C}{\delta_0}\int_{\R}  |\partial_y\tilde{v}^r|^2\cdot  \f{|v^{\beta}-\tilde{v}^{\beta}|^2}{v^{\beta}}   \, dy.
\end{align}

Now we  need to estimate the last term on RHS of \eqref{5.20}. Let $v_\ast\in (0,1)$ be suitably small constant such that
\begin{align}
	p(v|\tilde{v}) {\bf 1}_{\{v\leq v_\ast\}} &\geq [p(v) -p(\tilde{v}) - \gamma (\tilde{v})^{-\gamma}] {\bf 1}_{\{v\leq v_\ast\}}= [p(v) - (\gamma +1) p(\tilde{v})] {\bf 1}_{\{v\leq v_\ast\}} \nonumber\\
	&\geq \f12 p(v) {\bf 1}_{\{v\leq v_\ast\}}\geq \f12 v^{-\beta}\nonumber,
\end{align}
where we have used the condition $0\leq \beta=\gamma-\alpha\leq \gamma$. Hence we obtain
\begin{align}\label{5.22}
	\int_{\R}  |\partial_y\tilde{v}^r|^2\cdot  \f{|v^{\beta}-\tilde{v}^{\beta}|^2}{v^{\beta}}  {\bf 1}_{\{v\leq v_\ast\}} \, dy
	&\leq C\int_{\R}  |\partial_y\tilde{v}^r|^2 \cdot  v^{-\beta}  {\bf 1}_{\{v\leq v_\ast\}} \, dy\nonumber\\
	&\leq C\int_{\R}  |\partial_y\tilde{v}^r|^2 \cdot  p(v|\tilde{v})  {\bf 1}_{\{v\leq v_\ast\}} \, dy\nonumber\\
	&\leq C \v_2 \int_{\R}  |\partial_y\tilde{u}^r|  \cdot  p(v|\tilde{v})  \, dy=C \v_2 G^r(v),
\end{align}
where we have used \eqref{2.5}, \eqref{2.14} and $|\partial_y\tilde{v}^r|=\nu|\partial_{x}\tilde{v}^r|\leq C \f{\nu}{a} \v_2$ .

\smallskip

On the other hand, for $v\geq v^\ast\gg1$ and  noting $0\leq \beta\leq 1$, we have
\begin{align}\label{5.23}
	\int_{\R}  |\partial_y\tilde{v}^r|^2\cdot  \f{|v^{\beta}-\tilde{v}^{\beta}|^2}{v^{\beta}}  {\bf 1}_{\{v\geq v^\ast\}} \, dy
	&\leq C\int_{\R}  |\partial_y\tilde{v}^r|^2 \cdot   v^{\beta}   {\bf 1}_{\{v\geq v^\ast\}} \, dy\nonumber\\
	&\leq C\int_{\R}  |\partial_y\tilde{v}^r|^2\cdot   p(v|\tilde{v})  {\bf 1}_{\{v\geq v^\ast\}} \, dy\nonumber\\
	&\leq C \v_2 \int_{\R}  |\partial_y\tilde{u}^r|  \cdot  p(v|\tilde{v})  \, dy=C \v_2 G^r(v).
\end{align}
For $0<v_\ast \leq v\leq v^\ast$, it holds that
\begin{align}\label{5.24}
	&\int_{\R}  |\partial_y\tilde{v}^r|^2\cdot  \f{|v^{\beta}-\tilde{v}^{\beta}|^2}{v^{\beta}}  {\bf 1}_{\{v_\ast\leq v\leq v^\ast\}} \, dy\leq C(v_\ast,v^{\ast})\int_{\R}  |\partial_y\tilde{v}^r|^2 \cdot  |v-\tilde{v}|^2  {\bf 1}_{\{v_\ast\leq v\leq v^\ast\}} \, dy\nonumber\\
	&\leq C(v_\ast,v^{\ast}) \int_{\R}  |\partial_y\tilde{v}^r|^2 \cdot  p(v|\tilde{v}) {\bf 1}_{\{v_\ast\leq v\leq v^\ast\}}\, dy\leq C \v_2 \int_{\R}  |\partial_y\tilde{u}^r| \cdot   p(v|\tilde{v})  \, dy=C \v_2 G^r(v).
\end{align}
Combining \eqref{5.22}-\eqref{5.24}, one obtains
\begin{align*}
	\int_{\R}  |\partial_y\tilde{v}^r|^2\cdot  \f{|v^{\beta}-\tilde{v}^{\beta}|^2}{v^{\beta}}   \, dy\leq C \v_2 G^r(v),
\end{align*}
which, together with \eqref{5.20}, yields
\begin{align}\label{5.25}
	|B_4(v)|
	&\leq C \delta_{0} \Big[D(v)  + G^r(v)
	+ \big(G_{2}(v)-G_{2}(\bar{v})\big)+  \f{\v_1}{\lambda}G_2(\bar{v})\Big].
\end{align}

\smallskip

\noindent{\it Step 2. Proof of \eqref{5.10}.} Recall \eqref{3.29}, one has
\begin{align}\label{5.26}
	|B_5(U)|&\leq \left|\int_{\R} w_{X}\,   (p(v)-p(\tilde{v})) \cdot F_1\, dy\right| + \left|\int_{\R} w_{X}\, (h-\tilde{h})\cdot F_2  \, dy\right|\nonumber\\
	&=: B_{51}(v) + B_{52}(h).
\end{align}
Noting \eqref{3.7}, it holds that
\begin{align}\label{5.35}
	F_1&=\tilde{v}^{\beta}\, p'(\tilde{v})\, \partial_{yy}\tilde{v} + \big[\beta-(\gamma+1)\big]\,\tilde{v}^{\beta-1} p'(\tilde{v})\, |\partial_y \tilde{v}|^2\nonumber\\
	&\quad-(\tilde{v}_X^s)^{\beta} p'(\tilde{v}_X^s)\, \partial_{yy}\tilde{v}_X^s
	-\big[\beta-(\gamma+1)\big]\,(\tilde{v}_X^s)^{\beta-1} p'(\tilde{v}_X^s)\, |\partial_y \tilde{v}_X^s|^2\nonumber\\
	&=\partial_{yy}\tilde{v}_X^s\, [b_1(\tilde{v})-b_1(\tilde{v}_X^s)] + b_1(\tilde{v}) \, \partial_{yy}\tilde{v}^r + 2b_2(\tilde{v}) \, \partial_y \tilde{v}_X^s\, \partial_y \tilde{v}^r\nonumber\\
	&\quad + |\partial_y \tilde{v}_X^s|^2 \, [b_2(\tilde{v})-b_2(\tilde{v}_X^s)]
	+ b_2(\tilde{v})\, |\partial_y \tilde{v}^r|^2
\end{align}
where we have used the notations $b_1(s)=s^{\beta} p'(s)$, $b_2(s)=\big[\beta-(\gamma+1)\big] s^{\beta-1} p'(s)$.

\smallskip

A direct calculation shows that
\begin{align}\label{5.28}
	&\left|\int_{\R} w_{X} \cdot (p(v)-p(\tilde{v})) \cdot \partial_{yy}\tilde{v}_X^s\cdot [b_1(\tilde{v})-b_1(\tilde{v}_X^s)]\, dy\right|\nonumber\\
	&\leq C(v_-)\v_1 \int_{\Omega^c} |p(v)-p(\tilde{v})|\cdot |\partial_{y}\tilde{v}_X^s|\cdot |\tilde{v}^r-v_m| \, dy\nonumber\\
	&\quad +C(v_-)\v_1 \int_{\Omega} |p(v)-p(\tilde{v})|\cdot |\partial_{y}\tilde{v}_X^s|\cdot |\tilde{v}^r-v_m| \, dy\nonumber\\
	&\leq C(v_-)\v_1 \int_{\Omega^c} |p(v)-p(\bar{v})|\cdot |\partial_{y}\tilde{v}_X^s|\cdot |\tilde{v}^r-v_m| \, dy + C(v_-)\v_1 \int_{\R}  |\partial_{y}\tilde{v}_X^s|\cdot |\tilde{v}^r-v_m| \, dy\nonumber\\
	&\leq  C(v_-)\v_1 \sqrt{\f{\v_1}{\lambda}}\left(\int_{\Omega^c} |\partial_{y}w_X|\cdot|p(v)-p(\tilde{v})|^2 \, dy\right)^{\f12}
	\left(\int_{\Omega^c}  |\partial_{y}\tilde{v}_X^s|\cdot |\tilde{v}^r-v_m|^2 \, dy\right)^{\f12}\nonumber\\
	&\quad +  C(v_-)\v_1 \int_{\R}  |\partial_{y}\tilde{v}_X^s|\cdot |\tilde{v}^r-v_m| \, dy.
\end{align}
Applying \eqref{4.15} and using \eqref{5.16}, one obtains
\begin{align}\label{5.37}
	\int_{\Omega^c} |\partial_{y}w_X|\cdot|p(v)-p(\bar{v})|^2 \, dy&\leq C\sqrt{\f{\v_1}{\lambda}}\Big[D(v) + \delta_0\f{\v_1}{\lambda} G_2(v) + \delta_0 G^r(v)\Big]\nonumber\\
	&\quad + C\sqrt{\f{\v_1^{2-q}}{\lambda}}\Big[D(v) + \delta_0\f{\v_1}{\lambda} G_2(v) + \delta_0 G^r(v)\Big]^q,
\end{align}
where $1\leq q\leq \frac{2}{2\gamma-1}<2$. For the interaction term of viscous shock wave and the approximate rarefaction wave, it follows from \eqref{5.10-1} and Lemma \ref{lem2.3} that
\begin{align}\label{5.38}
	\int_{\R} |\partial_y \tilde{v}_X^s|\cdot |\tilde{v}^r-v_m| \, dy
	&\leq C(v_-)\v_2 \int_{-\infty}^0  |\partial_y \tilde{v}^s_X|\exp\Big(-\f{\nu}{a}(|y|+\lambda_2(v_m)\tau)\Big) dy\nonumber\\
	&\quad + C(v_-) \v_1^2 \v_2 \int_0^\infty  \exp\Big(-c_1\v_1(|y|+\f12|\sigma_1|\tau)\Big) dy\nonumber\\
	&\leq C(v_-)\v_1\v_2 \Big[\exp\big\{-\f{\nu}{C_\ast a}\tau\big\}
	+ \exp\big\{-\f{\v_1}{C_\ast}\tau\big\}\Big],
\end{align}
which, together with \eqref{5.28}-\eqref{5.37}, yields
\begin{align}\label{5.39}
	&\left|\int_{\R} w_{X}\cdot   (p(v)-p(\tilde{v})) \cdot \partial_{yy}\tilde{v}_X^s\cdot  [b_1(\tilde{v})-b_1(\tilde{v}_X^s)]\, dy\right|\nonumber\\
	&\leq  \v_1 \Big[D(v) + \delta_0\f{\v_1}{\lambda} G_2(v) + \delta_0 G^r(v)\Big] \nonumber\\
	&\quad+ C(v_-)\v_1\v_2 \Big[\exp\big\{-\f{\nu}{C_\ast a}\tau\big\}
	+ \exp\big\{-\f{\v_1}{C_\ast}\tau\big\}\Big].
\end{align}
Similarly, we have
\begin{align}\label{5.40}
	&\left|\int_{\R} w_{X}\,   (p(v)-p(\tilde{v})) \cdot |\partial_y \tilde{v}_X^s|^2 \, [b_2(\tilde{v})-b_{2}(\tilde{v}_X^s)]\, dy\right|\nonumber\\
	&\leq \v_1 \Big[D(v) + \delta_0\f{\v_1}{\lambda} G_2(v) + \delta_0 G^r(v)\Big] \nonumber\\
	&\quad+C(v_-)\v_1\v_2 \Big[\exp\big\{-\f{\nu}{C_\ast a}\tau\big\}
	+ \exp\big\{-\f{\v_1}{C_\ast}\tau\big\}\Big].
\end{align}

\smallskip

It follows from \eqref{2.3-3}, \eqref{2.14} and  \eqref{5.12-2} that
\begin{align}
	\int_{\R} |\partial_y\tilde{v}_X^s\cdot \partial_y\tilde{v}^r|dy
	&\leq \int_{-\infty}^0 |\partial_y\tilde{v}_X^s\cdot \partial_y\tilde{v}^r|dy
	+\int_0^{\infty} |\partial_y\tilde{v}_X^s\cdot \partial_y\tilde{v}^r|dy\nonumber\\
	&\leq C(v_-) \f{\nu \v_1\v_2}{a} \exp\Big\{-\f{\nu}{C_\ast a}\tau\Big\} + C(v_-)\v_1\v_2 \exp\Big\{-\f{\v_1 \tau}{C_\ast}\Big\},\nonumber
\end{align}
which, together with \eqref{5.37}, implies
\begin{align}\label{5.41}
	&\left|\int_{\R} w_{X}\cdot  (p(v)-p(\tilde{v})) \cdot 2b_2(\tilde{v}) \cdot  \partial_y \tilde{v}_X^s\cdot  \partial_y \tilde{v}^r\, dy\right|\nonumber\\
	&\leq C(v_-) \int_{\Omega^c} |p(v)-p(\tilde{v})|\cdot |\partial_{y}\tilde{v}_X^s|\cdot |\partial_y \tilde{v}^r| \, dy  +C(v_{-})\int_{\Omega} |p(v)-p(\tilde{v})|\cdot |\partial_{y}\tilde{v}_X^s|\cdot |\partial_y \tilde{v}^r| \, dy\nonumber\\
	&\leq  C(v_-)\sqrt{\f{\v_1}{\lambda}}\left(\int_{\Omega^c} |\partial_{y}w_X|\cdot|p(v)-p(\bar{v})|^2 \, dy\right)^{\f12}
	\left(\int_{\Omega^c}  |\partial_{y}\tilde{v}_X^s|\cdot |\partial_y \tilde{v}^r|^2 \, dy\right)^{\f12}\nonumber\\
	&\quad +  C(v_-) \int_{\Omega}  |\partial_{y}\tilde{v}_X^s|\cdot |\partial_y \tilde{v}^r| \, dy\nonumber\\
	&\leq  	C\left(\f{\v_1}{\lambda}\right)^{\f14} \Big[D(v) + \delta_0\f{\v_1}{\lambda} G_2(v) + \delta_0 G^r(v)\Big] \nonumber\\
	&\qquad +	C \f{\nu \v_1\v_2}{a} \exp\Big\{-\f{\nu}{C_\ast a}\tau\Big\} + C\v_1\v_2 \exp\Big\{-\f{\v_1 \tau}{C_\ast}\Big\},
\end{align}
where we have used $\v_1\leq \lambda^2$ in the last inequality.

\smallskip

We consider the term involving $|\partial_y\tilde{v}^r|^2$ on RHS of \eqref{5.35}. It follows from \eqref{2.12} that
\begin{align}\label{5.43}
	\int_{\R} |\partial_y\tilde{v}^r|^2 dy
	&= \nu \int_{\R} |\partial_x\tilde{v}^r|^2 dx
	\leq C\nu \left(\min\Big\{\f{\v_2}{\sqrt{a}},\, \sqrt{\v_2} t^{-\f12}\Big\}\right)^2\leq C\v_2 \min\Big\{\f{\nu\v_2}{a},\, \tau^{-1}\Big\}.
\end{align}
Let $k_1:=\f12 p(v_-)$ and $ k_2:=2p(v_-) +1$, we define
\begin{align}
	p(v_k)-p(\tilde{v})
	:=
	\begin{cases}
		-k_1, \qquad\qquad\,  \mbox{if}\,\, p(v)-p(\tilde{v})<-k_1,\\
		p(v)-p(\tilde{v}),\quad\mbox{if} -k_1\leq p(v)-p(\tilde{v})\leq k_2,\\
		k_2, \qquad\qquad\quad \mbox{if}\,\, p(v)-p(\tilde{v})> k_2.
	\end{cases}\nonumber
\end{align}
It is direct to check that
\begin{align}\label{5.45}
	\begin{split}
		|p(v)-p(\tilde{v})|&\leq C_1(v_-) p(v|\tilde{v}),\quad \mbox{if}\,\, p(v)-p(\tilde{v})\leq -k_1,\\
		|p(v)-p(\tilde{v})|&\leq C_1(v_-) p(v|\tilde{v}),\quad \mbox{if}\,\, p(v)-p(\tilde{v})\geq k_2,
	\end{split}
\end{align}
where $C_1(v_-)>0$ is some positive constant depending only on $v_-$.
Clearly it holds
\begin{align}\label{5.46}
	&\left|\int_{\R}  w_{X} \cdot  (p(v)-p(\tilde{v})) \cdot b_2(\tilde{v})\cdot  |\partial_y \tilde{v}^r|^2 dy\right|\nonumber\\
	&\leq C(v_-)\int_{\R}  |p(v)-p(\tilde{v})|\cdot  |\partial_y \tilde{v}^r|^2 {\bf 1}_{\{-k_1\leq p(v)-p(\tilde{v})\leq k_2\}} dy\nonumber\\
	&\quad + C(v_-)\int_{\R}  |p(v)-p(\tilde{v})|\cdot  |\partial_y \tilde{v}^r|^2 {\bf 1}_{\{ p(v)-p(\tilde{v})< -k_1\}} dy\nonumber\\
	&\quad + C(v_-)\int_{\R}  |p(v)-p(\tilde{v})|\cdot  |\partial_y \tilde{v}^r|^2 {\bf 1}_{\{p(v)-p(\tilde{v})> k_2\}} dy\nonumber\\
	&=:J_1 + J_2+ J_3.
\end{align}
For $J_1$, it follows from \eqref{5.43} that
\begin{align}\label{5.47}
	J_1&\leq  C(v_-)\int_{\R}  |p(v)-p(\tilde{v})|\cdot  |\partial_y \tilde{v}^r|^2 {\bf 1}_{\{-k_1\leq p(v)-p(\tilde{v})\leq k_2\}} dy\nonumber\\
	&\leq C(v_-)\int_{\R}  |p(v_k)-p(\tilde{v})|\cdot |\partial_y \tilde{v}^r|^2  dy\nonumber\\
	&\leq C(v_-)\|p(v_k)-p(\tilde{v})\|_{L^\infty} \int_{\R} |\partial_y \tilde{v}^r|^2  dy\nonumber\\
	&\leq C(v_-)\v_2 \min\Big\{\f{\nu\v_2}{a},\, \tau^{-1}\Big\}\, \|\partial_y(p(v_k)-p(\tilde{v}))\|_{L^2}^{\f12}\cdot \|p(v_k)-p(\tilde{v})\|_{L^2}^{\f12}\nonumber\\
	&\leq  \f{C(v_-)}{\delta_0} \|p(v_k)-p(\tilde{v})\|_{L^2}^{\f23} \left(\v_2 \min\Big\{\f{\nu\v_2}{a},\, \tau^{-1}\Big\}\right)^{\f43}
	+\delta_0 \|\partial_y(p(v_k)-p(\tilde{v}))\|_{L^2}^2  \nonumber\\
	&\leq  \f{C(v_-)}{\delta_0} \v_2^{\f43} \left(\min\Big\{\f{\nu\v_2}{a},\, \tau^{-1}\Big\}\right)^{\f43} \Big\{1+ \int_{\R}|p(v_k)-p(\tilde{v})|^2dy\Big\} \nonumber\\
	&\quad+C(v_-)\delta_0 \int_{\R} w_X\cdot  v_k^{\beta}\cdot  |\partial_y(p(v_k)-p(\tilde{v}))|^2dy  \nonumber\\
	&\leq \f{C(v_-)}{\delta_0} \v_2^{\f43} \left(\min\Big\{\f{\nu\v_2}{a},\, \tau^{-1}\Big\}\right)^{\f43} \Big\{1+ \int_{\R} Q(v_k|\tilde{v})dy\Big\} \nonumber\\
	&\quad+C(v_-)\delta_0 \int_{\R} w_X\cdot  v^{\beta}\cdot  |\partial_y(p(v)-p(\tilde{v}))|^2dy\nonumber\\
	&\leq \f{C(v_-)}{\delta_0} \v_2^{\f43} \left(\min\Big\{\f{\nu\v_2}{a},\, \tau^{-1}\Big\}\right)^{\f43}\Big(1+\int_{\R} \eta(U|\tilde{U})dy\Big)
	+ C(v_-) \delta_0 D(v).
\end{align}
For $J_2, J_3$, using \eqref{5.45}, one obtains
\begin{align}\label{5.48}
	J_2+J_3 &\leq  C(v_-)\int_{\R} p(v|\tilde{v})\, |\partial_y \tilde{v}^r|^2 dy \leq C(v_-) \f{\nu\v_2}{a} G^r(v).
\end{align}
Substituting \eqref{5.47}-\eqref{5.48} into \eqref{5.46}, one has
\begin{align}\label{5.49}
	&\left|\int_{\R}  w_{X}\cdot    (p(v)-p(\tilde{v})) \cdot b_2(\tilde{v})\cdot |\partial_y \tilde{v}^r|^2 dy\right|\nonumber\\
	&\leq C(v_-)\Big[  \delta_0 D(v) + \f{\nu\v_2}{a} G^r(v)\Big] + \f{C(v_-)}{\delta_0} \v_2^{\f43} \left(\min\Big\{\f{\nu\v_2}{a},\, \tau^{-1}\Big\}\right)^{\f43}\Big(1+\int_{\R} \eta(U|\tilde{U})dy\Big).
\end{align}

\smallskip

For the term involving $|\partial_{yy}\tilde{v}^r|^2$ on the right hand side of \eqref{5.35}, we notice from \eqref{2.12} that
\begin{align}\label{5.50}
	\int_{\R} |\partial_{yy}\tilde{v}^r|dy
	&\leq \nu\int_{\R} |\partial_{xx}\tilde{v}^r|dx\leq C(v_-)\, \nu  \, \min\{\f{\v_2}{a},\, \f1t\}\leq C(v_-)\,  \min\{\f{\nu  \v_2 }{a},\, \tau^{-1}\}.
\end{align}
On the other hand, a direct calculation shows that
\begin{align}\label{5.51}
	|\partial_{yy}\tilde{v}^r|\leq C(v_-)\f{\nu}{a} \, |\partial_y\tilde{v}^r|.
\end{align}
Clearly, it holds that
\begin{align}\label{5.52}
	&\left|\int_{\R}  w_{X}\cdot   (p(v)-p(\tilde{v})) \cdot b_2(\tilde{v})\cdot |\partial_{yy}\tilde{v}^r| dy\right|\nonumber\\
	&\leq C(v_-)\int_{\R}  |p(v)-p(\tilde{v})|\cdot |\partial_{yy}\tilde{v}^r| {\bf 1}_{\{-k_1\leq p(v)-p(\tilde{v})\leq k_2\}} dy\nonumber\\
	&\quad + C(v_-)\int_{\R}  |p(v)-p(\tilde{v})|\cdot|\partial_{yy}\tilde{v}^r| {\bf 1}_{\{ p(v)-p(\tilde{v})< -k_1\}} dy\nonumber\\
	&\quad + C(v_-)\int_{\R}  |p(v)-p(\tilde{v})|\cdot |\partial_{yy}\tilde{v}^r| {\bf 1}_{\{p(v)-p(\tilde{v})> k_2\}} dy\nonumber\\
	&=:J_4 + J_5+ J_6.
\end{align}
For $J_4$, it follows from \eqref{5.50} that
\begin{align}\label{5.53}
	J_4&\leq  C(v_-)\int_{\R}  |p(v)-p(\tilde{v})|\cdot |\partial_{yy} \tilde{v}^r| {\bf 1}_{\{-k_1\leq p(v)-p(\tilde{v})\leq k_2\}} dy\nonumber\\
	&\leq C(v_-)\int_{\R}  |p(v_k)-p(\tilde{v})|\cdot |\partial_{yy} \tilde{v}^r|  dy\nonumber\\
        &\leq C(v_-)\|p(v_k)-p(\tilde{v})\|_{L^\infty} \int_{\R} |\partial_{yy} \tilde{v}^r| dy\nonumber\\
	&\leq C(v_-) \min\Big\{\f{\nu\v_2}{a},\, \tau^{-1}\Big\}\, \|\partial_y(p(v_k)-p(\tilde{v}))\|_{L^2}^{\f12}\|p(v_k)-p(\tilde{v})\|_{L^2}^{\f12}\nonumber\\
	&\leq \f{C(v_-)}{\delta_0} \left(\min\Big\{\f{\nu\v_2}{a},\, \tau^{-1}\Big\}\right)^{\f43}\Big(1+\int_{\R} \eta(U|\tilde{U})dy\Big)
	+ C(v_-) \delta_0 D(v).
\end{align}
For $J_5, J_6$, using \eqref{5.45} and \eqref{5.51}, one can get
\begin{align}\label{5.53-1}
	J_5 +  J_6&\leq C(v_-)\int_{\R} p(v|\tilde{v})\cdot  |\partial_{yy} \tilde{v}^r| dy\leq C(v_-)\f{\nu}{a}\int_{\R} p(v|\tilde{v})\cdot  |\partial_{y} \tilde{v}^r| dy \leq C(v_-) \f{\nu}{a} G^r(v),
\end{align}
which, together with \eqref{5.52}-\eqref{5.53}, yields that
\begin{align}\label{5.54}
	&\left|\int_{\R}  w_{X}\cdot    (p(v)-p(\tilde{v})) \cdot b_2(\tilde{v})\cdot |\partial_{yy}\tilde{v}^r| dy\right|\nonumber\\
	&\leq C(v_-)\Big[  \delta_0 D(v) + \f{\nu}{a} G^r(v)\Big] + \f{C(v_-)}{\delta_0} \left(\min\Big\{\f{\nu\v_2}{a},\, \tau^{-1}\Big\}\right)^{\f43}\Big(1+\int_{\R} \eta(U|\tilde{U})dy\Big).
\end{align}

\smallskip

Combining \eqref{5.35}, \eqref{5.39}-\eqref{5.40}, \eqref{5.41}, \eqref{5.49} and \eqref{5.54}, one obtains
\begin{align}\label{5.55}
	B_{51}(v)
	&\leq C(v_-) \big[\delta_0 + \v_1 + \left(\f{\v_1}{\lambda}\right)^{\f14}\big] D(v) + C(v_-)\big[\v_1+ \delta_0\left(\f{\v_1}{\lambda}\right)^{\f14}+\f{\nu}{a}\big]\, G^r(v) \nonumber\\
	&\quad + C(v_-)\delta_0 \f{\v_1}{\lambda} G_2(v)  + \f{C(v_-)}{\delta_0} \left(\min\Big\{\f{\nu\v_2}{a},\, \tau^{-1}\Big\}\right)^{\f43}\Big(1+\int_{\R} \eta(U|\tilde{U})dy\Big) \nonumber\\
	&\quad +C\v_1\v_2\Big[  \exp\big\{-\f{\nu}{C_\ast a}\tau\big\} +    \exp\big\{-\f{\v_1}{C_\ast}\tau\big\}\Big].
\end{align}

\medskip

Next, we consider the estimate of $B_{52}(h)$. A direct calculation shows
\begin{align}\label{5.56}
	B_{52}(h) &\leq \left|\int_{\R} w_{X}\cdot  (h-\tilde{h})\cdot \partial_y \tilde{v}^r\cdot  [p'(\tilde{v})-p'(\tilde{v}^r)]  \, dy\right|\nonumber\\
	&\quad + \left|\int_{\R} w_{X}\cdot  (h-\tilde{h})\cdot \partial_y \tilde{v}^s_X \cdot [p'(\tilde{v})-p'(\tilde{v}^s_X)]  \, dy\right|.
\end{align}

It follows from \eqref{5.11-1} and \eqref{2.12}  that
\begin{align}\label{5.57}
	&\int_{\R} |\partial_y \tilde{v}^r|^2\cdot |\tilde{v}^s_X-v_m|^2 dy\leq \int_{-\infty}^0 |\partial_y \tilde{v}^r|^2\cdot |\tilde{v}^s_X-v_m|^2 dy
	+ \int_0^{\infty} |\partial_y \tilde{v}^r|^2\cdot |\tilde{v}^s_X-v_m|^2 dy\nonumber\\
	&\leq C(v_-)\v_1^2\v_2^2 \f{\nu^2}{a^2} \int_{-\infty}^0 \exp\Big(-\f{2\nu}{a}[|y|+\lambda_2(v_m)\tau]\Big) dy \nonumber\\
	&\quad + C(v_-)\v_1^2  \Big(\min\{\f{\nu\v_2}{a},\, \tau^{-1}\}\Big)^2 \int_0^{\infty} \exp\Big\{-2c_1\v_1 (|y|+\f12 |\sigma_1|\tau)\Big\} dy\nonumber\\
	&\leq C(v_-)\v_1^2\v_2^2 \f{\nu}{a} \exp\Big\{-\f{\nu}{C_\ast a}\tau\Big\}
	+ C(v_-)\v_1  \Big(\min\{\f{\nu\v_2}{a},\, \tau^{-1}\}\Big)^2  \exp\Big\{-\f{\v_1}{C_\ast}\tau\Big\}.
\end{align}
Then  the first term on RHS of \eqref{5.56} can be controlled as
\begin{align}\label{5.58}
	&\left|\int_{\R} w_{X}\cdot  (h-\tilde{h})\cdot \partial_y \tilde{v}^r\cdot [p'(\tilde{v})-p'(\tilde{v}^r)]  \, dy\right|\nonumber\\
	&\leq C(v_-)\left(\int_{\R} |h-\tilde{h}|^2\, dy\right)^{\f12}
	\left(\int_{\R} |\partial_y \tilde{v}^r|^2\cdot |v^s_X-v_m|^2dy\right)^{\f12}\nonumber\\
	&\leq C(v_-)\left[\v_1\v_2\sqrt{\f{\nu}{a}} \exp\Big\{-\f{\nu}{C_\ast a}\tau\Big\}
	+ \sqrt{\v_1}  \min\{\f{\nu\v_2}{a},\, \tau^{-1}\}  \exp\Big\{-\f{\v_1}{C_\ast}\tau\Big\}\right]\nonumber\\
	&\qquad\quad \times \Big(\int_{\R} |h-\tilde{h}|^2 dy\Big)^{\f12}\nonumber\\
	&\leq \Big[\f{\nu}{a}\exp\Big\{-\f{\nu}{C_\ast a}\tau\Big\} + \v_1 \exp\Big\{-\f{\v_1}{C_\ast}\tau\Big\}\Big]\int_{\R} |h-\tilde{h}|^2 dy\nonumber\\
	&\quad + C(v_-)\v_1^2\v_2^2 \exp\Big\{-\f{\nu}{C_\ast a}\tau\Big\}
	+ C(v_-) \Big(\min\{\f{\nu\v_2}{a},\, \tau^{-1}\} \Big)^2.
\end{align}

\smallskip

For the second term on RHS of \eqref{5.56}, it follows from \eqref{3.15-3} and \eqref{5.38}  that
\begin{align}
	&\left|\int_{\R} w_{X}\cdot  (h-\tilde{h})\cdot \partial_y \tilde{v}^s_X\cdot  [p'(\tilde{v})-p'(\tilde{v}^s_X)]  \, dy\right|\nonumber\\
	&\leq \left(\int_{\R} |\partial_y \tilde{v}^s_X|\cdot |h-\tilde{h}|^2 dy\right)^{\f12}\cdot \left(\int_{\R} |\partial_y \tilde{v}^s_X|\cdot |\tilde{v}^r-v_m|^2 dy \right)^{\f12}\nonumber\\
	&\leq C(v_-)\sqrt{\f{\v_1\v_2}{\lambda}} \left(\int_{\R} |\partial_yw_X|\cdot |h-\tilde{h}|^2 dy\right)^{\f12}\cdot \left(\int_{\R} |\partial_y \tilde{v}^s_X|\cdot |\tilde{v}^r-v_m| dy \right)^{\f12}\nonumber\\
	&\leq \delta_0 \v_1\v_2 \int_{\R} \f{|\sigma_1|}{2} |\partial_yw_X|\cdot |h-\tilde{h}|^2 dy
	+C(v_-) \f{\v_1\v_2}{\delta_0\lambda} \left(\exp\big\{-\f{\nu}{C_\ast a}\tau\big\}
	+   \exp\big\{-\f{\v_1}{C_\ast}\tau\big\}\right)\nonumber\\
	&\leq \delta_0 \v_1\v_2\, G^h(h)  +C(v_-) \f{\v_1\v_2}{\delta_0\lambda} \left(\exp\big\{-\f{\nu}{C_\ast a}\tau\big\}
	+   \exp\big\{-\f{\v_1}{C_\ast}\tau\big\}\right),\label{5.58-1}
\end{align}
which, together with \eqref{5.58}, yields
\begin{align}\label{5.60}
	B_{52}(h)&\leq \delta_0 \v_1\v_2 \, G^h(h)  +C(v_-) \f{\v_1\v_2}{\delta_0\lambda} \left(\exp\big\{-\f{\nu}{C_\ast a}\tau\big\}
	+   \exp\big\{-\f{\v_1}{C_\ast}\tau\big\}\right)\nonumber\\
	&\quad +\Big[\f{\nu}{a}\exp\big\{-\f{\nu}{C_\ast a}\tau\big\} + \v_1 \exp\big\{-\f{\v_1}{C_\ast}\tau\big\}\Big]\int_{\R} |h-\tilde{h}|^2 dy\nonumber\\
     &\quad +C(v_{-})\left(\min\{\frac{\nu \varepsilon_{0}}{a},\tau^{-1}\}\right)^2.
\end{align}
Then, using \eqref{5.60} and \eqref{5.55}, we conclude \eqref{5.10}.

\medskip

\noindent{\it Step 3. Proof of \eqref{5.11}. } Recall the definition of $B_{\delta}(U)$ in \eqref{3.29}. Using $\dis\f{\nu}{a}\leq \delta_0$,  \eqref{5.1}-\eqref{5.10} and Lemma \ref{lem3.4}, we have
\begin{align}
	|B_{\delta_1}(U)|
	&\leq |B_1(v)| + |B_2^-(U)| + 2 |B_2^+(v)| + |B_3(v)| + |B_4(v)| + |B_5(U)|\nonumber\\
	&\leq |B_1(v)-B_1(\bar{v})| + 2|B_2^+(v)-B_2^+(\bar{v})| + |B_1(\bar{v})| + 2|B_2^+(\bar{v})|\nonumber\\
	&\quad + |B_2^-(U)|+  |B_3(v)| + |B_4(v)| + |B_5(U)|\nonumber\\
	&\leq  C \delta_{0} \Big(D(v)+[G_{2}(v)-G_{2}(\bar{v})] +\f{\v_1}{\lambda}G_{2}(v) + G^r(v)\Big) + 2C_1^\ast \f{\v_1^2}{\lambda}\nonumber\\
	&\quad +\left(\delta_{1}+C \delta_{0}\right) G_{1}^{-}(U) +\delta_0 \v_1\v_2 \, G^h(h)\nonumber\\
	&\quad  +  f(\nu,a,\v_1,\v_2,\tau)  \int_{\R} \eta(U|\tilde{U}) dy   + g(\nu,a,\v_1,\v_2,\tau)\nonumber\\
	&\leq  C  \delta_{0} D(v) + C\delta_0 G^r(v)  +C  \delta_{0}G_{2}(v) + \left(\delta_{1}+C \delta_{0}\right) G_{1}^{-}(U)\nonumber\\
	&\quad +\delta_0 \v_1\v_2 \, G^h(h) + 2C_1^\ast \f{\v_1^2}{\lambda}  +  f(\nu,a,\v_1,\v_2,\tau) \int_{\R} \eta(U|\tilde{U}) dy   + g(\nu,a,\v_1,\v_2,\tau)\nonumber\\
	&\leq C \delta_0  \Big[D(v) + G^r(v)\Big] +  f(\nu,a,\v_1,\v_2,\tau)  \int_{\R} \eta(U|\tilde{U}) dy \nonumber\\
	&\quad  + 4C_1^\ast  \f{\v_1^2}{\lambda} +  g(\nu,a,\v_1,\v_2,\tau).\nonumber
\end{align}

\medskip

\noindent{\it Step 4. Proof of \eqref{5.13}.} For $\mathbb{Y}^s_5(U)$, using \eqref{5.38}, one has
\begin{align*}
	|\mathbb{Y}^s_5(v)|&\leq C(v_-)\int_{\R} |\partial_y\tilde{v}_X^s\cdot (\tilde{v}^r-v_m)|\cdot |v-\tilde{v}| dy\nonumber\\
	&\leq C(v_-)\int_{\R} |\partial_y\tilde{v}_X^s\cdot (\tilde{v}^r-v_m)|\cdot |v-\tilde{v}| {\bf 1}_{\{v\leq 1+2v_-\}} dy\nonumber\\
	&\quad +C(v_-)\int_{\R} |\partial_y\tilde{v}_X^s\cdot (\tilde{v}^r-v_m)|\cdot |v-\tilde{v}| {\bf 1}_{\{v\geq 1+2v_-\}} dy\nonumber\\
	&\leq C(v_-)\int_{\R} |\partial_y\tilde{v}_X^s\cdot (\tilde{v}^r-v_m)|{\bf 1}_{\{v\leq 1+2v_-\}} dy\nonumber\\
	&\quad +C(v_-)\int_{\R} |\partial_y\tilde{v}_X^s\cdot (\tilde{v}^r-v_m)|\cdot Q(v|\tilde{v}) {\bf 1}_{\{v\geq 1+2v_-\}} dy\nonumber\\
	&\leq C(v_-) \f{\v_1\v_2}{\lambda} \int_{\R} |\partial_yw_X|\cdot  Q(v|\tilde{v}) dy+
	C\v_1\v_2 \Big[ \exp\Big\{-\f{\nu}{C_\ast a}\tau\Big\}
	+ \exp\Big\{-\f{\v_1}{C_\ast}\tau\Big\}\Big]\nonumber\\
	&\leq C(v_-)\f{\v_1\v_2}{\lambda} G_2(v)+
	C(v_-)\v_1\v_2 \Big[ \exp\Big\{-\f{\nu}{C_\ast a}\tau\Big\}
	+ \exp\Big\{-\f{\v_1}{C_\ast}\tau\Big\}\Big],
\end{align*}
which, together with \eqref{4.102}, \eqref{5.16} and \eqref{3.31}, yields
\begin{align*}
	&|\mathbb{Y}^{g}(v)-\mathbb{Y}^{g}(\bar{v})|^2+|\mathbb{Y}^{b}(U)|^2+|\mathbb{Y}^{l}(U)|^2+|\mathbb{Y}^{s}(U)|^2 \nonumber\\
	&\leq C\f{\v_1^2}{\lambda} \bigg( \delta_{0}  \mathbf{D}(v)+\left(\mathbf{G}_{2}(v)-\mathbf{G}_{2}(\bar{v})\right)+\delta_{0} \sqrt{\f{\v_1}{\lambda}} \tilde{\mathbf{G}}_{2}(v)+\left(\frac{\v_1}{\lambda}\right)^{1 / 4} \mathbf{G}_{2}(\bar{v})\nonumber \\
	&\qquad+\mathbf{G}_{1}^{-}(U)+\left(\frac{\lambda}{\v_1}\right)^{1 / 4} \mathbf{G}_{1}^{+}(v)\bigg) + 4 |\mathbb{Y}^s_5(v)|^2\nonumber\\
	&\leq C\f{\v_1^2}{\lambda} \bigg( \delta_{0}  D(v)+\left(G_{2}(v)-G_{2}(\bar{v})\right)  +\delta_0 G^r(v) +\left(\frac{\v_1}{\lambda}\right)^{\f14} G_{2}(\bar{v})\nonumber \\
	&\qquad+G_{1}^{-}(U)+\left(\frac{\lambda}{\v_1}\right)^{\f14} G_{1}^{+}(v)\bigg) + 4 |\mathbb{Y}^s_5(v)|^2\nonumber\\
	&\leq C\f{\v_1^2}{\lambda}  \bigg\{ \delta_{0} \Big[D(v) + G^r(v)\Big]+\Big[G_{2}(v)-G_{2}(\bar{v})\Big] +\left(\frac{\v_1}{\lambda}\right)^{\f14} G_{2}(\bar{v})\nonumber \\
	&\qquad+G_{1}^{-}(U)+\left(\frac{\lambda}{\v_1}\right)^{\f14} G_{1}^{+}(U)\bigg\} + C\v_1^2\v_2^2 \Big[ \exp\Big\{-\f{\nu}{C_\ast a}\tau\Big\}
	+ \exp\Big\{-\f{\v_1}{C_\ast}\tau\Big\}\Big].
\end{align*}
Therefore the proof of Lemma \ref{lem5.2} is complete. $\hfill\Box$

\section{Proof of the Theorem \ref{thm3.1}}
In this section, we shall establish the following global estimate for the compressible Navier-Stokes equations \eqref{1.1}. The global estimate is crucial to prove our main Theorem \ref{thm1.1}.
\begin{theorem}\label{thm3.1}
	Recall Lemmas \ref{lem3.1} and \ref{lem3.3}. There exist small positive constants $\lambda, \v_0, \delta_0$ and $\delta_1$, such that if $0<\v_1\leq \v_0$, $0< \v_2\leq C(v_-)\v_1$,   $0<\nu\leq \delta_0$ and $\v_1\ll\lambda \ll \delta_0\ll \delta_1$,  then it holds
	\begin{align}\label{1.20}
		&\int_{\R} w_X(y-\sigma_1\tau) \,\eta(U|\tilde{U})(\tau) dy  + \int_0^\infty \Big(1+ |J^{bad}| + J^{good}\Big) (\tau) \, {\bf 1}_{\{Y(U)\geq \v_1^2\}} \, d\tau \nonumber\\
		&\quad +\int_0^{\infty} \delta_0\f{\v_1}{\lambda} |B_{\delta_1}(U)(\tau)|  \,  {\bf 1}_{\{Y(U)\leq \v_1^2\}}\, d\tau
		+ \int_0^{\infty} \bigg(\f14 \, G_1^+(U) + \f14 \, G_1^-(U) \nonumber\\
		&\quad + \delta_0 \f{\v_1}{\lambda} G_2(v)   + \f12\v_1\v_2 \, G^h(h) + \delta_0D(v) + \delta_0 G^r(v) + G_3(v)\bigg)(\tau) {\bf 1}_{\{Y(U)\leq \v_1^2\}} \, d\tau \nonumber\\
		&\leq C(U_-,\delta_0^{-1}) \bigg[\int_{\R} w_X(y)\,\eta(U|\tilde{U})(0) dy +   (\f{\nu\v_2}{a})^{\f16}  + \v_2 \f{a}{\nu}\bigg],
	\end{align}
	and
	\begin{align}\label{1.20-2}
		X(\tau)\leq \f12 |\sigma_1|\tau,\qquad \forall \,\tau\geq 0.
	\end{align}
	Moreover, there exists a function $\mathfrak{f}(\tau)$ such that
	\begin{align}\label{1.20-3}
		|\dot{X}(\tau)|
		&\leq C(U_-,\,\v_1^{-1})\, \Big\{\mathfrak{f}(\tau) + 1\Big\},
	\end{align}
	with
	\begin{align}\label{1.20-4}
		\int_0^\infty \mathfrak{f}(\tau) d\tau \leq C(U_-, \v_1^{-1}, \v_2^{-1}) \, \bigg[\int_{\R}  \eta(U|\tilde{U})(0) dy +   (\f{\nu\v_2}{a})^{\f16}  + \v_2 \f{a}{\nu}\bigg].
	\end{align}
\end{theorem}

\begin{remark}\label{rem5.2}
	If $\v_2=0$, then our global estimate \eqref{1.20} degenerate to the case of single shock wave, see \cite{Kang-Vasseur-2021-Invent}. On the other hand, since $a(\nu)/\nu\to +\infty$ as $\nu\to 0$, the estimate \eqref{1.20} is not uniformly bounded with respect to $\nu$. Thus the contraction property as in \cite{Kang-Vasseur-2021-JEMS} is invalid. Fortunately, it is enough to obtain stability and uniqueness of the Riemann solution in the inviscid limit in original coordinate $(t,x)$.
\end{remark}

\smallskip

Recall Lemmas \ref{lem3.1} and \ref{lem3.3}, we denote
\begin{align}
	\mathcal{F}(U)(\tau):=\dot{X}(\tau)\, Y(U) + J^{bad}-J^{good}\equiv \dot{X}(\tau)\, Y(U) + B_{\delta_1}(U)-G_{\delta_1}(U),\nonumber
\end{align}
where $\dot{X}(\tau)$ is the one define in \eqref{5.7-1}.

\smallskip

For the case of $Y(U)\geq \v_1^2$, then it follows from \eqref{5.7-1} that
\begin{align}\label{6.2}
	\mathcal{F}(U)(\tau)\, {\bf 1}_{\{Y(U)\geq \v_1^2\}}
	&\leq -2|J^{bad}|-1-J^{good}\leq 0.
\end{align}

For the case of $Y(U)\leq \v_1^2$, noting \eqref{5.7-1}, a direct calculation shows that
\begin{align}
	\dot{X}(\tau)\, Y(U) {\bf 1}_{\{Y(U)\leq \v_1^2\}}
	&=-\f{1}{\v_1^4} |Y(U)|^2 \, \Big(2|J^{bad}| + 1\Big) \, {\bf 1}_{\{0\leq Y(U)\leq \v_1^2\}} \nonumber\\
	&\quad + \f{\sigma_1}2  \f{1}{\v_1^2}\, |Y(U)|^2 \, {\bf 1}_{\{-\v_1^2\leq Y(U)\leq 0\}}  -\f12 \sigma_1 Y(U)\,{\bf 1}_{\{Y(U)\leq -\v_1^2\}} \nonumber\\
	&\leq -\f{1}{\v_1^4} |Y(U)|^2 \,  {\bf 1}_{\{0\leq Y(U)\leq \v_1^2\}}
	- \f{|\sigma_1|}{2}  \f{1}{\v_1^2}\, |Y(U)|^2 \, {\bf 1}_{\{-\v_1^2\leq Y(U)\leq 0\}} \nonumber\\
	&\quad  -\f12 |\sigma_1|\cdot |Y(U)|\,{\bf 1}_{\{Y(U)\leq -\v_1^2\}}=:\mathbb{G}_Y(U)\leq 0.\nonumber
\end{align}
Hence it holds
\begin{align}\label{6.4}
	\mathcal{F}(U)(\tau)\, {\bf 1}_{\{Y(U)\leq \v_1^2\}} &\leq \mathbb{G}_Y(U) + \Big(B_{\delta_1}(U)- G_{\delta_1}(U)\Big) {\bf 1}_{\{Y(U)\leq \v_1^2\}} \nonumber\\
	&= \mathbb{G}_Y(U) + \Big(B_{\delta_1}(U)- \big[(1-\v) \, G_1^+(U) + (1-\v) \, G_1^-(U) + G_2(v) \nonumber\\
	&\qquad\qquad\qquad\quad + G^r(v) + \v \, G^h(h) + D(v) + G_3(v)\big]\Big) {\bf 1}_{\{Y(U)\leq \v_1^2\}}.
\end{align}

To control the RHS of \eqref{6.4}, we need the following lemma.
\begin{lemma}\label{lem6.3}
	Let $\dis \f{\v_1}{\lambda}\ll \delta_0\ll \delta_1$ and $\v= \v_1\v_2$. If $Y(U)\leq \v_1^2$, then it holds
	\begin{align}\label{6.5}
		\mathbb{R}(U)&:=\mathbb{G}_Y(U) + {\bf 1}_{\{Y(U)\leq \v_1^2\}}\,\left(B_{\delta_1}(U) + \delta_0\f{\v_1}{\lambda} |B_{\delta_1}(U)|-\f{\v}{2} G^{h}(h) \right.\nonumber\\
		&\quad\left. -\f12 G_1^-(U) -\f12 G_1^+(U)-(1-\delta_0 \f{\v_1}{\lambda}) G_2(v)
		-(1-\delta_0) \Big[D(v)+ G^r(v)\Big]\right)\nonumber\\
		&\leq  3g(\nu,a,\v_1,\v_2,\tau) + 2f(\nu,a,\v_1,\v_2,\tau) \int_{\R} \eta(U|\tilde{U}) dy.
	\end{align}
\end{lemma}

\noindent{\bf Proof.} We divide the proof into four steps.

\noindent{\it Step 1.} We consider the case $\dis D(v)+G^r(v)\geq \f{2C_1^\ast}{\sqrt{\delta_0}} \f{\v_1^2}{\lambda}$, then we have from \eqref{5.11} that
\begin{align}
	\mathbb{R}(U)&\leq  2 |B_{\delta_1}(U)|\, {\bf 1}_{\{Y(U)\leq \v_1^2\}}
	-(1-\delta_0) \Big[D(v)+ G^r(v)\Big]\, {\bf 1}_{\{Y(U)\leq \v_1^2\}} \nonumber\\
	&\leq -(1-C\delta_0) \Big[D(v)+ G^r(v)\Big]\, {\bf 1}_{\{Y(U)\leq \v_1^2\}}\nonumber\\
	&\quad   +8 C_1^\ast \f{\v_1^2}{\lambda}  + 2f(\nu,a,\v_1,\v_2,\tau) \int_{\R} \eta(U|\tilde{U}) dy + 2g(\nu,a,\v_1,\v_2,\tau)\nonumber\\
	&\leq 2f(\nu,a,\v_1,\v_2,\tau) \int_{\R} \eta(U|\tilde{U}) dy + 2g(\nu,a,\v_1,\v_2,\tau).\nonumber
\end{align}

\smallskip

\noindent{\it Step 2.} Now we consider the other case
\begin{equation*}
	\dis D(v)+G^r(v)\leq \f{2C_1^\ast}{\sqrt{\delta_0}}\cdot \f{\v_1^2}{\lambda}.
\end{equation*}
Under the assumption $Y(U)\leq \v_1^2$, one obtains from Lemma \ref{lem3.4} that
\begin{align}\label{6.8}
	|Y(U)|\leq \tilde{C}_1 \f{\v_1^2}{\lambda}.
\end{align}
Hence it holds that
\begin{align}\label{6.9}
	-\f{1}{2} |\sigma_1|\cdot |Y(U)| \, {\bf 1}_{\{Y(U)\leq -\v_1^2\}}&= -\f{1}{2} |\sigma_1|\cdot |Y(U)| \, {\bf 1}_{\big\{-\tilde{C}_1 \f{\v_1^2}{\lambda}\leq Y(U)\leq -\v_1^2\big\}} \nonumber\\
	&\leq -\f{1}{2} |\sigma_1|\cdot |Y(U)|^{2} \f{1}{|Y(U)|} \, {\bf 1}_{\big\{-\tilde{C}_1 \f{\v_1^2}{\lambda}\leq Y(U)\leq -\v_1^2\big\}} \nonumber\\
	&\leq -\f{|\sigma_1|\lambda}{2\tilde{C}_1 \v_1^2}\cdot |Y(U)|^{2}  \, {\bf 1}_{\big\{-\tilde{C}_1 \f{\v_1^2}{\lambda}\leq Y(U)\leq -\v_1^2\big\}}.
\end{align}

Noting $\dis \f{\v_1}{\lambda}\ll \delta_0\ll \delta_1\ll 1$, then it follows from \eqref{6.8}-\eqref{6.9} that
\begin{align*}
	\mathbb{G}_Y(U)&\leq -\f{1}{\v_1^4} |Y(U)|^2 \,  {\bf 1}_{\{0\leq Y(U)\leq \v_1^2\}}
	- \f{|\sigma_1|}{2}  \f{1}{\v_1^2}\, |Y(U)|^2 \, {\bf 1}_{\{-\v_1^2\leq Y(U)\leq 0\}} \nonumber\\
	&\quad  -\f{|\sigma_1|\lambda}{2\tilde{C}_1 \v_1^2}\cdot |Y(U)|^{2}  \, {\bf 1}_{\big\{-\tilde{C}_1 \f{\v_1^2}{\lambda}\leq Y(U)\leq -\v_1^2\big\}}\nonumber\\
	&\leq -\f{|\sigma_1|\lambda}{2\tilde{C}_1 \v_1^2}\cdot |Y(U)|^{2}  \, {\bf 1}_{\big\{-\tilde{C}_1 \f{\v_1^2}{\lambda}\leq Y(U)\leq \v_1^2\big\}}\nonumber\\
	&\leq -\f{4}{\delta_1 \v_1}\cdot |Y(U)|^{2}  \, {\bf 1}_{\big\{-\tilde{C}_1 \f{\v_1^2}{\lambda}\leq Y(U)\leq \v_1^2\big\}},
\end{align*}
which yields
\begin{align}\label{6.10}
	\mathbb{R}(U)&\leq -\f{4}{\delta_1 \v_1}\cdot |Y(U)|^{2}  \, {\bf 1}_{\big\{-\tilde{C}_1 \f{\v_1^2}{\lambda}\leq Y(U)\leq \v_1^2\big\}} + \bigg(B_{\delta_1}(U) + \delta_0\f{\v_1}{\lambda} |B_{\delta_1}(U)| \nonumber\\
	&\quad  -\f{\v}{2} G^{h}(h) -\f12 G_1^-(U) -\f12 G_1^+(U)-(1-\delta_0 \f{\v_1}{\lambda}) G_2(v)\nonumber\\
	&\quad-(1-\delta_0) \Big[D(v)+ G^r(v)\Big]\bigg)\, {\bf 1}_{\big\{-\tilde{C}_1 \f{\v_1^2}{\lambda}\leq Y(U)\leq \v_1^2\big\}}.
\end{align}

Noting \eqref{5.2-1}-\eqref{5.2-2} and \eqref{3.36}, we have
\begin{align}
	Y(U)=[\mathbb{Y}^g(v)-\mathbb{Y}^g(\bar{v})] + [\mathbb{Y}^g(\bar{v})-\mathcal{Y}_g(\bar{v})] + \mathcal{Y}_g(\bar{v}) + \mathbb{Y}^b(U) + \mathbb{Y}^l(U) + \mathbb{Y}^s(U),\nonumber
\end{align}
which implies
\begin{align}\label{6.12}
	-4 |Y(U)|^2&\leq - |\mathcal{Y}_g(\bar{v})|^2 + 4\Big\{|\mathbb{Y}^g(v)-\mathbb{Y}^g(\bar{v})|^2 + |\mathbb{Y}^g(\bar{v})-\mathcal{Y}_g(\bar{v})|^2\nonumber\\
	&\qquad\qquad\qquad\qquad + |\mathbb{Y}^b(U)|^2 + |\mathbb{Y}^l(U)|^2 + |\mathbb{Y}^s(U)|^2\Big\}.
\end{align}
Using \eqref{6.12}, one obtains
\begin{align}\label{6.13}
	-\f{4}{\delta_1 \v_1}\cdot |Y(U)|^{2}  \, {\bf 1}_{\big\{-\tilde{C}_1 \f{\v_1^2}{\lambda}\leq Y(U)\leq \v_1^2\big\}}
	&\leq  \Big(- \f{1}{\delta_1 \v_1}|\mathcal{Y}_g(\bar{v})|^2 + \hat{\mathbb{Y}}\Big)\, {\bf 1}_{\big\{-\tilde{C}_1 \f{\v_1^2}{\lambda}\leq Y(U)\leq \v_1^2\big\}},
\end{align}
where we have denoted
\begin{align}\label{6.14}
	\hat{\mathbb{Y}}&:= \f{4}{\delta_1 \v_1}\Big\{|\mathbb{Y}^g(v)-\mathbb{Y}^g(\bar{v})|^2  + |\mathbb{Y}^b(U)|^2 + |\mathbb{Y}^l(U)|^2 + |\mathbb{Y}^s(U)|^2+ |\mathbb{Y}^g(\bar{v})-\mathcal{Y}_g(\bar{v})|^2\Big\}.
\end{align}

\medskip

Noting \eqref{3.29} and \eqref{3.40}, we have
\begin{align}
	B_{\delta_1}(U)&= B_1(\bar{v}) + \frac{1}{1-\v} \,B_{2}^+(\bar{v}) + [B_1(v)-B_1(\bar{v})]+ \frac{1}{1-\v} \, [B_{2}^+(v)-B_{2}^+(\bar{v})] \nonumber\\
	&\quad+ B_{2}^-(U) +B_3(v) +B_4(v) + B_5(U)\nonumber\\
	&=\mathcal{I}_1(\bar{v})+ \mathcal{I}_2(\bar{v}) + [B_1(\bar{v})-\mathcal{I}_1(\bar{v})] +  [\frac{1}{1-\v} \,B_{2}^+(\bar{v})-\mathcal{I}_2(\bar{v})]\nonumber\\
	&\quad + [B_1(v)-B_1(\bar{v})]+ \frac{1}{1-\v} \, [B_{2}^+(v)-B_{2}^+(\bar{v})]\nonumber\\
	&\quad+ B_{2}^-(U) +B_3(v) +B_4(v) + B_5(U),\nonumber
\end{align}
which yields that
\begin{align}\label{6.16}
	&B_{\delta_1}(U) + \delta_0\f{\v_1}{\lambda} |B_{\delta_1}(U)|\leq \mathcal{I}_1(\bar{v})+ \mathcal{I}_2(\bar{v}) + \delta_0\f{\v_1}{\lambda}\big(|\mathcal{I}_1(\bar{v})|+|\mathcal{I}_2(\bar{v})|\big)
	+ (1+\delta_0\f{\v_1}{\lambda})\, \hat{\mathbb{B}},
\end{align}
with
\begin{align}
	\hat{\mathbb{B}}:=&|B_1(\bar{v})-\mathcal{I}_1(\bar{v})| + \big|\frac{1}{1-\v} \,B_{2}^+(\bar{v})-\mathcal{I}_2(\bar{v})\big|\nonumber\\
	& + |B_1(v)-B_1(\bar{v})|+ \frac{1}{1-\v} \, |B_{2}^+(v)-B_{2}^+(\bar{v})|\nonumber\\
	&+ |B_{2}^-(U)| +|B_3(v)| +|B_4(v)| + |B_5(U)|.\nonumber
\end{align}
Substituting \eqref{6.16} and \eqref{6.13} into \eqref{6.10}, one gets
\begin{align}
	\mathbb{R}(U)&\leq  - \f{1}{\delta_1 \v_1}|\mathcal{Y}_g(\bar{v})|^2 \,  {\bf 1}_{\big\{-\tilde{C}_1 \f{\v_1^2}{\lambda}\leq Y(U)\leq \v_1^2\big\}} \nonumber\\
	&\quad   + \Big(\mathcal{I}_1(\bar{v})+ \mathcal{I}_2(\bar{v}) + \delta_0\f{\v_1}{\lambda}\big(|\mathcal{I}_1(\bar{v})|+|\mathcal{I}_2(\bar{v})|\big) \Big) \,{\bf 1}_{\big\{-\tilde{C}_1 \f{\v_1^2}{\lambda}\leq Y(U)\leq \v_1^2\big\}} \nonumber\\
	&\quad  -\bigg(\f{\v}{2} G^{h}(h) +\f12 G_1^-(U) +\f12 G_1^+(U)+(1-\delta_0 \f{\v_1}{\lambda}) G_2(v)\nonumber\\
	&\qquad\qquad+(1-\delta_0) \Big[D(v)+ G^r(v)\Big]\bigg){\bf 1}_{\big\{-\tilde{C}_1 \f{\v_1^2}{\lambda}\leq Y(U)\leq \v_1^2\big\}}\nonumber\\
	&\quad + \big[\hat{\mathbb{Y}}  + 2 \hat{\mathbb{B}} \big] \, {\bf 1}_{\big\{-\tilde{C}_1 \f{\v_1^2}{\lambda}\leq Y(U)\leq \v_1^2\big\}}.\nonumber
\end{align}

Noting  $\dis \f{\v_1}{\lambda}\leq \delta_0\ll \delta_1$ and Lemmas \ref{lem3.4-1} and \ref{lem3.5}, we can obtain
\begin{align}\label{6.19}
	\mathbb{R}(U)&\leq  \bigg\{- \f{1}{\delta_1 \v_1}|\mathcal{Y}_g(\bar{v})|^2 + \mathcal{I}_1(\bar{v}) +  \delta_1 |\mathcal{I}_1(\bar{v})| + \mathcal{I}_2(\bar{v}) + \delta_1 \f{\v_1}{\lambda} |\mathcal{I}_2(\bar{v})|\nonumber\\
	&\quad  -(1-\delta_1\f{\v_1}{\lambda}) \, G_2(\bar{v}) - (1-\delta_1) D(\bar{v})
	\bigg\}  {\bf 1}_{\big\{-\tilde{C}_1 \f{\v_1^2}{\lambda}\leq Y(U)\leq \v_1^2\big\}} \nonumber\\
	&\quad  -\bigg(\f{\v}{2} G^{h}(h) +\f12 G_1^-(U) +\f12 G_1^+(U)+(1-\delta_1 \f{\v_1}{\lambda}) [G_2(v) - G_2(\bar{v})]\nonumber\\
	&\quad + (\delta_1-\delta_0) \f{\v_1}{\lambda}  G_2(v) + (\delta_1-\delta_0) D(v)  +(1-\delta_0)  G^r(v) \bigg){\bf 1}_{\big\{-\tilde{C}_1 \f{\v_1^2}{\lambda}\leq Y(U)\leq \v_1^2\big\}}\nonumber\\
	&\quad + \big[\hat{\mathbb{Y}}  + 2 \hat{\mathbb{B}} \big] \, {\bf 1}_{\big\{-\tilde{C}_1 \f{\v_1^2}{\lambda}\leq Y(U)\leq \v_1^2\big\}}\nonumber\\
	&\leq -\bigg(\f{\v}{2} G^{h}(h) +\f12 G_1^-(U) +\f12 G_1^+(U)+(1-\delta_1 \f{\v_1}{\lambda}) [G_2(v) - G_2(\bar{v})]\nonumber\\
	&\quad + (\delta_1-\delta_0) \f{\v_1}{\lambda}  G_2(v) + (\delta_1-\delta_0) D(v) +(1-\delta_0)  G^r(v) \bigg){\bf 1}_{\big\{-\tilde{C}_1 \f{\v_1^2}{\lambda}\leq Y(U)\leq \v_1^2\big\}}\nonumber\\
	&\quad + \big[\hat{\mathbb{Y}}  + 2 \hat{\mathbb{B}} \big] \, {\bf 1}_{\big\{-\tilde{C}_1 \f{\v_1^2}{\lambda}\leq Y(U)\leq \v_1^2\big\}}.
\end{align}

\noindent{\it Step 3. Estimate on $\hat{\mathbb{Y}}$.} Noting \eqref{3.36} and \eqref{5.2-3}, one has that
\begin{align*}
	|\mathbb{Y}^g(\bar{v})-\mathcal{Y}_g(\bar{v})|
	&\leq C\int_{\Omega^c} |\partial_yw_X|\cdot |p(\bar{v})-p(\tilde{v})|^2dy
	+ C\int_{\Omega^c} |\partial_yw_X|\, Q(\bar{v}|\tilde{v})dy\nonumber\\
	&\quad + C\f{\v_1}{\lambda}\int_{\Omega^c} |\partial_yw_X\,[p(\bar{v})-p(\tilde{v})]| dy\nonumber\\
	&\leq C\f{\v_1}{\lambda} \Big( \int_{\Omega^c} |\partial_yw_X| dy\Big)^{\f12} \Big( \int_{\Omega^c} |\partial_yw_X|\cdot |p(\bar{v})-p(\tilde{v})|^2 dy\Big)^{\f12}\nonumber\\
	&\quad + C\int_{\Omega^c} |\partial_yw_X|\, Q(\bar{v}|\tilde{v})dy\nonumber\\
	&\leq C\f{\v_1}{\sqrt{\lambda}} \Big(\int_{\Omega^c} |\partial_yw_X|\, Q(\bar{v}|\tilde{v})dy\Big)^{\f12},
\end{align*}
which, together with \eqref{4.17} and \eqref{5.16}, yields that
\begin{align}\label{6.20}
	|\mathbb{Y}^g(\bar{v})-\mathcal{Y}_g(\bar{v})|^2
	& \leq C\f{\v_1^2}{\lambda }\int_{\Omega^c} |\partial_yw_X|\, Q(\bar{v}|\tilde{v})dy\leq C\f{\v_1^2}{\lambda} \sqrt{\f{\v_1}{\lambda}}
	\Big[D(v) + \delta_0 \f{\v_1}{\lambda} G_2(v) + \delta_0 G^r(v)\Big].
\end{align}
Combining \eqref{6.14}, \eqref{6.20} and \eqref{5.13}, one obtains
\begin{align}\label{6.21}
	\hat{\mathbb{Y}}
	&\leq C\f{\v_1}{\delta_1\lambda}  \bigg\{ \delta_{0}  \Big[D(v) + G^r(v)\Big]+\Big[G_{2}(v)-G_{2}(\bar{v})\Big] +\left(\frac{\v_1}{\lambda}\right)^{\f14} G_{2}(\bar{v})\nonumber \\
	&\quad+G_{1}^{-}(U)+\left(\frac{\lambda}{\v_1}\right)^{\f14} G_{1}^{+}(v)\bigg\} + C\f{\v_1\v_2^2}{\delta_1} \Big[ \exp\big\{-\f{\nu}{C_\ast a}\tau\big\}
	+ \exp\big\{-\f{\v_1}{C_\ast}\tau\big\}\Big]\nonumber\\
	&\leq C\f{\v_1}{\delta_1\lambda}  \bigg\{ \delta_{0}  \Big[D(v) + G^r(v)\Big]+\Big[G_{2}(v)-G_{2}(\bar{v})\Big] +\left(\frac{\v_1}{\lambda}\right)^{\f14} G_{2}(\bar{v})\nonumber \\
	&\quad+G_{1}^{-}(U)+\left(\frac{\lambda}{\v_1}\right)^{\f14} G_{1}^{+}(U)\bigg\} +g(\nu,a,\v_1,\v_2,\tau).
\end{align}

\noindent{\it Step 4. Estimate on $\hat{\mathbb{B}}$.} From the definitions of $B_1(\bar{v}), B_2^+(\bar{v})$ in \eqref{3.29}, and $\mathcal{I}_1(\bar{v})$ and $\mathcal{I}_2(\bar{v})$ in \eqref{3.40}, we get
\begin{align}\label{6.22}
	|B_1(\bar{v})-\mathcal{I}_1(\bar{v})|=0,
\end{align}
and
\begin{align}\label{6.23}
	\big|\frac{1}{1-\v} \,B_{2}^+(\bar{v})-\mathcal{I}_2(\bar{v})\big|
	&=\f{1}{2|\sigma_1|} \left|\f{1}{1-\v}\int_{\Omega}  \partial_yw_X \, [p(\bar{v})-p(\tilde{v})]^2 dy - \int_{\R}  \partial_yw_X \, [p(\bar{v})-p(\tilde{v})]^2 dy\right|\nonumber\\
	&\leq C\v \int_{\R}  |\partial_yw_X| \,  |p(\bar{v})-p(\tilde{v})|^2 dy + C\int_{\Omega^c}  |\partial_yw_X| \, |p(\bar{v})-p(\tilde{v})|^2 dy\nonumber\\
	&\leq C\v \int_{\R}  |\partial_yw_X| \,  Q(\bar{v}|\tilde{v}) dy + C\int_{\Omega^c}|\partial_yw_X| Q(\bar{v}|\tilde{v}) dy\nonumber\\
	&\leq C\v G_2(\bar{v}) + C\sqrt{\f{\v_1}{\lambda}}
	\Big[D(v) + \delta_0 \f{\v_1}{\lambda} G_2(v) + \delta_0 G^r(v)\Big].
\end{align}
Combining \eqref{6.22}-\eqref{6.23}, \eqref{5.1}, \eqref{5.3} and \eqref{5.10}, one has
\begin{align}\label{6.24}
	\hat{\mathbb{B}}&\leq  C\v G_2(\bar{v}) + C\sqrt{\f{\v_1}{\lambda}}
	\Big[D(v) + \delta_0 \f{\v_1}{\lambda} G_2(v) + \delta_0 G^r(v)\Big] + \left(\delta_{1}+C \delta_{0}\right) G_{1}^{-}(U)\nonumber\\
	&\quad + C  \delta_{0}\Big(D(v)+[G_{2}(v)-G_{2}(\bar{v})] +\f{\v_1}{\lambda}G_{2}(v) + G^r(v)\Big) + \delta_0 \v G^h(h)\nonumber\\
	&\quad + C \delta_0  D(v)  + f(\nu,a,\v_1,\v_2,\tau) \int_{\R} \eta(U|\tilde{U}) dy + g(\nu,a,\v_1,\v_2,\tau) \nonumber\\
	&\leq C \delta_{0} \Big[ D(v)+G^r(v)\Big]+ C \delta_0 [G_{2}(v)-G_{2}(\bar{v})] +C\delta_0\f{\v_1}{\lambda}G_{2}(v) + \delta_0 \v G^h(h) \nonumber\\
	&\quad +\left(\delta_{1}+C \delta_{0}\right) G_{1}^{-}(U) + f(\nu,a,\v_1,\v_2,\tau) \int_{\R} \eta(U|\tilde{U}) dy+ g(\nu,a,\v_1,\v_2,\tau).
\end{align}

Substituting \eqref{6.24} and \eqref{6.21} into \eqref{6.19}, one obtains
\begin{align}
	\mathbb{R}(U)
	&\leq -\bigg( [\f12-C\delta_0] \v G^{h}(h) + [\f12-\delta_1-C\delta_0] G_1^-(U) +[\f12-\delta_0] G_1^+(U)\nonumber\\
	&\quad  +(1-\delta_1 \f{\v_1}{\lambda}-C\delta_0) [G_2(v) - G_2(\bar{v})] + (\delta_1-C\delta_0) \f{\v_1}{\lambda}  G_2(v) \nonumber\\
	&\quad + (\delta_1-C \delta_0) D(v)  +(1-C \delta_0)  G^r(v) \bigg){\bf 1}_{\big\{-\tilde{C}_1 \f{\v_1^2}{\lambda}\leq Y(U)\leq \v_1^2\big\}}\nonumber\\
	&\quad + 2 f(\nu,a,\v_1,\v_2,\tau) \int_{\R} \eta(U|\tilde{U}) dy + 3 g(\nu,a,\v_1,\v_2,\tau)\nonumber\\
	&\leq  2f(\nu,a,\v_1,\v_2,\tau) \int_{\R} \eta(U|\tilde{U}) dy + 3 g(\nu,a,\v_1,\v_2,\tau),\nonumber
\end{align}
which yields \eqref{6.5}. Therefore the proof of Lemma \ref{lem6.3} is complete. $\hfill\Box$

\begin{lemma}\label{lem6.2}
	It holds that
	\begin{align}\label{6.26}
		&\f{d}{d\tau} \int_{\R} w_X\,\eta(U|\tilde{U}) dy + \Big(1+ |J^{bad}| + J^{good}\Big) \, {\bf 1}_{\{Y(U)\geq \v_1^2\}} \nonumber\\
		&\quad +\delta_0\f{\v_1}{\lambda} |B_{\delta_1}(U)|  \,  {\bf 1}_{\{Y(U)\leq \v_1^2\}}
		+ \bigg(\f14 \, G_1^+(U) + \f14 \, G_1^-(U) + \delta_0 \f{\v_1}{\lambda} G_2(v) \nonumber\\
		&\quad   + \f12\v  \, G^h(h) + \delta_0D(v) + \delta_0 G^r(v) + G_3(v)\bigg) {\bf 1}_{\{Y(U)\leq \v_1^2\}}\nonumber\\
		&\leq  3 g(\nu,a,\v_1,\v_2,\tau) +  2f(\nu,a,\v_1,\v_2,\tau) \int_{\R} \eta(U|\tilde{U}) dy.
	\end{align}
\end{lemma}

\noindent{\bf Proof.} It follows from \eqref{6.4} and \eqref{6.5} that
\begin{align}
	&\mathcal{F}(U)(\tau)\, {\bf 1}_{\{Y(U)\leq \v_1^2\}} \nonumber\\
	&\leq \mathbb{G}_Y(U) + \Big(B_{\delta_1}(U)- \big[(1-\v) \, G_1^+(U) + (1-\v) \, G_1^-(U) + G_2(v) \nonumber\\
	&\qquad\qquad\qquad\quad + G^r(v) + \v \, G^h(h) + D(v) + G_3(v)\big]\Big) {\bf 1}_{\{Y(U)\leq \v_1^2\}}\nonumber\\
	&\leq \mathbb{R}(U) -\delta_0\f{\v_1}{\lambda} |B_{\delta_1}(U)|  \,  {\bf 1}_{\{Y(U)\leq \v_1^2\}}
	- \bigg(\f14 \, G_1^+(U) + \f14 \, G_1^-(U) + \delta_0 \f{\v_1}{\lambda} G_2(v) \nonumber\\
	&\qquad\qquad  + \f12\v \, G^h(h) + \delta_0D(v) + \delta_0 G^r(v) + G_3(v)\bigg) {\bf 1}_{\{Y(U)\leq \v_1^2\}}\nonumber\\
	&\leq -\delta_0\f{\v_1}{\lambda} |B_{\delta_1}(U)|  \,  {\bf 1}_{\{Y(U)\leq \v_1^2\}}
	- \bigg(\f14 \, G_1^+(U) + \f14 \, G_1^-(U) + \delta_0 \f{\v_1}{\lambda} G_2(v) \nonumber\\
	&\quad   + \f12\v_1\v_2 \, G^h(h) + \delta_0D(v) + \delta_0 G^r(v) + G_3(v)\bigg) {\bf 1}_{\{Y(U)\leq \v_1^2\}}\nonumber\\
	&\quad +3g(\nu,a,\v_1,\v_2,\tau) + 2f(\nu,a,\v_1,\v_2,\tau) \int_{\R} \eta(U|\tilde{U}) dy,\nonumber
\end{align}
which, together with \eqref{6.2}, gives \eqref{6.26}. Therefore the proof of Lemma \ref{lem6.2} is complete. $\hfill\Box$

\


\noindent{\bf Proof of Theorem \ref{thm3.1}.} Noting \eqref{5.12}, we have
\begin{align}\label{6.29}
	&\int_{0}^\infty f(\nu,a,\v_1,\v_2,\tau) d\tau\nonumber\\
	&=C(v_-)\int_0^\infty \Big[ \f{1}{\delta_0} (\min\{\f{\nu\v_2}{a},\tau^{-1}\})^{\f43}
	+ \f{\nu}{a} \exp\{-\f{\nu}{C_\ast a}\tau\} + \v_1 \exp\{-\f{\v_1}{C_\ast}\tau\} \Big] d\tau\nonumber\\
	&\leq C(v_-) \Big[\f{1}{\delta_0} (\f{\nu\v_2}{a})^{\f13}+ 2C_\ast \Big]\leq C(v_-,\delta_{0}^{-1}),
\end{align}
and
\begin{align}\label{6.30}
	&\int_{0}^\infty g(\nu,a,\v_1,\v_2,\tau) d\tau\nonumber\\
	&\leq C(v_{-})\int_0^\infty  \Big[\f{1}{\delta_0} (\min\{\f{\nu\v_2}{a},\tau^{-1}\})^{\f43}
	+ \f{\v_1\v_2}{\delta_0\lambda}\Big(\exp\{-\f{\nu}{C_\ast a}\tau\} +   \exp\{-\f{\v_1}{C_\ast}\tau\}\Big) \Big]d\tau\nonumber\\
	&\leq C(v_{-}) \bigg\{\f{1}{\delta_0} (\f{\nu\v_2}{a})^{\f16}+ \f{\v_1\v_2}{\delta_0\lambda} [\f{C_\ast a}{\nu}+\f{C_\ast}{\v_1}]\bigg\}\nonumber\\
	&\leq C(U_-, \delta_0^{-1})\Big\{  (\f{\nu\v_2}{a})^{\f16}  + \v_2 \f{a}{\nu}\Big\}.
\end{align}
We point out that the global integrability \eqref{6.29} is uniform in $\nu$, which is crucial to get \eqref{1.20}. While, the estimate \eqref{6.30} is not uniform in $\nu$ since $a/\nu\to \infty$ as $\nu\to 0$. Fortunately, it is enough for us to take inviscid limit in original coordinate $(t,x)$.

Applying Gr\"{o}nwall's inequality to \eqref{6.26}, and using \eqref{6.29}-\eqref{6.30}, one obtains
\begin{align}\label{6.31}
	&\int_{\R} w_X\, \eta(U|\tilde{U})(\tau) dy  + \int_0^\infty \Big(1+ |J^{bad}(U)| + J^{good}(U)\Big) (\tau) \, {\bf 1}_{\{Y(U)\geq \v_1^2\}} \, d\tau \nonumber\\
	&\quad +\int_0^{\infty} \delta_0\f{\v_1}{\lambda} |B_{\delta_1}(U)(\tau)|  \,  {\bf 1}_{\{Y(U)\leq \v_1^2\}}\, d\tau
	+ \int_0^{\infty} \bigg(\f14 \, G_1^+(v) + \f14 \, G_1^-(U) \nonumber\\
	&\quad + \delta_0 \f{\v_1}{\lambda} G_2(v)   + \f12\v\, G^h(h) + \delta_0D(v) + \delta_0 G^r(v) + G_3(v)\bigg)(\tau) {\bf 1}_{\{Y(U)\leq \v_1^2\}} \, d\tau \nonumber\\
	&\leq C\Big[\int_{\R} w_X\, \eta(U|\tilde{U})(0) dy+ 3\int_{0}^\infty g(\nu,a,\v_1,\v_2,\tau) d\tau \Big] \exp\Big\{2\int_{0}^\infty f(\nu,a,\v_1,\v_2,\tau) d\tau \Big\}\nonumber\\
	&\leq \hat{C}(U_-,\delta_0^{-1}) \bigg[\int_{\R} w_X\, \eta(U|\tilde{U})(0) dy +   (\f{\nu\v_2}{a})^{\f16}  + \v_2 \f{a}{\nu}\bigg],
\end{align}
which yields \eqref{1.20}.

\smallskip

By similar arguments as in \cite{Kang-Vasseur-2021-Invent}, we have from \eqref{3.7-1} that
\begin{align*}
	|J^{bad}(U)|&\leq |J^{bad}(U)(\tau)| \, {\bf 1}_{\{Y(U)\geq \v_1^2\}} + |B_{\delta_1}(U)(\tau)|  \,  {\bf 1}_{\{Y(U)\leq \v_1^2\}}\nonumber\\
	&\quad + |J^{good}(U)-G_{\delta_1}(U)|\,  {\bf 1}_{\{Y(U)\leq \v_1^2\}}\nonumber\\
	&\leq  |J^{bad}(U)(\tau)| \, {\bf 1}_{\{Y(U)\geq \v_1^2\}} + |B_{\delta_1}(U)(\tau)|  \,  {\bf 1}_{\{Y(U)\leq \v_1^2\}}\nonumber\\
	&\quad + C\int_{\R} |\partial_y w_X|\, \big[|h-\tilde{h}|^2 + Q(v|\tilde{v})\big] dy\nonumber\\
	&\leq  |J^{bad}(U)(\tau)| \, {\bf 1}_{\{Y(U)\geq \v_1^2\}} + |B_{\delta_1}(U)(\tau)|  \,  {\bf 1}_{\{Y(U)\leq \v_1^2\}} \nonumber\\ &\quad + C(v_-) \Big[G^h(h) +  G_2(v)\Big],
\end{align*}
which, together with \eqref{5.7-1}, yields that
\begin{align}
	|\dot{X}(\tau)|&\leq \max\Big\{\f{1}{\v_1^2} (2|J^{bad}(U)| + 1),\, \f12 |\sigma_1|\Big\}\leq C(U_-, \v_1^{-1}) \, \Big\{ \mathfrak{f}(\tau) + 1\Big\},\nonumber
\end{align}
with
\begin{align*}
	\mathfrak{f}(\tau):&= |J^{bad}(U)| \, {\bf 1}_{\{Y(U)\geq \v_1^2\}} + |B_{\delta_1}(U)|  \,  {\bf 1}_{\{Y(U)\leq \v_1^2\}}+ C(v_-) \Big[G^h(h) +  G_2(v)\Big].
\end{align*}
It follows from \eqref{6.31} that
\begin{align}
	\int_0^\infty \mathfrak{f}(\tau) d\tau \leq \hat{C}(U_-, \v_1^{-1}, \v_2^{-1}) \, \bigg[\int_{\R}  \eta(U|\tilde{U})(0) dy +   (\f{\nu\v_2}{a})^{\f16}  + \v_2 \f{a}{\nu}\bigg].\nonumber
\end{align}
Hence we conclude \eqref{1.20-3}-\eqref{1.20-4}. Therefore the proof of Theorem \ref{thm3.1} is complete. $\hfill\Box$

\

\section{Proof of the main theorem: Theorem \ref{thm1.1}}

Now we go back to the original coordinate $(t,x)$ to prove Theorem \ref{thm1.1}. By using similar arguments as in \cite[Section 5.1]{Kang-Vasseur-2021-Invent}, and noting \eqref{1.17-1} and Lemmas \ref{lem2.3}-\ref{lem2.2}, we can construct a sequence of smooth approximate initial datum $(\mathfrak{v}_0^\nu, \mathfrak{u}_0^\nu)$ satisfying \eqref{1.17-2}, the details of construction are omitted here for simplicity of presentation. We divide the rest proof into several subsections.

\smallskip

\subsection{Uniform estimate.} Recall the definitions of $(\tilde{\mathfrak{v}}^s,\tilde{\mathfrak{u}}^s)$ and $(\tilde{v}^r, \tilde{u}^r)$ in \eqref{2.2} and \eqref{2.5} respectively. We define
\begin{align}\label{7.1-1}
	\begin{aligned}
		\tilde{\mathfrak{v}}^\nu_X(t,x)&:=\tilde{\mathfrak{v}}_X^s(t,x)+ \tilde{v}^r(t,x)-v_m,\\
		\tilde{\mathfrak{u}}^\nu_X(t,x)&:=\tilde{\mathfrak{u}}_X^s(t,x)+ \tilde{u}^r(t,x)-u_m,\\
		\tilde{\mathfrak{h}}_{X}^{\nu}(t,x)&:=\tilde{\mathfrak{u}}_{X}^{\nu}+\nu \frac{\gamma}{\alpha}\big(p(\tilde{\mathfrak{v}}_{X}^{s})^{\frac{\alpha}{\gamma}}\big)_{x},
	\end{aligned}
\end{align}
where
\begin{align}\label{7.1}
	X_{\nu}(t):= \nu X_{\nu}(\tau) \equiv\nu \, X(\f{t}{\nu}),
\end{align}
and
\begin{align}\label{7.1-2}
	\begin{split}
		\dis \tilde{\mathfrak{v}}_X^s(t,x)&=\tilde{\mathfrak{v}}^s(t,x-X_\nu(t))\equiv\tilde{v}^s\left(\f{x-\sigma_1t-X_\nu(t)}{\nu}\right),\\
		\dis \tilde{\mathfrak{u}}_X^s(t,x)&=\tilde{\mathfrak{u}}^s(t,x-X_\nu(t))\equiv\tilde{u}^s\left(\f{x-\sigma_1t-X_\nu(t)}{\nu}\right).
	\end{split}
\end{align}
Recall \eqref{1.18}, one has
\begin{align}
	\begin{aligned}
	\eta\big((\mathfrak{v}^{\nu},\mathfrak{h}^{\nu})\,|\, (\tilde{\mathfrak{v}}_{X}^{\nu},\tilde{\mathfrak{h}}_{X}^{\nu})\big)&=\frac{1}{2}|\mathfrak{h}^{\nu}-\tilde{\mathfrak{h}}_{X}^{\nu}|^2+Q(\mathfrak{v}^{\nu}|\tilde{\mathfrak{v}}_{X}^{\nu})\\
	&\equiv \frac{1}{2}|\mathfrak{u}^{\nu}+\nu\frac{\gamma}{\alpha}(p(\mathfrak{v}^{\nu})^{\frac{\alpha}{\gamma}})_{x}-\tilde{\mathfrak{u}}_{X}^{\nu}-\nu\frac{\gamma}{\alpha}(p(\mathfrak{v}_{X}^{s})^{\frac{\alpha}{\gamma}})_{x}|^2+Q(\mathfrak{v}^{\nu}|\tilde{\mathfrak{v}}_{X}^{\nu}).
	\end{aligned}\nonumber
\end{align}
Then, for the solution $(\mathfrak{v}^\nu,\mathfrak{v}^\nu)(t,x)$ of compressible Navier-Stokes equations \eqref{1.1},  it follows from Theorem \ref{thm3.1} that
\begin{align}\label{7.13}
	&\sup_{t\geq0}\int_{\R} \eta\big((\mathfrak{v}^{\nu},\mathfrak{h}^{\nu})\,|\, (\tilde{\mathfrak{v}}_{X}^{\nu},\tilde{\mathfrak{h}}_{X}^{\nu})\big)\, dx
	+\int_0^\infty \int_{\R} |\partial_x\tilde{u}^r(t,x)|\,\, p(\mathfrak{v}^\nu\,|\,\tilde{\mathfrak{v}}_X^\nu)(t,x)\, dxdt \nonumber\\
	&\quad + \int_0^\infty \int_{\R} |\partial_x \tilde{\mathfrak{v}}_X^s(t,x)|\cdot Q(\mathfrak{v}^\nu\,|\, \tilde{\mathfrak{v}}_X^\nu)(t,x)\, dxdt\nonumber\\
	&\quad + \nu \int_0^\infty \int_{\R} (\mathfrak{v}^\nu)^{\gamma-\alpha}\, |\partial_x(p(\mathfrak{v}^\nu)-p(\tilde{\mathfrak{v}}_X^\nu))(t,x)|^2\, dxdt\nonumber\\
	&\leq \hat{C}(U_-, \v_1^{-1}, \v_2^{-1}) \,\left\{ \int_{\R} \eta\big((\mathfrak{v}^{\nu},\mathfrak{h}^{\nu})\,|\, (\tilde{\mathfrak{v}}_{X}^{\nu},\tilde{\mathfrak{h}}_{X}^{\nu})\big)(0,x)\, dx
	+ \nu\,(\f{\nu\v_2}{a})^{\f16}  + \v_2\, a(\nu) \right\}\nonumber\\
	&\leq \hat{C}(U_-, \v_1^{-1}, \v_2^{-1})  \Big\{\mathcal{E}_0 + \nu\,(\f{\nu\v_2}{a})^{\f16}  + \v_2\, a(\nu)\Big\}.
\end{align}

\subsection{Convergence of $X_\nu$ as $\nu\to 0+$.} It follows from \eqref{1.20-3} and \eqref{7.1} that
\begin{align}
	\left|\f{d}{dt} X_{\nu}(t)\right|= |X'(\tau)|
	\leq  C(U_-, \v_1^{-1}) \, \Big\{ \mathfrak{f}(\tau) + 1\Big\},\nonumber
\end{align}
which, together with \eqref{1.20-4}, yields that
\begin{align}
	&\int_0^T\left|\f{d}{dt} X_{\nu}(t)\right|dt
	\leq   C(U_-, \v_1^{-1}) \, \Big\{ \int_0^T\mathfrak{f}(\tau) dt+ T\Big\}\nonumber\\
	&\leq C(U_-, \v_1^{-1}) \, \Big\{ \nu\int_0^\infty\mathfrak{f}(\tau) d\tau+ T\Big\}\nonumber\\
	&\leq \hat{C}(U_-, \v_1^{-1}, \v_2^{-1}) \, \bigg[\nu\int_{\R}  \eta(U|\tilde{U})(0) dy +  \nu (\f{\nu\v_2}{a})^{\f16}  +\nu \v_2 \f{a}{\nu} + T \bigg]\nonumber\\
	&\leq \hat{C}(U_-, \v_1^{-1}, \v_2^{-1}) \, \bigg[ \int_{\R}  \eta(U^{\nu}|\tilde{U}^{\nu})(0) dx +  \nu   + a(\nu)\,\v_2   + T \bigg],\nonumber
\end{align}
and
\begin{align}\label{7.5}
	|X_{\nu}(t)|\leq  \hat{C}(U_-, \v_1^{-1}, \v_2^{-1}) \, \bigg[ \int_{\R}  \eta(U^{\nu}|\tilde{U}^{\nu})(0) dx +  \nu   + a(\nu)\,\v_2   + t \bigg].
\end{align}
Since   $\dis \f{d}{dt} X_{\nu} \in L^1(0,T)$ and  $ X_{\nu} \in L^1(0,T)$, by compactness of BV space, there exits a  function $X_\infty(\cdot) \in BV(0,T)$ such that
\begin{align}\label{7.6}
	X_{\nu}\to X_\infty\quad \mbox{in}\,\, L^{1}(0,T),\quad \mbox{up to a subsequence as }\,\, \nu\to 0+.
\end{align}
It follows from \eqref{7.5} that $X_{\infty}(t)\in BV(0,T)$ and
\begin{align}
	|X_\infty (t)|\leq	\hat{C}(U_-, \v_1^{-1}, \v_2^{-1}) \, (\mathcal{E}_{0}+t)\qquad \text{for a.e. } t\in [0,T].\label{7.10-1}
\end{align}

\smallskip

Using \eqref{1.20-2}, one obtains
\begin{align}
	X_{\nu}(t)=\nu \, X(\f{t}{\nu})\, \leq -\f12 \sigma_1t\qquad \text{for } t\in [0,T],\nonumber
\end{align}
which, together with \eqref{7.6}, yields that
\begin{align}\label{7.12}
	X_\infty(t)\leq -\f{1}{2}\sigma_1 \, t\qquad \text{for a.e. }t\in [0,T].
\end{align}
Hence we have
\begin{align}\label{7.12-1}
	\sigma_1 t + X_\infty(t)\leq \f{1}{2}\sigma_1 \, t\leq 0\qquad \text{for a.e. }t\in [0,T].
\end{align}
That means the shifted 1-shock wave and 2-rarefaction wave are still detached.

\subsection{Weak Convergence of $(\mathfrak{v}^\nu, \mathfrak{u}^{\nu})$ as $\nu\to 0+$: Proof of \eqref{1.17-3}.}
 Even though the influence of rarefaction wave should be considered, the similar arguments in \cite[Section 5.2.2]{Kang-Vasseur-2021-Invent} are still valid. Indeed, using similar arguments as in \cite[Section 5.2.2]{Kang-Vasseur-2021-Invent} and \eqref{7.13}, $\eqref{A.1}_2$ (we omit the details of proof here), there exist limits $$\mathfrak{v}_\infty\in L^\infty(0,T; L^\infty(\R) + \mathcal{M}(\R)) \quad \mbox{and}\quad  \mathfrak{u}_\infty\in L^\infty(0,T; L^2_{\rm loc}(\R))$$ such that
\begin{align}\label{7.14}
	\begin{split}
		&\mathfrak{v}^\nu \rightharpoonup \mathfrak{v}_\infty \quad \mbox{in}\,\, \mathcal{M}_{\rm loc}((0,T)\times\R),\\
		&\mathfrak{h}^\nu \rightharpoonup \mathfrak{u}_\infty \quad \mbox{in}\,\, L^\infty(0,T; L^2_{\rm loc}(\R)),\\
		&\mathfrak{u}^\nu \rightharpoonup \mathfrak{u}_\infty \quad \mbox{in}\,\, \mathcal{M}_{\rm loc}((0,T)\times\R).
	\end{split}
\end{align}

\subsection{Proof of \eqref{1.21}.} We denote
$$
M=3\max\{v_{-},v_{+}^{-1}\},\quad L=2\max\{\|\tilde{\mathfrak{h}}_{X}^{\nu}\|_{L^{\infty}}, |u_{+}|, |u_{-}|\},
$$
then it is clear that $\tilde{\mathfrak{v}}_{X}^{\nu}\in [2{M}^{-1},1/2M]$ and $\tilde{\mathfrak{h}}_{X}^{\nu}\in [-1/2L,1/2L]$. For later use, we define
\begin{equation}\label{7.15}
	\underline{\mathfrak{v}}^{\nu}=\psi_1(\mathfrak{v}^{\nu}),\quad \text{ with }\psi_{1}(x)=\left\{
	\begin{aligned}
		&x,\qquad\,\,\,\text{if }{M}^{-1}\leq x \leq {M}\\
		&M^{-1},\quad \text{if }x<M^{-1}\\
		&M,\qquad \,\text{if }x>M\\
	\end{aligned}
	\right.,
\end{equation}
and
\begin{equation}\label{7.16}
	\underline{\mathfrak{h}}^{\nu}=\psi_2(\mathfrak{h}^{\nu}),\qquad \text{ with }\psi_2(x)=\left\{
	\begin{aligned}
		&x,\qquad\,\,\text{if }|x|\leq L\\
		&-L,\quad \text{if }x<-L\\
		&L,\qquad \,\,\text{if }x>L\\
	\end{aligned}
	\right.\quad.
\end{equation}

For any fixed $t\in (0,T)$ and $\kappa>0$, let $\phi_{\kappa}(s)=\frac{1}{\kappa}\phi(\frac{s}{\kappa})$ be a standard mollifier, where $\phi:\R\to \R$ is a nonnegative smooth function with compact support such that $\int_{\R}\phi(s)ds=1$ and $\operatorname{supp}\phi=[-1,1]$.
For any $\vartheta>0$, it follows from \eqref{7.13} that there exists $\nu_{*}$ such that for all $\nu<\nu_{*}$
\begin{equation}\label{7.17}
	\sup_{t\in [0,T]}\int_{\R}\eta\big((\mathfrak{v}^{\nu},\mathfrak{h}^{\nu})(t,x)|(\tilde{\mathfrak{v}}_{X}^{\nu},\tilde{\mathfrak{h}}_{X}^{\nu})(t,x)\big)dxdt<C\mathcal{E}_{0}+\vartheta.
\end{equation}
Consider the following functional
$$
L_{\nu}(t)= \int_{0}^{T}\phi_{\kappa}(t-s)\int_{\R}\eta\big((\underline{v}^{\nu},\underline{h}^{\nu})(s,x)|(\tilde{\mathfrak{v}}_{X}^{\nu},\tilde{\mathfrak{h}}_{X}^{\nu})(s,x)\big)dxds,\quad \text{for }0<\kappa<\min\{t, T-t\}.
$$
Using \eqref{A.1-1}, \eqref{7.15}-\eqref{7.16} and the fact $\int_{\R}\phi(s)ds=1$, we have
\begin{equation}\label{7.17-1}
L_{\nu}(t)\leq \int_{0}^{T}\phi_{\kappa}(t-s)\int_{\R}\eta\big((\mathfrak{v}^{\nu},\mathfrak{h}^{\nu})(s,x)|(\tilde{\mathfrak{v}}_{X}^{\nu},\tilde{\mathfrak{h}}_{X}^{\nu})(s,x)\big)dxds\leq C\mathcal{E}_{0}+\vartheta.
\end{equation}

Similar as \cite[Lemma 5.2]{Kang-Vasseur-2021-Invent}, we have following lemma.
\begin{lemma}\label{lem7.2}
	Denote
	$$
	R_{\nu}(t):=\int_{0}^{T}\phi_{\kappa}(t-s)\int_{\R}\eta\big((\underline{v}^{\nu},\underline{h}^{\nu})(s,x)|(\bar{\mathfrak{v}}_{X_{\infty}},\bar{\mathfrak{u}}_{X_{\infty}})(s,x)\big)dxds,\quad\text{for }0<\kappa<\min\{t, T-t\},
	$$
	with $\bar{\mathfrak{v}}_{X_{\infty}}(t,x)$ and $\bar{\mathfrak{u}}_{X_{\infty}}(t,x)$ defined in \eqref{1.23}.
	Then it holds that
	\begin{equation}\label{7.18}
	\Big\vert L_{\nu}(t)-R_{\nu}(t)\Big\vert \rightarrow 0\quad \text{up to a subsequence as }\nu \to 0.
	\end{equation}
\end{lemma}

\noindent\textbf{Proof.} Since $\underline{\mathfrak{h}}^{\nu}, \tilde{\mathfrak{h}}_{X}^{\nu}$ and $\bar{\mathfrak{u}}_{X_{\infty}}$ are bounded, we have
\begin{equation*}
	\begin{aligned}
	\Big\vert |\underline{\mathfrak{h}}^{\nu}-\tilde{\mathfrak{h}}_{X}^{\nu}|^2-|\underline{\mathfrak{h}}^{\nu}-\bar{\mathfrak{u}}_{X_{\infty}}|^2\Big\vert&\leq C\Big\vert\tilde{\mathfrak{h}}_{X}^{\nu}-\bar{\mathfrak{u}}_{X_{\infty}}\Big\vert.
	\end{aligned}
\end{equation*}
For later use, we denote
	$$
\begin{aligned}
	\bar{\mathfrak{v}}_{X}(t,x):=\mathfrak{v}^{s}(x-\sigma_{1}t-X_{\nu}(t))+\mathfrak{v}^{r}(t,x)-v_{m},\\
	\bar{\mathfrak{u}}_{X}(t,x):=\mathfrak{u}^{s}(x-\sigma_{1}t-X_{\nu}(t))+\mathfrak{u}^{r}(t,x)-u_{m}.
\end{aligned}
$$
Noting
\begin{equation*}
\begin{aligned}
\Big\vert\tilde{\mathfrak{h}}_{X}^{\nu}-\bar{\mathfrak{u}}_{X_{\infty}}\Big\vert&\leq \Big\vert\tilde{\mathfrak{h}}_{X}^{\nu}-\bar{\mathfrak{u}}_{X}\Big\vert+\Big\vert\bar{\mathfrak{u}}_{X}-\bar{\mathfrak{u}}_{X_{\infty}}\Big\vert\\
&\leq \Big\vert (\tilde{\mathfrak{u}}_{X}^{\nu})^{s}-\mathfrak{u}^{s}(x-\sigma_{1}t-X_{\nu}(t))\Big\vert+\Big\vert\tilde{u}^{r}-\mathfrak{u}^{r}\Big\vert\\
&\quad +\Big\vert\mathfrak{u}^{s}(x-\sigma_{1}t-X_{\nu}(t))-\mathfrak{u}^{s}(x-\sigma_{1}t-X_{\infty}(t))\Big\vert+\Big\vert \frac{\gamma \nu}{\alpha}\Big(\big(p(\tilde{\mathfrak{v}}_{X}^{s})\big)^{\frac{\alpha}{\gamma}}\Big)_{x}\Big\vert\\
&:=I_1+I_2+I_3+I_4.
\end{aligned}
\end{equation*}
For $I_1$, it follows from \eqref{7.1-2} and Lemma \ref{lem2.3} that
\begin{equation}\label{7.19}
\begin{aligned}
\|I_1\|_{L^1(\R)}&=\int_{\R}\Big\vert\tilde{u}^{s}\left(\frac{x-\sigma_1 t-X_{\nu}(t)}{\nu}\right)-\mathfrak{u}^{s}(x-\sigma_1t-X_{\nu}(t))\Big\vert\,dx\\
&=\nu\int_{\R}\Big\vert \tilde{u}^{s}\left(y-\frac{\sigma_1 t+X_{\nu}(t)}{\nu}\right)-\mathfrak{u}^{s}(\nu y-\sigma_1t-X_{\nu}(t))\Big\vert\,dy\\
&=\nu\int_{-\infty}^{\frac{\sigma_1 t+X_{\nu}(t)}{\nu}}\Big\vert \tilde{u}^{s}\left(y-\frac{\sigma_1 t+X_{\nu}(t)}{\nu}\right)-u_{-}\Big\vert\,dy\\
&\quad +\nu\int_{\frac{\sigma_1 t+X_{\nu}(t)}{\nu}}^{\infty}\Big\vert \tilde{u}^{s}\left(y-\frac{\sigma_1 t+X_{\nu}(t)}{\nu}\right)-u_{m}\Big\vert\,dy\\
&\leq C\nu.
\end{aligned}
\end{equation}
For $I_2$ and $I_4$, it follows from \eqref{7.1-2}, Lemmas \ref{lem2.3} and \ref{lem2.2}  that
\begin{equation}\label{7.20}
\|I_{2}\|_{L^1(\R)}\leq C\v_2a(\nu),\qquad|I_4\|_{L^1(\R)}\leq C\nu\|\partial_{x}\tilde{v}^{s}(x)\|_{L^1(\R)}\leq C\v_{1}\nu.
\end{equation}
For $I_{3}$, a direct calculation shows that $\|I_3\|_{L^1(\R)}=|u_{m}-u_{-}|\cdot |X_{\nu}(t)-X_{\infty}(t)|$, then using \eqref{7.6}, we have
\begin{equation}\label{7.21}
	\int_{0}^{T}\phi_{\kappa}(t-s)\|I_{3}\|_{L^1(\R)}ds\lesssim \int_{0}^{T}\phi_{\kappa}(t-s)|X_{\nu}(s)-X_{\infty}(s)|\,ds\rightarrow 0\quad \text{as }\nu \to 0.
\end{equation}
Combining \eqref{7.19}, \eqref{7.20} and \eqref{7.21} together, we obtain
\begin{equation}\label{7.22}
	\lim\limits_{\nu\to 0}\int_{0}^{T}\phi_{\kappa}(t-s)\int_{\R}	\Big\vert |\underline{h}^{\nu}-\tilde{\mathfrak{h}}_{X}^{\nu}|^2-|\underline{h}^{\nu}-\bar{\mathfrak{u}}_{X_{\infty}}|^2\Big\vert\,dxds=0.
\end{equation}
Similarly, since $\underline{\mathfrak{v}}^{\nu}$ and $\bar{\mathfrak{v}}$ are bounded, using the definition $Q(\cdot|\cdot)$ and Taylor's expansion, we have
\begin{equation}\label{7.23}
\begin{aligned}
	&\Big\vert Q\big(\underline{\mathfrak{v}}^{\nu}|\tilde{\mathfrak{v}}_{X}^{\nu}\big)-Q\big(\underline{\mathfrak{v}}^{\nu}|\bar{\mathfrak{v}}_{X_{\infty}}\big)\Big\vert\\
	&\leq \Big\vert Q(\tilde{\mathfrak{v}}_{X}^{\nu})-Q(\bar{\mathfrak{v}}_{X_{\infty}})\Big\vert+\Big\vert Q'(\tilde{\mathfrak{v}}_{X}^{\nu})\Big\vert\cdot \Big\vert \tilde{\mathfrak{v}}_{X}^{\nu}- \bar{\mathfrak{v}}_{X_{\infty}}\Big\vert\\ &\quad+\Big(|\underline{\mathfrak{v}}^{\nu}|+|\bar{\mathfrak{v}}_{X_{\infty}}|\Big)\cdot \Big\vert Q'(\tilde{\mathfrak{v}}_{X}^{\nu})-Q'(\bar{\mathfrak{v}}_{X_{\infty}})\Big\vert\\
	&\leq C\Big\vert \tilde{\mathfrak{v}}_{X}^{\nu}-\bar{\mathfrak{v}}_{X_{\infty}}\Big\vert,
\end{aligned}
\end{equation}
which, together with \eqref{7.1-2}, Lemmas \ref{lem2.3}-\ref{lem2.2} and similar calculations as in \eqref{7.19}-\eqref{7.21}, yields that
\begin{equation}\label{7.24}
	\lim\limits_{\nu\to 0}\int_{0}^{T}\phi_{\kappa}(t-s)\int_{\R}\Big\vert Q\big(\underline{\mathfrak{v}}^{\nu}|\tilde{\mathfrak{v}}_{X}^{\nu}\big)-Q\big(\underline{\mathfrak{v}}^{\nu}|\bar{\mathfrak{v}}_{X_{\infty}}\big)\Big\vert\,dxds=0.
\end{equation}
Combining \eqref{7.24} and \eqref{7.22}, we get \eqref{7.18}. Therefore the proof of Lemma \ref{lem7.2} is complete. $\hfill\square$

Now we consider
\begin{equation}\label{7.25}
	\begin{aligned}
&\int_{0}^{T}\phi_{\kappa}(t-s)\int_{\R}\eta\big((\mathfrak{v}^{\nu},\mathfrak{h}^{\nu})(s,x)|(\bar{\mathfrak{v}}_{X_{\infty}},\bar{\mathfrak{u}}_{X_{\infty}})(s,x)\big)\,dxds\\
&=\frac{1}{2}\iint_{\{(s,x)|s\in [0,T],\,\,\mathfrak{h}^{\nu}(s)\in [-L,L]\}}\phi_{\kappa}(t-s)\Big\vert \underline{\mathfrak{h}}^{\nu}-\bar{\mathfrak{u}}_{X_{\infty}}\Big\vert^2\,dxds\\
&\quad +\frac{1}{2}\iint_{\{(s,x)|s\in [0,T],\,\,\mathfrak{h}^{\nu}(s)\notin [-L,L]\}}\phi_{\kappa}(t-s)\Big\vert \mathfrak{h}^{\nu}-\bar{\mathfrak{u}}_{X_{\infty}}\Big\vert^2\,dxds\\
&\quad +\iint_{\{(s,x)|s\in [0,T],\,\,\mathfrak{v}^{\nu}\in [M^{-1},M]\}}\phi_{\kappa}(t-s)Q(\underline{\mathfrak{v}}^{\nu}|\bar{\mathfrak{v}}_{X_{\infty}})\,dxds\\
&\quad +\iint_{\{(s,x)|s\in [0,T],\,\,\mathfrak{v}^{\nu}\notin [M^{-1},M]\}}\phi_{\kappa}(t-s)Q(\mathfrak{v}^{\nu}|\bar{\mathfrak{v}}_{X_{\infty}})\,dxds\\
&:=J_{1}+J_2+J_3+J_4.
\end{aligned}
\end{equation}
Using \eqref{7.17-1}, we have
\begin{equation}\label{7.26}
		J_{1}+J_{3}\leq R_{\nu}=(R_{\nu}-L_{\nu})+L_{\nu}\leq (R_{\nu}-L_{\nu})+C\mathcal{E}_{0}+\vartheta.
\end{equation}
For $J_{2}$, since $\bar{\mathfrak{u}}_{X_{\infty}}, \tilde{\mathfrak{h}}_{X}^{\nu}\in [-L/2,L/2]$, a direct calculation shows that
$$
\begin{aligned}
\Big\vert \mathfrak{h}^{\nu}-\bar{\mathfrak{u}}_{X_{\infty}}\Big\vert&\leq \Big\vert \mathfrak{h}^{\nu}-\tilde{\mathfrak{h}}_{X}^{\nu}\Big\vert+\Big\vert \tilde{\mathfrak{h}}_{X}^{\nu}-\bar{\mathfrak{u}}_{X_{\infty}}\Big\vert\\
&\leq 3\Big\vert \mathfrak{h}^{\nu}-\tilde{\mathfrak{h}}_{X}^{\nu}\Big\vert\qquad \text{for }\mathfrak{h}^{\nu}\notin [-L,L].
\end{aligned}
$$
Thus using \eqref{7.17}, one has
\begin{equation}\label{7.27}
J_{2}\leq \frac{9}{2}\iint_{\{(s,x)|s\in [0,T],\,\,\mathfrak{h}^{\nu}(s)\notin [-L,L]\}}\phi_{\kappa}(t-s)\Big\vert\mathfrak{h}^{\nu}-\tilde{\mathfrak{h}}_{X}^{\nu}\Big\vert^2\,dxds\leq C(\mathcal{E}_{0}+\vartheta).
\end{equation}
Similarly, since $\bar{\mathfrak{v}}_{X_{\infty}}, \tilde{\mathfrak{v}}_{X}^{\nu}\in [2M^{-1},1/2M]$, it follows from \eqref{7.17} that
\begin{equation}\label{7.28}
\begin{aligned}
J_4&\leq C\iint_{\{(s,x)|s\in [0,T],\,\,\mathfrak{v}^{\nu}\notin [M^{-1},M]\}}\phi_{\kappa}(t-s)\max\{(\mathfrak{v}^{\nu})^{-\gamma+1},\mathfrak{v}^{\nu}\}\,dxds\\
&\leq C\iint_{\{(s,x)|s\in [0,T],\,\,\mathfrak{v}^{\nu}\notin [M^{-1},M]\}}\phi_{\kappa}(t-s)Q(\mathfrak{v}^{\nu}|\tilde{\mathfrak{v}}_{X}^{\nu})\,dxds\leq C(\mathcal{E}_{0}+\vartheta).
\end{aligned}
\end{equation}
Hence, substituting \eqref{7.26}-\eqref{7.28} into \eqref{7.25}, we obtain
\begin{equation}\label{7.29}
\begin{aligned}
	&\int_{0}^{T}\phi_{\kappa}(t-s)\int_{\R}\frac{1}{2}\Big\vert \mathfrak{h}^{\nu}(s,x)-\bar{\mathfrak{u}}_{X_{\infty}}(s,x)\Big\vert^2\,dxds\\
	&+\int_{0}^{T}\phi_{\kappa}(t-s)\int_{\R}Q\big(\mathfrak{v}^{\nu}(s,x)|\bar{\mathfrak{v}}_{X_{\infty}}(s,x)\big)\,dxds\leq |R_{\nu}-L_{\nu}|+C(\mathcal{E}_{0}+\vartheta).
\end{aligned}
\end{equation}
Now, we shall prove that the LHS of \eqref{7.29} is lower semi-continuous with respect to the weak convergence in \eqref{7.14}.

Since the $L^2$ norm is weakly semi-continuous, using $\eqref{7.14}_{3}$, we immediately obtain
\begin{equation}\label{7.30}
\begin{aligned}
&\int_{0}^{T}\phi_{\kappa}(t-s)\int_{\R}\frac{1}{2}\Big\vert \mathfrak{u}_{\infty}(s,x)-\bar{\mathfrak{u}}_{X_{\infty}}(s,x)\Big\vert^2\,dxds\\
&\leq \liminf_{\nu\to 0}\int_{0}^{T}\phi_{\kappa}(t-s)\int_{\R}\frac{1}{2}\Big\vert \mathfrak{h}^{\nu}(s,x)-\bar{\mathfrak{u}}_{X_{\infty}}(s,x)\Big\vert^2\,dxds\\
&\leq C(\mathcal{E}_{0}+\vartheta)
\end{aligned}
\end{equation}
As pointed out in \cite{Kang-Vasseur-2021-Invent}, since $\mathfrak{v}_{\infty}\in L^{\infty}(0,T; L^{\infty}(\R)+\mathcal{M}(\R))$ is a measure in space, we may use the generalized functional $dQ(\mathfrak{v}^{\nu}|\bar{\mathfrak{v}}_{X_{\infty}})$ to handle the measure $\mathfrak{v}_{\infty}$. To this end, for any positive measure $dv=dv_{a}+dv_{s}=m_{a}dtdx+dv_{s}$ with $m_{a}\in L_{\mathrm{loc}}^1(\R_{+}^2)$ being the Randon-Nikodym's derivative of $v_{a}$ (by Radon-Nikodym's theorem), we define
$$
dQ(v|\bar{\mathfrak{v}}_{X_{\infty}}):=Q(m_{a}|\bar{\mathfrak{v}}_{X_{\infty}})\,dtdx+|Q'(\overline{V}(t,x))|\,dv_{s},
$$
with $\overline{V}(t,x)$ defined in \eqref{7.31}.

For later use, we denote
$$
\Omega:=\{(t,x)\in (0,T)\times \R|x<\sigma_{1}t+X_{\infty}(t)\}.
$$
Since $X_{\infty}\in BV(0,T)$, the Lebesgue measure of $\partial \Omega$ (the boundary of $\Omega$) is zero in $\R^2$. The complement of $\overline{\Omega}$ (the closure of $\Omega$) in $(0,T)\times \R$ is defined as follows:
$$
\left(\overline{\Omega}\right)^{c}=\{(t,x)\times (0,T)\times \R|x>\sigma_{1}t+X_{\infty}(t)\}.
$$
It follows from \eqref{1.23} and \eqref{7.12-1} that
\begin{equation}
\bar{\mathfrak{v}}_{X_{\infty}}(t,x)=\left\{
\begin{aligned}
	&v_{-},\quad \text{for }(t,x)\in \Omega,\\
	&\mathfrak{v}^{r},\quad\,\,\text{for }(t,x)\in \left(\overline{\Omega}\right)^{c}.
	\end{aligned}
\right.\nonumber
\end{equation}
Similar as \cite[Lemma 5.3]{Kang-Vasseur-2021-Invent}, we have following lemma.
\begin{lemma}\label{lem7.3}
Let $\Phi: (0,T)\times \R\to \R_{+}$ be any compactly supported nonnegative function and $\{v^{k}\}_{k=1}^{\infty}$ be a sequence of positive measures in $L^{\infty}((0,T)\times \R)+\mathcal{M}((0,T)\times \R)$ such that for some constant $C_{0}>0$ (independent of $k$),
$$
\iint_{(0,T)\times \R}\Phi(t,x)\,dQ\left(v^{k}|\bar{\mathfrak{v}}_{X_{\infty}}\right)(t,x)\leq C_{0},
$$
where
$$
dQ\left(v^{k}|\bar{\mathfrak{v}}_{X_{\infty}}\right)=Q\left(h_{a}^{k}|\bar{\mathfrak{v}}_{X_{\infty}}\right)dtdx+|Q'(\overline{V}(t,x))|dv_{s}^{k}(t,x),
$$
with $v^{k}=v_{a}^{k}+v_{s}^{k}$, $dv_{a}^{k}=m_{a}^{k}dtdx$ and $m_{a}^{k}\in L_{\mathrm{loc}}^1(\R_{+}^2)$ being the Randon-Nikodym's derivative of $v_{a}^k$ $\mathrm{(}$by Radon-Nikodym's theorem$\mathrm{)}$, and $\bar{V}(t,x)$ defined in \eqref{7.31}. Then there exists  a limit $v_{\infty}\in L^{\infty}((0,T)\times \R)+\mathcal{M}((0,T)\times \R)$ such that $v^{k}\rightharpoonup v_{\infty}$ in $\mathcal{M}_{\mathrm{loc}}((0,T)\times \R)$, and
$$
\begin{aligned}
&\iint_{(0,T)\times \R}\Phi(t,x)\,dQ\left(v_{\infty}|\bar{\mathfrak{v}}_{X_{\infty}}\right)(t,x)\\
&\leq \liminf_{k\to \infty}\iint_{(0,T)\times \R}\Phi(t,x)\,dQ\left(v^{k}|\bar{\mathfrak{v}}_{X_{\infty}}\right)(t,x)\leq C_{0}.
\end{aligned}
$$
\end{lemma}

\noindent\textbf{Proof.} We divide the proof into three steps.

\textbf{Step 1.} Since $v^{k}$ are positive measures in $L^{\infty}((0,T)\times \R)+\mathcal{M}((0,T)\times \R)$, then by using Radon-Nikodym's theorem, there exist positive functions $m_{a}^{k}\in L_{\mathrm{loc}}^1(\R_{+}^2)$ and $dv_{s}^{k}$ (singular part of $v_{k}$) such that
$$
dv^{k}(t,x)=m_{a}^{k}(t,x)dtdx+dv_{s}^{k}(t,x).
$$
From the definition of $Q(\cdot |\cdot)$, for any $\epsilon>0$, there exists a constant $N>0$ with $N>3\max\{v_{-},v_{+}^{-1}\}$ such that for all $v\geq N$,
\begin{equation}\label{7.33}
(|Q'(\bar{\mathbf{v}}_{\infty})|+\epsilon)v\geq Q(v|\bar{\mathbf{v}}_{\infty})\geq (|Q'(\bar{\mathbf{v}}_{\infty})|-\epsilon)v,
\end{equation}
where we have used the simplified notation $\bar{\mathbf{v}}_{\infty}:=\bar{\mathfrak{v}}_{X_{\infty}}$. We define the $N-\mathrm{truncation}$ of $m_{a}^{k}$ as follows:
$$
m_{N}^{k}:=\inf(m_{a}^{k},N).
$$
Also, we denote
$$
Q_{N}(v):=\left\{
\begin{aligned}
&Q(v),\quad \text{if }v\geq N^{-1},\\
&Q'(N^{-1})(v-N^{-1})+Q(N^{-1}),\quad \text{if }v\leq N^{-1}.
\end{aligned}
\right.
$$
Since $N^{-1}<\frac{v_{+}}{3}<\bar{\mathbf{v}}_{\infty}$, then $Q_{N}'(\bar{\mathbf{v}}_{\infty})=Q'(\bar{\mathbf{v}}_{\infty})$. For any $v_1, v_2\geq 0$, we consider the following relative functional
$$
Q_{N}(v_{1}|v_{2}):=Q_{N}(v_1)-Q_{N}(v_2)-Q_{N}'(v_2)(v_1-v_2).
$$
We claim that
\begin{equation}\label{7.34}
dQ_{N}(v^{k}|\bar{\mathbf{v}}_{\infty})\geq Q_{N}(m_{N}^{k}|\bar{\mathbf{v}}_{\infty})dtdx+(|Q_{N}'(\overline{V})|-\epsilon)(dv^{k}-m_{N}^{k}dtdx)-2\epsilon dv^{k},
\end{equation}
hence it means $dQ_{N}(v^{k}|\bar{\mathbf{v}}_{\infty})-\big[Q_{N}(m_{N}^{k}|\bar{\mathbf{v}}_{\infty})dtdx+(|Q_{N}'(\overline{V})|-\epsilon)(dv^{k}-m_{N}^{k}dtdx)-2\epsilon dv^{k}\big]$ is a nonnegative measure. Indeed, this can be verified as follows:

{\it{Case 1.}} If $h_{a}^{k}\leq N$, then it holds that $h_{N}^{k}=h_{a}^{k}$, and
\begin{equation}\label{7.35}
	\begin{aligned}
	dQ_{N}(v^{k}|\bar{\mathbf{v}}_{\infty})&=Q_{N}(m_{a}^{k}|\bar{\mathbf{v}}_{\infty})dtdx+|Q_{N}'(\overline{V})|dv_{s}^{k}\\
	&=Q(m_{N}^{k}|\bar{\mathbf{v}}_{\infty})dtdx+|Q_{N}'(\overline{V})|dv_{s}^{k}
	\\&\geq Q_{N}(m_{N}^{k}|\bar{\mathbf{v}}_{\infty})dtdx+(|Q_{N}'(\overline{V})|-\epsilon)(dv^{k}-m_{N}^{k}dtdx),
	\end{aligned}
\end{equation}
where we have used the fact that the measure $dv^{k}-m_{N}^{k}dtdx=(m_{a}^{k}-m_{N}^{k})dtdx+dv_{s}^{k}=dv_{s}^{k}$ is positive in the last inequality.

{\it{Case 2.}} If $m_{a}^{k}>N$, then $m_{N}^{k}=N$. Using \eqref{7.33} to obtain
\begin{equation}\label{7.36}
dQ_{N}(v^{k}|\bar{\mathbf{v}}_{\infty})=Q_{N}(m_{a}^{k}|\bar{\mathbf{v}}_{\infty})dtdx+|Q_{N}'(\overline{V})|dv_{s}^{k}\geq (|Q_{N}'(\bar{\mathbf{v}}_{\infty})|-\epsilon)m_{a}^{k}dtdx+|Q_{N}'(\overline{V})|dv_{s}^{k}.
\end{equation}
Since $\bar{\mathbf{v}}_{\infty}=\overline{V}$ a.e. for Lebesgue measure $dtdx$ on $(0,T)\times \R$, thus, using \eqref{7.33} and \eqref{7.36}, we have
\begin{equation}\label{7.37}
\begin{aligned}
dQ_{N}(v^{k}|\bar{\mathbf{v}}_{\infty})&\geq (|Q_{N}'(\overline{V})|-\epsilon)dv^{k}\\
&=(|Q_{N}'(\overline{V})|+\epsilon)N\,dtdx+(|Q_{N}'(\overline{V})|-\epsilon)(dv^{k}-N\,dtdx)-2\epsilon N\,dtdx\\
&\geq (|Q_{N}'(\bar{\mathbf{v}}_{\infty})|+\epsilon)N\,dtdx+(|Q_{N}'(\overline{V})|-\epsilon)(dv^{k}-m_{N}^{k}dtdx)-2\epsilon dv^{k}\\
&\geq Q_{N}(N|\bar{\mathbf{v}}_{\infty})dtdx+(|Q_{N}'(\overline{V})|-\epsilon)(dv^{k}-m_{N}^{k}dtdx)-2\epsilon dv^{k}\\
&\geq Q_{N}(m_{N}^{k}|\bar{\mathbf{v}}_{\infty})+(|Q_{N}'(\overline{V})|-\epsilon)(dv^{k}-m_{N}^{k}dtdx)-2\epsilon dv^{k}.
\end{aligned}
\end{equation}

Hence, it follows from \eqref{7.34} that
\begin{equation}\label{7.38}
\begin{aligned}
	C_{0}&\geq \limsup_{k\to \infty}\iint_{(0,T)\times \R}\Phi(t,x)dQ_{N}\left(v^{k}|\bar{\mathbf{v}}_{\infty}\right)(t,x)\\
	&\geq \limsup_{k\to \infty}\iint_{(0,T)\times \R}\Phi(t,x)\Big[Q_{N}(m_{N}^{k}|\bar{\mathbf{v}}_{\infty})dtdx+(|Q_{N}'(\overline{V})|-\epsilon)(dv^{k}-m_{N}^{k}dtdx)-2\epsilon dv^{k}\Big].
\end{aligned}
\end{equation}

\textbf{Step 2.} We denote $\Omega_{1}:=\left(\overline{\Omega}\right)^{c}$, and define
$$
\Omega_{1}^{\delta}:=\{(t,x)\in \Omega_{1}|d((t,x)|\Omega_{1}^{c})>\delta\},\quad \text{for any } \delta>0,
$$
where $d((t,x)|\Omega_{1}^{c})$ is the distance function between the point $(t,x)$ and the closed set $\Omega_{1}^{c}$. Then we take a positive smooth function $\psi_{1}^{\delta}$ such that
\begin{equation}\label{7.39}
\psi_{1}^{\delta}=\left\{
\begin{aligned}
	&1,\quad \text{on }\big(\Omega_{1}^{\delta}\big)^{c},\\
	&0,\quad \text{on } \Omega_{1}^{2\delta}.
\end{aligned}
\right.
\end{equation}
Since $v_{-}>v_{m}>v_{+}$, we have $|Q_{N}'(v_{-})|\leq |Q_{N}'(\mathfrak{v}^{r})(t,x)|$ for all $(t,x)\in (0,T)\times \R$, then using \eqref{7.39}, one obtains
\begin{align}
	&\iint_{(0,T)\times \R}\Phi(t,x)(|Q_{N}'(\overline{V})|-\epsilon)\,(dv^k-m_{N}^{k}dtdx)\nonumber\\
	&=\iint_{\Omega_{1}^{c}}\Phi(t,x)(|Q_{N}'(v_{-})|-\epsilon)\,(dv^k-m_{N}^{k}dtdx)+\iint_{\Omega_{1}\backslash \Omega_{1}^{\delta}}\Phi(t,x)(|Q_{N}'(\mathfrak{v}^{r})|-\epsilon)\,(dv^k-m_{N}^{k}dtdx)\nonumber\\
		&\quad+\iint_{\Omega_{1}^{\delta}\backslash \Omega_{1}^{2\delta}}\Phi(t,x)(|Q_{N}'(\mathfrak{v}^{r})|-\epsilon)\,(dv^k-m_{N}^{k}dtdx) +\iint_{\Omega_{1}^{2\delta}}\Phi(t,x)(|Q_{N}'(\mathfrak{v}^{r})|-\epsilon)\,(dv^k-m_{N}^{k}dtdx)\nonumber\\
		&\geq \iint_{\Omega_{1}^{c}}\Phi(t,x)(|Q_{N}'(v_{-})|-\epsilon)\,(dv^k-m_{N}^{k}dtdx)+\iint_{\Omega_{1}\backslash \Omega_{1}^{\delta}}\Phi(t,x)\psi_{1}^{\delta}(|Q_{N}'(v_{-})|-\epsilon)\,(dv^k-m_{N}^{k}dtdx)\nonumber\\
		&\quad +\iint_{\Omega_{1}^{\delta}\backslash \Omega_{1}^{2\delta}}\Phi(t,x)\psi_{1}^{\delta}(|Q_{N}'(v_{-})|-\epsilon)\,(dv^k-m_{N}^{k}dtdx)\nonumber\\
		&\quad +\iint_{\Omega_{1}^{\delta}\backslash \Omega_{1}^{2\delta}}\Phi(t,x)(1-\psi_{1}^{\delta})(|Q_{N}'(\mathfrak{v}^{r})|-\epsilon)\,(dv^k-m_{N}^{k}dtdx)\nonumber\\
		&\quad +\iint_{\Omega_{1}^{2\delta}}\Phi(t,x)(1-\psi_{1}^{\delta})(|Q'(\mathfrak{v}^{r})|-\epsilon)\,(dv^k-m_{N}^{k}dtdx)\nonumber\\
		&=(|Q_{N}'(v_{-})|-\epsilon)\iint_{(0,T)\times \R} \Phi(t,x)\psi_{1}^{\delta}\,(dv^k-m_{N}^{k}dtdx)\nonumber\\
		&\quad  +\iint_{(0,T)\times \R}\Phi(t,x)(1-\psi_1^{\delta})(|Q_{N}'(\mathfrak{v}^{r})|-\epsilon)\,(dv^k-m_{N}^{k}dtdx).\label{7.40}
\end{align}
Hence, substituting \eqref{7.40} into \eqref{7.38}, we obtain
\begin{equation}\label{7.41}
\begin{aligned}
	C_{0}&\geq \limsup_{k\to \infty}\Big[\iint_{\Omega}\Phi(t,x)Q_{N}(m_{N}^{k}|v_{-})\,dtdx+\iint_{\Omega_{1}}\Phi(t,x)Q_{N}(m_{N}^{k}|\mathfrak{v}^{r})\,dtdx\\
	&\quad +(|Q_{N}'(v_{-})|-\epsilon)\iint_{(0,T)\times \R} \Phi(t,x)\psi_{1}^{\delta}\,(dv^{k}-m_{N}^{k}dtdx)\\
	&\quad  +\iint_{(0,T)\times \R}\Phi(t,x)(1-\psi_1^{\delta})(|Q_{N}'(\mathfrak{v}^{r})|-\epsilon)\,(dv^{k}-m_{N}^{k}dtdx)-2\epsilon\iint_{(0,T)\times \R}\Phi(t,x)dv^{k}\Big].
\end{aligned}
\end{equation}
Since $|m_{N}^{k}|\leq N$ for all $k$, there exists $v_{\ast}$ such that
$$
m_{N}^{k}\overset{\ast}{ \rightharpoonup}m_{\ast}\quad \text{in }L^{\infty}((0,T)\times \R)\quad \text{as }k\to \infty.
$$
Recall
\begin{equation}\label{7.42}
Q_{N}(m_{N}^{k}|\mathfrak{v}^{r})=Q_{N}(m_{N}^{k})-Q_{N}'(\mathfrak{v}^{r})m_{N}^{k}-\big[Q_{N}(\mathfrak{v}^{r})-\mathfrak{v}^{r}Q_{N}'(\mathfrak{v}^{r})\big].
\end{equation}
By the definition of $Q_{N}(v)$, the function $v\mapsto Q_{N}(v)$ is convex, and the weak lower semi-continuity of convex functions implies that
\begin{equation}\label{7.43}
\liminf_{k\to \infty}\iint_{\Omega_{1}}\Phi(t,x)Q_{N}(m_{N}^{k})\,dtdx\geq \liminf_{k\to \infty}\iint_{\Omega_{1}}\Phi(t,x)Q_{N}(m_{*})\,dtdx.
\end{equation}
Moreover, since $Q_{N}'(\mathfrak{v}^{r})\Phi(t,x)\in L^1((0,T)\times \R)$, then it holds
\begin{equation}\label{7.44}
\lim\limits_{k\to \infty}\iint_{\Omega_{1}}\Phi(t,x)Q_{N}'(\mathfrak{v}^{r})m_{N}^{k}\,dtdx=\lim\limits_{k\to \infty}\iint_{\Omega_{1}}\Phi(t,x)Q_{N}'(\mathfrak{v}^{r})m_{*}\,dtdx.
\end{equation}
Hence, it follows from \eqref{7.42}-\eqref{7.44} that
\begin{equation}\label{7.45}
\liminf_{k\to \infty}\iint_{\Omega_{1}}\Phi(t,x)Q_{N}(m_{N}^{k}|\mathfrak{v}^{r})\,dtdx\geq \liminf_{k\to \infty}\iint_{\Omega_{1}}\Phi(t,x)Q_{N}(m_{*}|\mathfrak{v}^{r})\,dtdx.
\end{equation}
Similarly, for $Q(m_{N}^{k}|v_{-})$, it holds
\begin{equation}\label{7.46}
	\liminf_{k\to \infty}\iint_{\Omega}\Phi(t,x)Q_{N}(m_{N}^{k}|v_{-})\,dtdx\geq \liminf_{k\to \infty}\iint_{\Omega}\Phi(t,x)Q_{N}(m_{*}|v_{-})\,dtdx.
\end{equation}
Also, since $v^{k}\rightharpoonup v_{\infty}$ in $\mathcal{M}_{\mathrm{loc}}((0,T)\times \R)$, thus
\begin{equation}\label{7.47}
v^{k}-v_{N}^{k}\rightharpoonup v_{\infty}-v_{\ast}\quad \text{in }\mathcal{M}_{\mathrm{loc}}((0,T)\times \R).
\end{equation}
Taking $k\to \infty$ in \eqref{7.41}, using \eqref{7.45}-\eqref{7.47} and the fact that $Q'(\mathfrak{v}^{r})$ is continuous, we get
\begin{equation}\label{7.48}
\begin{aligned}
	C_{0}&\geq \Big[\iint_{\Omega}\Phi(t,x)Q_{N}(m_{\ast}|v_{-})\,dtdx+\iint_{\Omega_{1}}\Phi(t,x)Q_{N}(m_{\ast}|\mathfrak{v}^{r})\,dtdx\\
	&\quad +(|Q_{N}'(v_{-})|-\epsilon)\iint_{(0,T)\times \R}\Phi(t,x)\psi_{1}^{\delta}\,d(v_{\infty}-v_{\ast})\\
	&\quad +\iint_{(0,T)\times \R}\Phi(t,x)(1-\psi_{1}^{\delta})(|Q_{N}'(\mathfrak{v}^{r})|-\epsilon)\,\,d(v_{\infty}-v_{\ast})\\
	&\quad -2\epsilon\iint_{(0,T)\times \R}\Phi(t,x)\,dv_{\infty}\Big]:=\mathcal{R}.
\end{aligned}
\end{equation}

\textbf{Step 3.} Applying Radon-Nikodym's theorem to the positive measure $v_{\infty}$, we have
\begin{equation}\label{7.49}
dv_{\infty}=dv_{a}+dv_{s}\quad \text{with } dv_{a}=m_{a}dtdx,
\end{equation}
where $m_{a}\in L_{\mathrm{loc}}^{1}(\R_{+}^2)$ is the Radon-Nikodym's derivative of $v_{a}$ and the singular part $dv_{s}$ is a positive measures. Since the measure $v^{k}-v_{N}^{k}$ is positive, it follows from \eqref{7.47} and \eqref{7.49} that $v_{\infty}-v_{\ast}$ is nonnegative. Moreover, by uniqueness of Radon-Nikodym's theorem, it holds
$$
dv_{\infty}-dv_{\ast}=dv_{\infty}-m_{\ast}dtdx=(m_{a}-m_{\ast})\,dtdx+dv_{s},
$$
We rewrite $\mathcal{R}$ in \eqref{7.48} as
$$
\mathcal{R}=\mathcal{R}_{1}+\mathcal{R}_{2}+\mathcal{R}_{3},
$$
with
$$
\begin{aligned}
&\mathcal{R}_{1}:=|Q_{N}'(v_{-})|\iint_{(0,T)\times \R}\Phi \psi_{1}^{\delta}\,dv_{s}+\iint_{(0,T)\times \R}\Phi\cdot(1-\psi_{1}^{\delta})|Q_{N}'(\mathfrak{v}^{r})|\,dv_{s},\\
&\begin{aligned}
\mathcal{R}_{2}:=&\iint_{\Omega}\Phi Q_{N}(m_{\ast}|v_{-})\,dtdx+|Q_{N}'(v_{-})|\iint_{(0,T)\times \R}\Phi\psi_{1}^{\delta}(m_{a}-m_{\ast})\,dtdx\\
	&+\iint_{\Omega_{1}}\Phi Q_{N}(m_{\ast}|\mathfrak{v}^{r})\,dtdx+\iint_{(0,T)\times \R}\Phi\cdot (1-\psi_{1}^{\delta})|Q_{N}'(\mathfrak{v}^{r})|(m_{a}-m_{\ast})\,dtdx,
\end{aligned}\\
&\mathcal{R}_{3}:=-3\epsilon\iint_{(0,T)\times \R}\Phi\,dv_{\infty}+\epsilon\iint_{(0,T)\times \R}\Phi m_{\ast}\,dtdx.
\end{aligned}
$$
Since $\overline{\Omega}\subset (\Omega_{1}^{\delta})^{c}$, it follows from \eqref{7.39} that
\begin{equation}\label{7.50}
\mathcal{R}_{1}\geq |Q_{N}'(v_{-})|\iint_{\overline{\Omega}}\Phi(t,x)\,dv_{s}+\iint_{\Omega_{1}^{2\delta}}\Phi(t,x) |Q_{N}'(\mathfrak{v}^{r})|\,dv_{s}.
\end{equation}
Noting $dv_{s}$ is a positive measure, an application of monotone  convergence theorem shows that
\begin{equation}\label{7.51}
	\begin{aligned}
	\liminf_{\delta\to 0}\mathcal{R}_{1}&\geq |Q_{N}'(v_{-})|\iint_{\overline{\Omega}}\Phi\,dv_s+\lim_{\delta\to 0}\iint_{\Omega_{1}^{2\delta}}\Phi |Q_{N}'(\mathfrak{v}^{r})|\,dv_{s}\\
	&\geq Q_{N}'(v_{-})|\iint_{\overline{\Omega}}\Phi\,dv_s+\iint_{\Omega_{1}}\Phi |Q_{N}'(\mathfrak{v}^{r})|\,dv_{s}\\
	&=\iint_{(0,T)\times \R}\Phi |Q_{N}'(\overline{V})|\,dv_{s}.
	\end{aligned}
\end{equation}
For $\mathcal{R}_{2}$, it follows from \eqref{7.39} that
\begin{equation}\label{7.52}
\begin{aligned}
\mathcal{R}_{2}&\geq \iint_{\Omega}\Phi \Big[Q_{N}(m_{*}|v_{-})+|Q_{N}'(v_{-})|(m_{a}-m_{\ast})\Big]\,dtdx\\
&\quad + \iint_{\Omega_{1}^{2\delta}}\Phi \Big[Q_{N}(m_{*}|\mathfrak{v}^{r})+|Q_{N}'(\mathfrak{v}^{r})|(m_{a}-m_{\ast})\Big]\,dtdx.
\end{aligned}
\end{equation}
For any positive functions $w_{1}, w_{2}, w_{3}\geq 0$, a direct calculation shows that
\begin{equation}\label{7.53}
	Q_{N}(w_{1}+w_{2}|w_{3})-Q_{N}(w_{1}|w_{3})-|Q_{N}'(w_3)|w_{2}=Q_{N}(w_{1}+w_{2})-Q_{N}(w_1)\leq 0.
\end{equation}
Substituting \eqref{7.53} into \eqref{7.52}, we obtain
\begin{equation}\label{7.54}
\mathcal{R}_{2}\geq \iint_{\Omega}\Phi Q_{N}(m_{a}|v_{-})\,dtdx+\iint_{\Omega_{1}^{2\delta}}\Phi Q_{N}(m_{a}|\mathfrak{v}^{r})\,dtdx,
\end{equation}
which, monotone  convergence theorem, yields that
\begin{equation}\label{7.55}
\liminf_{\delta\to 0}\mathcal{R}_{2}\geq \iint_{(0,T)\times \R}\Phi Q_{N}(m_{a}|\bar{\mathbf{v}})\,dtdx.
\end{equation}
Combining \eqref{7.55}, \eqref{7.51} and \eqref{7.48}, we have
\begin{equation}\label{7.56}
\iint_{(0,T)\times \R}\Phi Q_{N}(m_{a}|\bar{\mathbf{v}})\,dtdx+\iint_{(0,T)\times \R}\Phi |Q_{N}'(\overline{V})|\,dv_{s}\leq C_{0}+3\epsilon\iint_{(0,T)\times \R}\Phi\,dv_{\infty}.
\end{equation}
Taking $N\to \infty$ in \eqref{7.56} and using Fatou's Lemma, one has
\begin{equation}\label{7.57}
	\iint_{(0,T)\times \R}\Phi Q(m_{a}|\bar{\mathbf{v}})\,dtdx+\iint_{(0,T)\times \R}\Phi |Q'(\overline{V})|\,dv_{s}\leq C_{0}+3\epsilon\iint_{(0,T)\times \R}\Phi\,dv_{\infty}.
\end{equation}
Taking $\epsilon\to 0$ in \eqref{7.57}, we conclude
\begin{equation}\label{7.58}
	\begin{aligned}
		&\iint_{(0,T)\times \R}\Phi(t,x)\,dQ\left(v_{\infty}|\bar{\mathfrak{v}}_{X_{\infty}}\right)(t,x)\\
		&=\iint_{(0,T)\times \R}\Phi(t,x) Q(m_{a}|\bar{\mathbf{v}})\,dtdx+\iint_{(0,T)\times \R}\Phi |Q'(\overline{V})|\,dv_{s}\leq C_{0}.
	\end{aligned}
\end{equation}
Therefore, the proof of Lemma \ref{lem7.3} is complete. $\hfill\square$

\bigskip

With the help of Lemma \ref{lem7.3}, we now prove \eqref{1.21}. For any $R>0$, we define a positive smooth function $\psi_{0}^{R}(x)$ with satisfying
\begin{equation}
\psi_{0}^{R}(x)=\left\{
\begin{aligned}
	&1,\quad \text{for }|x|\leq R,\\
	&0,\quad \text{for }|x|\geq 2R.
\end{aligned}
\right.\nonumber
\end{equation}
It follows from \eqref{7.29} that
\begin{equation}\label{7.60}
\iint_{(0,T)\times \R}\phi_{\kappa}(t-s)\psi_{0}^{R}(x)Q\big(\mathfrak{v}^{\nu}(s,x)|\bar{\mathfrak{v}}_{X_{\infty}}\big)\,dsdx\leq |R_{\nu}-L_{\nu}|+C(\mathcal{E}_{0}+\vartheta).
\end{equation}
which, together with \eqref{7.14} and Lemmas \ref{lem7.2}-\ref{lem7.3}, yields
\begin{equation}\label{7.61}
\begin{aligned}
&\iint_{(0,T)\times \R}\phi_{\kappa}(t-s)\psi_{0}^{R}(x)\,dQ\big(\mathfrak{v}_{\infty}(s,x)|\bar{\mathfrak{v}}_{X_{\infty}}\big)\leq
C(\mathcal{E}_{0}+\vartheta).
\end{aligned}
\end{equation}
Here, by Radon-Nikodym's theorem, $\mathfrak{v}_{\infty}=\mathfrak{v}_{a}+\mathfrak{v}_{s}$ with $d\mathfrak{v}_{a}=\mathfrak{m}_{a}\,dtdx$, and
$$
dQ\big(\mathfrak{v}_{\infty}|\bar{\mathfrak{v}}_{X{\infty}}\big)(t,x)=Q\big(\mathfrak{m}_{a}|\bar{\mathfrak{v}}_{X{\infty}}\big)\,dtdx+|Q'(\overline{V}(t,x))|\,d\mathfrak{v}_{s}(t,x),
$$
with $\overline{V}(t,x)$ defined in \eqref{7.31}.

Therefore, taking $R\to \infty$ in \eqref{7.61} and using Fatou's lemma, we obtain
\begin{equation}\label{7.62}
\iint_{(0,T)\times \R}\phi_{\kappa}(t-s)\,dQ\big(\mathfrak{v}_{\infty}(s,x)|\bar{\mathfrak{v}}_{X_{\infty}}\big)\leq C(\mathcal{E}_{0}+\vartheta)
\end{equation}
Combining \eqref{7.62} and \eqref{7.30}, we conclude
\begin{equation}\label{7.63}
\begin{aligned}
&\frac{1}{2}\iint_{(0,T)\times \R}\phi_{\kappa}(t-s)\Big\vert \mathfrak{u}_{\infty}(s,x)-\bar{\mathfrak{u}}_{X_{\infty}}(s,x)\Big\vert^2\,dxds\\
&\quad +\iint_{(0,T)\times \R}\phi_{\kappa}(t-s)\,dQ\big(\mathfrak{v}_{\infty}(s,x)|\bar{\mathfrak{v}}_{X_{\infty}}\big)\leq C(\mathcal{E}_{0}+\vartheta).
\end{aligned}	
\end{equation}
Taking $\kappa \to 0$ in \eqref{7.63} and using Lebesgue point theorem,
we get
$$
dQ\big(\mathfrak{v}_{\infty}|\bar{\mathfrak{v}}_{X_{\infty}}\big)\in L^{\infty}(0,T;\mathcal{M}(\R)),
$$
and for a.e. $t\in (0,T)$,
$$
\int_{\R}\frac{1}{2}\Big\vert \mathfrak{u}_{\infty}(t,x)-\bar{\mathfrak{u}}_{X_{\infty}}(t,x)\Big\vert^2\,dx+\left(\int_{x\in \R}dQ\big(\mathfrak{v}_{\infty}(t,\cdot)|\bar{\mathfrak{v}}_{X_{\infty}}(t,\cdot)\big)\right)\leq C(\mathcal{E}_{0}+\vartheta).
$$
Since $\vartheta>0$ is arbitrary, we deduce that
$$
\int_{\R}\frac{1}{2}\Big\vert \mathfrak{u}_{\infty}(t,x)-\bar{\mathfrak{u}}_{X_{\infty}}(t,x)\Big\vert^2\,dx+\left(\int_{\R}dQ\big(\mathfrak{v}_{\infty}(t,\cdot)|\bar{\mathfrak{v}}_{X_{\infty}}(t,\cdot)\big)\right)\leq C\mathcal{E}_{0},
$$
which completes the proof of \eqref{1.21}.

\subsection{Proof of \eqref{1.24}} Motivated by \cite{Kang-Vasseur-2021-Invent}, we first show weak continuity of $\mathfrak{v}_{\infty}$, that is, for any $\varphi\in C_{0}(\R)$ and $t_{0}\in (0,T)$, we shall prove
\begin{equation}\label{7.64}
\begin{aligned}
&\lim_{t\to t_{0}}\langle \varphi, \mathfrak{v}_{\infty}(t,dx)\rangle=\langle \varphi, \mathfrak{v}_{\infty}(t_{0},dx)\rangle,\\
&\lim_{t\to 0}\langle \varphi, \mathfrak{v}_{\infty}(t,dx)\rangle=\langle \varphi, \mathfrak{v}_{0}(x)\rangle.
\end{aligned}
\end{equation}

It follows from \eqref{7.13} that
$$
\mathfrak{h}^{\nu}=\mathfrak{u}^{\nu}-\nu\frac{\gamma}{\alpha} (p(\mathfrak{v})^{\frac{\alpha}{\gamma}})_{x}\quad \text{is uniformly bounded in } L^{\infty}(0,T;L_{\mathrm{loc}}^2(\R)).
$$
Noting
\begin{equation*}
\begin{aligned}
	&\nu\frac{\gamma}{\alpha}\Big\vert \big(p(\mathfrak{v}^{\nu})^{\frac{\alpha}{\gamma}}\big)_{x}\Big\vert=\nu \mathfrak{v}^{\beta}\big\vert p(\mathfrak{v}^{\nu})_{x}\big\vert\leq \nu (\mathfrak{v}^{\nu})^{\beta}\big\vert \big(p(\mathfrak{v}^{\nu})-p(\tilde{\mathfrak{v}}_{X}^{\nu})\big)_{x}\big\vert+\nu (\mathfrak{v}^{\nu})^{\beta}\big\vert \big(p(\tilde{\mathfrak{v}}_{X}^{\nu})\big)_{x}\big\vert,
\end{aligned}
\end{equation*}
which, together with \eqref{7.13} and the fact that $\mathfrak{v}^{\beta}\leq C\big(1+Q(\mathfrak{v}^{\nu}|\tilde{\mathfrak{v}}_{X}^{\nu})\big)$, implies that
\begin{equation}\label{7.66}
\nu \frac{\gamma}{\alpha}\left(p(\mathfrak{v}^{\nu})\right)_{x}\quad \text{is uniformly bounded in }L^2(0,T;L_{\mathrm{loc}}^1(\R)).
\end{equation}
Thus
$$
\mathfrak{u}^{\nu}=\mathfrak{h}^{\nu}+\nu\frac{\gamma}{\alpha} (p(\mathfrak{v})^{\frac{\alpha}{\gamma}})_{x}\quad \text{is uniformly bounded in } L^{2}(0,T;L_{\mathrm{loc}}^1(\R)),
$$
which, together with $\mathfrak{v}_{t}^{\nu}-\mathfrak{u}_{x}^{\nu}=0$, implies that
$$
\mathfrak{v}_{t}^{\nu}\quad \text{is uniformly bounded in }L^{2}(0,T;W_{\mathrm{loc}}^{-1,1}(\R)).
$$
Using \eqref{7.13} and the fact that $\mathfrak{v}\leq C\big(1+Q(\mathfrak{v}^{\nu}|\tilde{\mathfrak{v}}_{X}^{\nu})\big)$, one obtains $\mathfrak{v}^{\nu}$ is uniformly bounded in $L^{\infty}(0,T; L_{\mathrm{loc}}^1(\R))$, then it follows from Aubin-Lions Lemma that
$$
\mathfrak{v}^{\nu}\to \mathfrak{v}_{\infty}\quad \text{(up to subsequence) in }C([0,T]; W_{\mathrm{loc}}^{-s,1}(\R)),\quad s>0,
$$
which, together with $\eqref{1.17-2}_1$, completes the proof \eqref{7.64}.
\medskip

Now, we are going to complete the proof \eqref{1.24}. Since $X_{\infty}\in BV(0,T)$, there exists a positive constant $r=r(T)$ such that $\|X_{\infty}\|_{L^{\infty}(0,T)}=r$. Let $t_{0}\in (0,T)$ be a point so that \eqref{7.10-1}, \eqref{7.12}-\eqref{7.12-1} and \eqref{1.21} hold. Let $\psi(x)$ be a nonnegative smooth function satisfying
\begin{equation}\label{7.67}
\psi(x)=\left\{
\begin{aligned}
	&1,\quad \text{if }x\in [-r+\sigma_{1}T, \frac{\sigma_{1}t_0}{2}],\\
	&0,\quad \text{if }x\leq -r-t_{0}+\sigma_{1}T\text{ or }x\geq 0,
\end{aligned}
\right.\qquad \|\psi'(x)\|_{L^{\infty}(\R)}\leq \max\{\frac{2}{t_{0}},\frac{4}{|\sigma_{1}|t_{0}}\}.
\end{equation}
Similar as in \cite{Kang-Vasseur-2021-Invent}, we take a nonnegative smooth function $\theta(x)$ satisfying
$$
\theta(s)=\theta(-s), \quad \int_{\R}\theta(x)\,dx=1,\quad \operatorname{supp}\theta(x)=[-1,1],
$$
and denote
$$
\theta_{\zeta}(s):=\frac{1}{\zeta}\theta\big(\frac{s-\zeta}{\zeta}\big)\quad \text{for any }\zeta>0.
$$
For any $t\in (t_0,T)$ at which \eqref{7.10-1}, \eqref{7.12}-\eqref{7.12-1} and \eqref{1.21} hold, we define a nonnegative smooth function
$$
\varphi_{t,\zeta}(s):=\int_{0}^{s}\big(\theta_{\zeta}(z-t_{0})-\theta_{\zeta}(z-t)\big)\,dz,\quad \text{for }\zeta\in (0,\min\{t-t_{0},T-t\}).
$$
Since $\partial_{t}\mathfrak{v}_{\infty}-\partial_{x}\mathfrak{u}_{\infty}=0$ holds in the sense of distributions, one has
\begin{equation}\label{7.68}
\iint_{(0,T)\times \R}\varphi_{t.\zeta}'(s)\psi(x)\,d\mathfrak{v}_{\infty}(s,x)-\iint_{(0,T)\times \R}\varphi_{t,\zeta}(s)\psi'(x)\mathfrak{u}_{\infty}\,dsdx=0.
\end{equation}
We decompose the left hand side of \eqref{7.68} into three parts:
\begin{equation}
\begin{aligned}
&I_{1}^{\zeta}:=\iint_{(0,T)\times \R}\theta_{\zeta}(s-t_0)\psi(x)\,d\mathfrak{v}_{\infty}(s,x),\\
&I_{2}^{\zeta}:=-\iint_{(0,T)\times \R}\theta_{\zeta}(s-t)\psi(x)\,d\mathfrak{v}_{\infty}(s,x),\\
&I_{3}^{\zeta}:=-\iint_{(0,T)\times \R}\varphi_{t,\zeta}(s)\psi'(x)\mathfrak{u}_{\infty}dxds.
\end{aligned}\nonumber
\end{equation}
It follows from \eqref{7.64} that
$$
\lim_{\zeta\to 0}I_{1}^{\zeta}=\int_{\R}\psi(x)\,\mathfrak{v}_{\infty}(t_{0},dx),\quad \lim_{\zeta\to 0}I_{2}^{\zeta}=-\int_{\R}\psi(x)\,\mathfrak{v}_{\infty}(t,dx).
$$
Moreover, a direct calculation shows that
$$
\lim_{\zeta\to 0}I_{3}^{\zeta}=-\int_{t_{0}}^{t}\int_{\R}\psi'(x)\mathfrak{u}_{\infty}(s,x)dsdx.
$$
Hence, taking $\zeta\to 0$ in \eqref{7.68} to have
\begin{equation}\label{7.70}
\int_{\R}\psi(x)\big(\mathfrak{v}_{\infty}(t,dx)-\mathfrak{v}_{\infty}(t_{0},dx)\big)+\int_{t_{0}}^{t}\int_{\R}\psi'(x)\mathfrak{u}_{\infty}(s,x)\,dxds=0.
\end{equation}
For later use, we denote
$$
\begin{aligned}
&
\begin{aligned}
J_{1}:&=\int_{\R}\psi(x)\big(\mathfrak{v}_{\infty}(t,dx)-\mathfrak{v}_{\infty}(t_{0},dx)\big)\\
&=\int_{\R}\psi(x)\big(\mathfrak{v}_{\infty}(t,dx)-\bar{\mathfrak{v}}_{X_{\infty}}(t,x)\,dx\big) +\int_{\R}\psi(x)\big(\bar{\mathfrak{v}}_{X_{\infty}}(t,x)-\bar{\mathfrak{v}}_{X_{\infty}}(t_{0},x)\big)\,dx\\
&\quad +\int_{\R}\psi(x)\big(\bar{\mathfrak{v}}_{X_{\infty}}(t_{0},x)dx-\mathfrak{v}_{X_{\infty}}(t_{0},dx)\big)\\
&=J_{11}+J_{12}+J_{13},
\end{aligned}\\
&
\begin{aligned}
J_{2}:&=\int_{t_{0}}^{t}\int_{\R}\psi'(x)\mathfrak{u}_{\infty}(s,x)dxds\\
&=\int_{t_{0}}^{t}\int_{\R}\psi'(x)\big(\mathfrak{u}_{\infty}(s,x)-\bar{\mathfrak{u}}_{X_{\infty}}(s,x)\big)dxds +\int_{t_{0}}^{t}\int_{\R}\psi'(x)\bar{\mathfrak{u}}_{X_{\infty}}(s,x)\,dxds\\
&=J_{21}+J_{22}.
\end{aligned}
\end{aligned}
$$
It follows from \eqref{7.12-1} that
\begin{equation}\label{7.71}
-r+\sigma_{1}T\leq \sigma_{1}t+X_{\infty}(t)\leq \frac{\sigma_{1}}{2}t_{0}<0\qquad \text{for a.e. }t\in [t_{0},T).
\end{equation}
Recalling \eqref{1.23} and using \eqref{7.67} and \eqref{7.71}, a direct calculation shows that
\begin{equation}\label{7.72}
\begin{aligned}
	&J_{12}=(v_{-}-v_{m})[X_{\infty}(t)-X_{\infty}(t_{0})+\sigma_{1}(t-t_0)],\\
	&J_{22}=(u_{-}-u_{m})(t-t_{0}).
\end{aligned}
\end{equation}
Noting $\sigma_{1}=-\frac{u_{m}-u_{-}}{v_{m}-v_{-}}$, we have
\begin{equation}\label{7.73}
	J_{12}+J_{22}=(v_{-}-v_{m})(X_{\infty}(t)-X_{\infty}(t_{0})).
\end{equation}

On the other hand, for $J_{11}$, it follows from \eqref{1.21}, \eqref{A.1}, \eqref{7.67}  that
\begin{equation}\label{7.74}
	\begin{aligned}
    |J_{11}|&\leq \int_{-r-t_{0}+\sigma_{1}T}^{0}\Big\vert\mathfrak{m}_{a}(t,x)-\bar{\mathfrak{v}}_{X_{\infty}}(t,x)\Big\vert\mathbf{1}_{\{\mathfrak{m}_{a}\leq 3v_{-}\}}\,dx\\
    &\quad +\int_{-r-t_{0}+\sigma_{1}T}^{0}\Big\vert\mathfrak{m}_{a}(t,x)-\bar{\mathfrak{v}}_{X_{\infty}}(t,x)\Big\vert\mathbf{1}_{\{\mathfrak{m}_{a}\geq 3v_{-}\}}\,dx+\int_{\R}\psi(x)\mathfrak{v}_{s}(t,dx)\\
    &\leq \frac{1}{\sqrt{\tilde{c}_{1}}}\int_{-r-t_{0}+\sigma_{1}T}^{0}\sqrt{Q\big(\mathfrak{m}_{a}(t,x)|\bar{\mathfrak{v}}_{X_{\infty}}\big)}\,dx\\
    &\quad +\frac{1}{\tilde{c}_{2}}\int_{-r-t_{0}+\sigma_{1}T}^{0}Q\big(\mathfrak{m}_{a}(t,x)|\bar{\mathfrak{v}}_{X_{\infty}}\big)\,dx +\frac{1}{|Q'(v_{-})|}\int_{\R}\psi(x)|Q'(\overline{V})|\mathfrak{v}_{s}(t,dx)\\
    &\leq C(T)\sqrt{\mathcal{E}_{0}}+C\mathcal{E}_{0}.
     \end{aligned}
\end{equation}
Similarly, we have
\begin{equation}\label{7.75}
	\begin{aligned}
		|J_{13}|\leq C(T)\sqrt{\mathcal{E}_{0}}+C\mathcal{E}_{0}.
     \end{aligned}
\end{equation}
For $J_{21}$, using \eqref{7.67}, \eqref{1.21}, one has
\begin{equation}\label{7.76}
	\begin{aligned}
		|J_{21}|&\leq \|\psi'(x)\|_{L^{\infty}}\int_{t_{0}}^{t}\Big(\int_{-r-t_{0}+\sigma_{1}T}^{-r+\sigma_{1}T}+\int_{\frac{\sigma_{1}t_{0}}{2}}^{0}\Big)\Big\vert \mathfrak{u}_{\infty}(s,x)-\bar{\mathfrak{u}}_{X_{\infty}}(s,x)\Big\vert\,dxds\leq \frac{C(T)}{\sqrt{t_{0}}}\sqrt{\mathcal{E}_{0}}.
	\end{aligned}
\end{equation}
Combining \eqref{7.73}-\eqref{7.76} and \eqref{7.70}, we get
\begin{equation*}
|X_{\infty}(t)-X_{\infty}(t_0)|\leq C(T)\mathcal{E}_{0}+C(T)(1+\frac{1}{\sqrt{t_{0}}})\sqrt{\mathcal{E}_{0}}\qquad \text{for a.e. }t\in (t_{0},T),\nonumber
\end{equation*}
which, together with \eqref{7.10-1} , yields that
\begin{equation}\label{7.78}
	\begin{aligned}
	|X_{\infty}(t)|&\leq \hat{C}\cdot (\mathcal{E}_{0}+t_{0})+ C(T)\mathcal{E}_{0}+C(T)(1+\frac{1}{\sqrt{t_{0}}})\sqrt{\mathcal{E}_{0}}\\
	&\leq Ct_{0}+C(T)\mathcal{E}_{0}+C(T)(1+\frac{1}{\sqrt{t_{0}}})\sqrt{\mathcal{E}_{0}}\qquad \text{for a.e. }t\in (0,T),
	\end{aligned}
\end{equation}
where $C$ is independent of $t_{0}$, and $C(T)$ is a constant depending only on $T$.

If $\mathcal{E}_{0}^{\frac{1}{3}}\geq T$, by the density of $t_{0}$, we can take $t_{0}\in (0,T)$ satisfying $T/2<t_{0}<T$ in \eqref{7.78} to obtain
$$
|X_{\infty}(t)|\leq C(T)(\mathcal{E}_{0}+\mathcal{E}_{0}^{\frac{1}{3}}).
$$
If $0<\mathcal{E}_{0}^{\frac{1}{3}}<T$, by the density of $t_{0}$, we can take $t_{0}\in (0,T)$ satisfying $\mathcal{E}_{0}^{\frac{1}{3}}/2<t_{0}<\mathcal{E}_{0}^{\frac{1}{3}}$ in \eqref{7.78} to obtain
$$
|X_{\infty}(t)|\leq C(T)(\mathcal{E}_{0}+\mathcal{E}_{0}^{\frac{1}{3}}).
$$
Finally, if $\mathcal{E}_{0}=0$, then it follows from \eqref{7.78} that $|X_{\infty}(t)|\leq Ct_{0}$ for almost everywhere $t_{0}\in (0,T)$. Taking $t_{0}\to 0$, one has $X_{\infty}(t)\equiv 0\leq C(T)(\mathcal{E}_{0}+\mathcal{E}_{0}^{\frac{1}{3}})$. Therefore the proof of \eqref{1.24} is complete.

\subsection{Proof of \eqref{1.21-1}.} By direct calculations, one has from \eqref{1.21} that
\begin{align}\label{7.79}
&\int_{\R} \f12 |\mathfrak{u}_{\infty}-\bar{\mathfrak{u}}|^2 dx+\int_{\R}Q(\mathfrak{m}_{a}\,|\,\bar{\mathfrak{v}})\,dx+\int_{\R}|Q'(\overline{V}(t,x))|\mathfrak{v}_{s}(t,dx)\nonumber\\
&\leq \frac{1}{2}\int_{\R}|\mathfrak{u}_{\infty}(t,x)-\bar{\mathfrak{u}}_{X_{\infty}}(t,x)|^2\,dx+\frac{1}{2}\int_{\R}|\bar{\mathfrak{u}}_{X_{\infty}}(t,x)-\bar{\mathfrak{u}}(t,x)|^2\,dx\nonumber\\
&\quad +\int_{\R}Q(\mathfrak{m}_{a}\,|\,\bar{\mathfrak{v}}_{X_{\infty}})\,dx+\int_{\R}\big[Q(\mathfrak{m}_{a}\,|\,\bar{\mathfrak{v}})-Q(\mathfrak{m}_{a}\,|\,\bar{\mathfrak{v}}_{X_{\infty}})\big]\,dx+\int_{\R}|Q'(\bar{V}(t,x)|\mathfrak{v}_{s}(t,dx)\nonumber\\
&\leq C\mathcal{E}_{0}+\frac{1}{2}\int_{\R}|\bar{\mathfrak{u}}_{X_{\infty}}(t,x)-\bar{\mathfrak{u}}(t,x)|^2\,dx+\int_{\R}\big\vert Q(\mathfrak{m}_{a}\,|\,\bar{\mathfrak{v}})-Q(\mathfrak{m}_{a}\,|\,\bar{\mathfrak{v}}_{X_{\infty}})\big\vert\,dx.
\end{align}
For the second term on RHS of \eqref{7.79}, noting \eqref{1.13} and \eqref{1.23}, we obtain from \eqref{1.24} that
\begin{align}\label{7.80}
\frac{1}{2}\int_{\R}|\bar{\mathfrak{u}}_{X_{\infty}}(t,x)-\bar{\mathfrak{u}}(t,x)|^2\,dx\leq \frac{1}{2}|u_{-}-u_{m}|^2\cdot|X_{\infty}(t)|\leq C(T)(\mathcal{E}_{0}+\mathcal{E}_{0}^{\frac{1}{3}}).
\end{align}
For the last term on RHS of \eqref{7.79}, it follows from \eqref{1.21}, \eqref{1.24} and \eqref{A.1} that
\begin{align}\label{7.81}
&\int_{\R}\big\vert Q(\mathfrak{m}_{a}\,|\,\bar{\mathfrak{v}})-Q(\mathfrak{m}_{a}\,|\,\bar{\mathfrak{v}}_{X_{\infty}})\big\vert\,dx\nonumber\\
&\leq \frac{1}{\gamma-1}\int_{\R}\big\vert (\bar{\mathfrak{v}}_{X_{\infty}})^{-\gamma+1}-(\bar{\mathfrak{v}})^{-\gamma+1}\big\vert\,dx+\int_{\R}\big\vert(\bar{\mathfrak{v}}_{X_{\infty}})^{-\gamma}-(\bar{\mathfrak{v}})^{-\gamma}\big\vert\cdot \big\vert \mathfrak{m}_{a}-\bar{\mathfrak{v}}_{X_{\infty}}\big\vert\,dx+\int_{\R}(\bar{\mathfrak{v}})^{-\gamma}\big\vert \bar{\mathfrak{v}}_{X_{\infty}}-\bar{\mathfrak{v}}\big\vert\,dx\nonumber\\
&\leq C(|(v_{m}^{-\gamma+1}-(v_{-})^{-\gamma+1}|+|v_{m}-v_{-}|) |X_{\infty}(t)|\nonumber\\
&\quad +\int_{\R}\big\vert(\bar{\mathfrak{v}}_{X_{\infty}})^{-\gamma}-(\bar{\mathfrak{v}})^{-\gamma}\big\vert\cdot \big\vert \mathfrak{m}_{a}-\bar{\mathfrak{v}}_{X_{\infty}}\big\vert\mathbf{1}_{\{\mathfrak{m}_{a}\geq 3v_{-}\}}\,dx\nonumber\\
&\quad +\int_{\R}\big\vert(\bar{\mathfrak{v}}_{X_{\infty}})^{-\gamma}-(\bar{\mathfrak{v}})^{-\gamma}\big\vert\cdot \big\vert \mathfrak{m}_{a}-\bar{\mathfrak{v}}_{X_{\infty}}\big\vert\mathbf{1}_{\{\mathfrak{m}_{a}\leq 3v_{-}\}}\,dx\nonumber\\
&\leq C(T)(\mathcal{E}_{0}+\mathcal{E}_{0}^{\frac{1}{3}})+C\int_{\R}Q(\mathfrak{m}_{a}\,|\,\bar{\mathfrak{v}}_{X_{\infty}})\,dx+C\int_{\R}\big\vert(\bar{\mathfrak{v}}_{X_{\infty}})^{-\gamma}-(\bar{\mathfrak{v}})^{-\gamma}\big\vert\cdot \sqrt{Q(\mathfrak{m}_{a}\,|\,\bar{\mathfrak{v}}_{X_{\infty}})}\,dx\nonumber\\
&\leq C(T)(\mathcal{E}_{0}+\mathcal{E}_{0}^{\frac{1}{3}})+C\mathcal{E}_{0}+C\sqrt{\mathcal{E}}_{0}\cdot |X_{\infty}(t)|^{\frac{1}{2}}\leq C(T)(\mathcal{E}_{0}+\mathcal{E}_{0}^{\frac{1}{3}}).
\end{align}
Substituting \eqref{7.80}--\eqref{7.81} into \eqref{7.79}, we conclude \eqref{1.21-1}. Therefore the proof of Theorem \ref{thm1.1} is complete. $\hfill\square$


\bigskip
\noindent{\bf Acknowledgments.}
Feimin Huang's research is partially supported by National Key R\&D Program of China No. 2021YFA1000800, and National Natural Sciences Foundation of China No. 12288201. Yi Wang's research is supported by National Key R\&D Program of China, NSFC (Grant No. 12171459, 12288201, 12090014) and CAS Project for Young Scientists in Basic Research, Grant No. YSBR-031. Yong Wang's research is partially supported by National Natural Science Foundation of China No. 12022114, 12288201, and CAS Project for Young Scientists in Basic Research, Grant No. YSBR-031.

\end{document}